%% file: deform-1.tex
\documentclass{amsart}
\include{headerstemplate}

\begin{document}
\title[Deformations and automorphisms: globalization]{Deformations and 
automorphisms: a framework for globalizing local tangent and obstruction 
spaces}
\author{Brian Osserman}
\begin{abstract} Building on Schlessinger's work, we define a framework
for studying geometric deformation problems which allows us to systematize
the relationship between the local and global tangent and obstruction spaces
of a deformation problem. Starting from Schlessinger's functors of Artin
rings, we proceed in two steps: we replace functors to sets by categories
fibered in groupoids, allowing us to keep track of automorphisms, and we 
work with deformation problems naturally associated to a scheme $X$, and 
which naturally localize on $X$, so that we can formalize the local 
behavior. The first step is already carried out by Rim in the context
of his homogeneous groupoids, but we develop the theory substantially
further. In this setting, many statements known for a range of specific 
deformation problems can be proved in full generality, under very general
stack-like hypotheses.
\end{abstract}
\thanks{This paper was partially supported by a fellowship from the
National Science Foundation.}
\maketitle

\tableofcontents

\section{Introduction}

Given a formal deformation problem, two of the most basic questions one can 
ask are: to what extent is it representable, and what are its tangent and 
obstruction spaces? In \cite{sc2}, Schlessinger gave an elementary and 
concise answer to the first question. The second question appears to be 
much more involved, with the most definitive work on the subject being 
Illusie's \cite{il1}. We propose a new framework which allows us to describe 
precisely the relationship between the global tangent and obstruction spaces 
and the local ones. While far less ambitious in scope, our approach is 
relatively elementary and allows us to treat a wide range of deformation 
problems uniformly and transparently. All the statements which we prove 
are well known in examples, but we show that they are in fact formal 
consequences of relatively mild hypotheses in a generality comparable to
that of Schlessinger's work. As an added bonus, we are able to replace 
Schlessinger's conditions with more natural descent-theoretic ones, and we 
ultimately obtain results on the representing scheme.

\subsection{Statements}

A basic example of the sort of intuitive statement which we wish to
be able to state as a general theorem is the following: if one has a 
deformation problem associated to a scheme $X$, which has a sheaf $\cA$ 
on $X$ of infinitesimal automorphisms, and if (as in the case of deformations 
of smooth varieties, or vector bundles on a fixed variety) both deformations 
and obstructions are locally trivial, then the tangent space is 
$H^1(X,\cA)$, and obstructions lie naturally in $H^2(X,\cA)$. The idea is
of course that we obtain our deformations by gluing together copies of the 
trivial deformation on open covers, and such gluings are controlled by 
$\cA$. 

In order to formalize such statements, we need to modify 
Schlessinger's context of functors of Artin rings in two ways. First,
in order to be able to work with infinitesimal automorphisms, we are led 
to replace functors to sets with groupoids. Under additional stack-type 
gluing conditions on the deformation problem, we call such objects 
{\bf deformation stacks}. Second, in Schlessinger's work there is no base 
scheme $X$ in the picture on which one can formulate the concepts of local 
or global. Accordingly, we consider problems associated to a scheme $X$, and 
which localize naturally on $X$: i.e., instead of associating deformations
to every Artin ring, we associate deformations to every pair $(U,A)$ of an 
open set on $X$ and an Artin ring. Under further gluing conditions, these
objects will be called {\bf geometric deformation stacks}, or 
{\bf gd-stacks}. These will naturally carry two sheaves of $k$-vector spaces 
on $X$: $\cA$, the sheaf of infinitesimal automorphisms, and $\cT$, the
sheaf of local first-order deformations.

We prove the following two theorems (see \S \ref{review} below for a 
review of Schlessinger's terminology, and \S \ref{sec:gd-stacks} for new 
definitions):

\begin{thm}\label{main-desc} Let $\cS$ be a gd-stack. 
Then the tangent space $T_{\cS}$ of $\cS$ fits into an exact sequence of 
$k$-vector spaces
$$0 \to H^1(X,\cA) \to T_{\cS} \to H^0(X,\cT) \to H^2(X,\cA),$$
and if we are given a local obstruction sheaf $\cOb$ for $\cS$, 
we have successive obstructions lying in 
$H^0(X,\cOb)$, $H^1(X,\cT)$, and $H^2(X,\cA)/H^0(X,\cT)$.
\end{thm}

Note that as a special case of the theorem, we obtain a precise version of 
our earlier assertion for the case that deformations and obstructions are 
locally trivial, when we can set $\cT=\cOb=0$; see Corollary 
\ref{cor:h1h2}.

As a consequence of Theorem \ref{main-desc}, we draw conclusions about the 
existence and properties of a hull, under fairly general circumstances.

\begin{thm}\label{main-rep} Let $\cS$ be a gd-stack on a scheme $X$. Then:

\begin{ilist}
\itm the associated functor $F_{\cS_{(X,\us)}}$ satisfies Schlessinger's 
(H1) and (H2), and satisfies (H4) if and only if for each tiny extension
$A' \to A$ in $\Art(\Lambda,k)$, and each object $\eta \in A'$, the natural
map
$$\Aut(\eta) \to \Aut(\eta|_A)$$
is surjective;
\itm if $X$ is proper, and the sheaves $\cA$ and $\cT$ both
carry the structure of coherent $\cO_X$-modules, then $F_{\cS_{(X,\us)}}$
satisfies Schlessinger's (H3), so has a hull $R$;
\itm if further we have a local obstruction sheaf $\cOb$ for $\cS$,
and it carries the structure of a coherent $\cO_X$-module, then
$$h^0(X,\cT)+h^1(X,\cA)- h^0(X,\cOb) - h^1(X,\cT)- h^2(X,\cA) 
\leq \dim R -\dim \Lambda \leq \dim T_{\cS},$$ 
and if the first inequality is an equality and $\Lambda$ is regular, $R$ 
is a local complete intersection ring. If we have
$$h^0(X,\cOb) = h^1(X,\cT) = 
\dim T_{\cS} + h^2(X,\cA)-h^1(X,\cA)-h^0(X,\cT) = 0,$$ 
then $R$ is smooth over
$\Lambda$.
\end{ilist}
\end{thm}

Although our conditions for a gd-stack are formally stronger than
Schlessinger's, it appears to be the case that (at least when 
the problem can be naturally defined for some $X$ and its open subsets)
any deformation problem which satisfies (H1) and (H2) does so
because it is associated to a gd-stack, and that moreover one will verify 
the conditions for a gd-stack in the process of checking (H1) and (H2).
This is born out by a number of examples, treated in \S 
\ref{subsec:first-exs} below.

We also mention that in the locally unobstructed case, the form of the
tangent and obstruction spaces are such that one expects that they arise
as the hypercohomology of a two-term complex. We explore this further in 
\cite{os15}.

Although this work is self-contained and was developed largely 
independently, it is closely related to work of Grothendieck and Rim 
as follows: Grothendieck's work on Exal in \cite{gr4} laid the framework 
for treating deformations from a groupoid point of view, and exploiting 
certain additive structures on categories in order to do so. Meanwhile, 
Rim's theory of homogeneous groupoids \cite{ri1}, which is formally 
equivalent to our deformation stacks, began the process of treating 
deformation problems systematically in the context of groupoids. The bulk 
of the present paper is \S \ref{sec:def-stacks}, which develops the theory 
of deformation stacks further, and may be viewed as synthesizing and 
expanding on the ideas of Grothendieck and Rim, using categorical torsor 
structures to prove very general statements on the structure of liftings 
of objects and automorphisms over small ring extensions. Finally, it should 
be emphasized that everything we do has been very well known in specific 
examples for some time; our main contribution is to provide a framework in 
which the arguments can be systematized, so that they apply formally to 
many deformations problems at once. Schlessinger's original paper was a 
valuable source of such statements made in more specific contexts, as were 
the lecture notes \cite{ha2} of Robin Hartshorne.

\subsection*{A word on the use of stacks} 

We wish to underline that although we use concepts from the theory of stacks,
no reader should be intimidated by this fact; indeed, we make no use of
algebraic stacks or the \'etale topology, and one should think of stacks in 
this context as being nothing more than functors which remember 
automorphisms and satisfy certain natural gluing conditions (which we will
restate below). In fact, we use no results at all from the theory of stacks, 
and our presentation is entirely elementary and self-contained.

\subsection*{Acknowledgements} 

I would like to thank David Eisenbud, Ravi Vakil, Charles Weibel, Robin 
Hartshorne, and Luc Illusie, and particularly Martin Olsson and Max 
Lieblich, for many patient and illuminating conversations.

\subsection{Review of Schlessinger}\label{review}

For the convenience of the reader and in order to assign terminology 
intended to motivate the relationship to deformation stacks, we briefly 
review Schlessinger's criteria. Given a
field $k$, and a complete Noetherian ring $\Lambda$ with residue field $k$,
we denote by $\Art(\Lambda,k)$ the category of local Artin $\Lambda$-algebras
with residue field $k$. We will use the notational convention that $\epsilon$
is always a square-zero element. For any $k$-vector space $V$, we will also 
denote by $k[V]$ the algebra having additive group $k \oplus V$,
with square-zero multiplication for elements of $V$. Both $k[\epsilon]$ and
$k[V]$ always denote rings endowed with the ``trivial'' $\Lambda$-algebra
structure (i.e., the one factoring through $\Lambda \twoheadrightarrow k$).
Finally, for any $A \in \Art(\Lambda,k)$ we denote by $\pi:A \to k$ the
residue field map.

$\Lambda$ is frequently either $k$ or, if $k$ is perfect of characteristic 
$p$, and one wants to work in mixed characteristic, the Witt vectors 
$W(k)$. The latter case is universal: every complete Noetherian local ring 
with residue field $k$ is canonically an algebra over $W(k)$ (see 
Proposition 10 of II, \S 5 of \cite{se1}). When working with families over 
a base space, one often takes $\Lambda$ to be a complete local ring of the 
base.

Schlessinger considered functors of the following type.

\begin{defn} A covariant functor $F:\Art(\Lambda,k) \to \Set$ is a 
{\bf predeformation
functor} if $F(k)$ consists of a single element. If $F$ is a predeformation
functor, we say that $T_F:=F(k[\epsilon])$ is its {\bf tangent space}.
\end{defn}

\begin{notn}\label{n:triv-objs}
Given a predeformation functor $F$, we use $\zeta_0$ to denote the unique 
object of $F(k)$, and $\zeta_V$ (respectively, $\zeta_{\epsilon}$) to
denote the object induced on $k[V]$ (respectively, $k[\epsilon]$) by
$\zeta_0$ under the structure map.

Given a morphism $f:A \to A'$ in $\Art(\Lambda,k)$, and an object 
$\eta \in F(A)$, we will denote the object of $F(A')$ induced by $\eta$
under $f$ by $f_*(\eta)$ or, when there is no ambiguity, by $\eta|_{A'}$.
\end{notn}

We will deviate slightly from Schlessinger's terminology.

\begin{defn} A surjection $p: A' \to A$ in $\Art(\Lambda,k)$ with kernel
$I$ is {\bf small} if $I \cdot \fm_{A'}=0$.
We say that $p$ is {\bf tiny} if it is small,
and if further $I$ is a principal ideal.
\end{defn}

Note that the kernel of a small extension is necessarily square-zero and a 
$k$-vector space, and that a tiny extension has kernel isomorphic to $k$.
Schlessinger used ``small'' for what we call ``tiny''; our usage of
``small'' follows Huybrechts and Lehn (2.A.5 of \cite{h-l}), as well as 
Fantechi and G\"ottsche (Definition 6.1.9 of \cite{f-g1}). The more 
general notion of small is important for certain applications of 
obstruction theory; see in particular Theorem \ref{dim-ests} below.

Given morphisms 
$A' \to A$, $A'' \to A$ in $\Art(\Lambda,k)$, there is a natural map 
\begin{equation}\label{func-nat-map} \cF(A' \times _A A'') \to \cF(A') \times
_{\cF(A)} \cF(A'').\end{equation}

Schlessinger's conditions are the following.

\begin{itemize}
\item[(H1)] The map (\ref{func-nat-map}) is surjective whenever $A'' \to A$
is a tiny extension.
\item[(H2)] The map (\ref{func-nat-map}) is a bijection when $A=k$, and
$A''=k[\epsilon]$.
\item[(H3)] The tangent space $T_F$ is finite-dimensional over $k$.
\item[(H4)] The map (\ref{func-nat-map}) is a bijection whenever
$A''=A'$ and $A' \to A$ is a tiny extension.
\end{itemize}

We remark that if (H2) is satisfied, then $T_F$ can be naturally given the
structure of a $k$-vector space. We further remark that every reasonable
deformation problem seems to satisfy (H1) and (H2), while (H3) is typically
satisfied where proper schemes are involved, and (H4) is a substantially 
stronger condition, closely tied to the behavior of automorphisms. We 
therefore make the following definition.

\begin{defn} We say that a predeformation functor $F$ is a {\bf deformation
functor} if it satisfies Schlessinger's (H1) and (H2).
\end{defn}

Denote by $\widehat{\Art}(\Lambda,k)$ the category of complete Noetherian 
local $\Lambda$-algebras with residue field $k$.
We recall that a predeformation functor $F$ may be extended to a functor
$\hat{F}: \widehat{\Art}(\Lambda,k) \to \Set$, simply by taking the appropriate 
limits over all $R/\fm_R^n$ for $R \in \widehat{\Art}(\Lambda,k)$. We recall 
the following basic definitions, the first one a direct extension of the 
notion of formal smoothness to morphisms of functors:

\begin{defn} If $F,F':\Art(\Lambda,k) \to \Set$ are functors, and we have
a morphism $\vp: F \to F'$, we say that $\vp$ is {\bf formally smooth} if 
for every
surjective map $A' \twoheadrightarrow A$ in $\Art(\Lambda,k)$, the canonical
map 
$$F(A') \to F(A) \times_{F'(A)} F'(A')$$
is surjective.
\end{defn}

The notions of representability are the following.

\begin{defn} We say that $F:\Art(\Lambda,k)\to \Set$ is 
{\bf prorepresentable} if
there is a pair $(R,\xi)$, with $R \in \widehat{\Art}(\Lambda,k)$, and
$\xi \in \hat{F}(R)$, such that the induced map 
$h_R|_{\Art(\Lambda,k)} \to F$ is an isomorphism of functors.

We say that a pair $(R,\xi)$ is a {\bf hull} of $F$ if the map
$h_R|_{\Art(\Lambda,k)} \to F$ is formally smooth, and induces an 
isomorphism $h_R(k[\epsilon]) \risom F(k[\epsilon])$.
\end{defn}

Schlessinger's basic theorem is as follows.

\begin{thm}\label{thm:schl} (Schlessinger, Theorem 2.11 of \cite{sc2}) Let 
$F$ be a predeformation functor. Then $F$ has a hull if and only if $F$ 
satisfies (H1), (H2), and (H3). Moreover, $F$ is prorepresentable if and 
only if it satisfies in addition (H4).
\end{thm}

\section{Deformation stacks}\label{sec:def-stacks}

We begin with a simple translation of Schlessinger's work on deformation
functors into a more stack-theoretic language, which, while imposing 
potentially stricter hypotheses, gives a more complete picture of the 
situation. This is equivalent to Rim's homogeneous groupoids, but our 
definition is more motivated by the descent conditions of stack theory,
and we examine this relationship more closely than is required to prove
our main results. We then proceed to give a number of examples of
deformation problems which naturally form deformation stacks, and 
finally to prove several technical results which ultimately play an 
important role in proving our main theorems. The main theorem in this
direction, which is vital to understanding obstructions, is that with no 
additional hypotheses, liftings of objects and automorphisms over small 
extensions can be studied in terms of the tangent space $T_{\cS}$ and
the infinitesimal automorphism group $A_{\cS}$ respectively.

\subsection{Definitions}\label{ds-defs}

To avoid technical issues, we assume throughout that we work with
small categories, so that every category has an 
associated set of isomorphism classes. This can be accomplished by,
for instance, working in a fixed universe in the sense of Grothendieck.
We begin by recalling the basic definitions relating to categories fibered 
in groupoids.

\begin{defn}\label{def:groupoid} Fix a category $C$. We say that a category 
$\cS$, together with
a (covariant) functor to $C$, is a {\bf category fibered in groupoids}
over $C$, if:
\begin{ilist}
\itm Given a morphism $T \to T'$ in $C$, and an object 
$\eta'$ in $\cS$ over $T'$, there exists $\eta$ in $\cS$ over $T$ and a 
morphism $\eta \to \eta'$ over the given morphism $T \to T'$.
\itm Given a diagram with the solid arrows below:
$$\xymatrix{{\eta}\ar@{..>}[dr]\ar[drr]& & \\
& {\eta'}\ar[r] & {\eta''} \\
{T} \ar[r] & {T'} \ar[r] & {T''}}$$
the dotted arrow, making the triangle commute, exists and is unique.
\end{ilist}
\end{defn}

The terminology is justified as follows: given a functor $\cS \to C$, 
and an object $T \in C$, we define the {\bf fiber category} $\cS_T$
to consist of the objects of $\cS$ lying over $T$, together with the
morphisms lying over the identity morphism $T \to T$. It then follows
from (ii) above that every fiber category $\cS_T$ of a category fibered in
groupoid $\cS$ is indeed a groupoid, in the sense that every morphism is
an isomorphism.

An observation which we will use implicitly many times when checking
compatibility of definitions is the following, which is an immediate
consequence of the definitions.

\begin{prop} Let $\cS \to C$ be a category fibered in groupoids. Then
every morphism $f:\eta' \to \eta''$ in $\cS$, lying over $T' \to T''$
in $C$, is a ``relative monomorphism,'' in the sense that if 
$g_1,g_2:\eta \to \eta'$ both lie over a given $T \to T'$ in $C$, and
$f \circ g_1 = f \circ g_2$, then $g_1=g_2$.
\end{prop}

The definition of a category fibered in groupoids implicitly incorporates 
notions of pullbacks of objects and morphisms. We begin with objects.

\begin{defn} Suppose we are given a morphism $f:T \to T'$ in $C$, and
$\eta' \in \cS_{T'}$.
The pair of $\eta$ and the morphism $\eta \to \eta'$ given by 
(i) above is the {\bf pullback} of $\eta'$ under $f$,
which we denote by $f^* \eta'$, or $\eta'|_{T}$. 
\end{defn}

By (ii) above we immediately see the following uniqueness statement.

\begin{prop} A pullback $f^*(\eta')$ is unique up to unique isomorphism.
\end{prop}

We emphasize that although we use object notation for pullbacks of 
objects, in order to obtain the above uniqueness it is necessary to 
include also the data of the morphism $\eta \to \eta'$. Even with
this data, pullbacks are not uniquely defined in an absolute sense, but
because they are unique up to unique isomorphism, in practice they may 
often be treated as uniquely defined. For instance, we may canonically 
identify morphisms to or from any two pullbacks of a given object. 
Sometimes we are required to explicitly consider the ``unique isomorphism'' 
between two pullbacks, as in the discussion of the cocycle condition below.

We can also define pullbacks of morphisms.

\begin{defn} Suppose we are given a morphism $f:T \to T'$ in $C$.
If we have two objects $\eta,\eta' \in S_{T'}$, and a morphism
$\vp:\eta \to \eta'$ over the identity morphism of $T'$, we also have the
{\bf pullback} $f^*(\vp):f^*(\eta) \to f^*(\eta')$ of $\vp$ under $f$,
a morphism lying over the identity morphism of $T$. This pullback, also
sometimes denoted by $\vp|_{T'}$, is defined as the morphism obtained by 
applying (ii) above to the composed morphism 
$f^*(\eta) \to \eta \overset{\vp}{\to} \eta'$ and $f^*(\eta') \to \eta'$.
\end{defn}

The following basic compatibility properties follow from the definitions.
More specifically, the first is immediate from (i), while the rest may be 
checked using the uniqueness hypothesis in (ii) of the definition of a 
category fibered in groupoids. 

\begin{prop}\label{prop:cfgs-basic} We have:
\begin{ilist}
\itm The pullback of a pullback is a pullback: given $T \to T' \to T''$
in $C$, and $\eta \in S_{T''}$, then $(\eta|_{T'})|_T$ is a pullback
of $\eta$ under the composition $T \to T''$.
\itm The isomorphisms between different pullbacks satisfy the cocycle 
condition: given $T \to T'$ and $\eta \in S_{T'}$, suppose that 
$\eta_1,\eta_2,\eta_3$ are three different pullbacks of $\eta$ to $T$,
and for each $i,j$ let
$\vp_{i,j}:\eta_i \risom \eta_j$ be the unique isomorphism identifying
$\eta_i$ with $\eta_j$. Then these satisfy the
{\bf cocycle condition} $\vp_{2,3} \circ \vp_{1,2} = \vp_{1,3}$.
\itm Pullback of morphisms commutes with composition: 
given $T \overset{f}{\to} T' \overset{g}{\to} T''$, and $\eta,\eta'$ in 
$\cS_{T''}$, and a morphism $\vp:\eta \to \eta'$ lying over the identity 
morphism of $T''$, we have $f^*(g^*(\vp))=(f \circ g)^*(\vp)$.
\itm Pullback of morphisms commutes with composition: given
$f:T \to T'$, and $\eta,\eta',\eta''$ in $\cS_{T'}$, and two morphisms 
$\eta \overset{\vp}{\to} \eta' \overset{\vp'}{\to} \eta''$ lying over the 
identity morphism of $T'$, we have $f^*(\vp' \circ \vp)=f^*(\vp') \circ
f^*(\vp)$.
\end{ilist}
\end{prop}

We now give a definition of deformation stack which is very close to
(and indeed, equivalent to; see Remark \ref{rem:rim} below) Rim's 
definition of homogeneous groupoid. For a definition which places
deformation stacks more visibly in the context of standard stack conditions,
see \S \ref{sec:groth-tops} below.

\begin{defn}\label{def:def-stack} A category $\cS$ fibered in groupoids over 
$\Art(\Lambda, k)^{\opp}$ is a {\bf deformation stack} if the fiber of 
$\cS$ over $k$ is the ``trivial'' groupoid: i.e., there exists a unique 
morphism between any two objects; and if for every square of the form
$$\xymatrix{{A} & {A'} \ar[l]^{p'} \\
{A''} \ar[u]^{p''} & {A' \times_A A''} \ar[u]  \ar[l]}$$
in $\Art(\Lambda,k)$ with $p''$ surjective, we have:
\begin{ilist}
\itm (``morphisms form a sheaf'') Given objects $\eta,\eta'$ in $\cS$
over $A' \times_A A''$, the natural map
$$\Mor(\eta,\eta') \to 
\Mor(\eta|_{A'},\eta'|_{A'}) \times_{\Mor(\eta|_A,\eta'|_A)} 
\Mor(\eta|_{A''},\eta'|_{A''})$$
is a bijection. Here all sets of morphisms are taken in the appropriate
fiber categories, or equivalently, are assumed to lie over the identity 
morphism of the underlying Artin ring.
\itm (``objects satisfy effective descent'') Given objects $\eta',\eta''$ in
$\cS$ over $A'$ and $A''$ respectively, and a morphism 
$\vp:\eta'|_A  \to \eta''|_A$ over $\id_A$, there exists an object 
$\tilde{\eta}$ over $A' \times_A A''$ such that $\eta'=\tilde{\eta}|_{A'}$,
$\eta''=\tilde{\eta}|_{A''}$, and $\vp$ is the unique isomorphism 
identifying both $\eta'|_A$ and $\eta''|_A$ as pullbacks of $\tilde{\eta}$
to $A$. 
\end{ilist}
\end{defn}

\begin{rem}\label{rem:rim} The definition of a deformation stack can be 
stated more
compactly by requiring that the groupoid version of the natural map 
\eqref{func-nat-map} be an equivalence of categories. This is equivalent
to our stack conditions because the conditions on morphisms and objects
are equivalent to the natural functor being fully faithful and essentially
surjective, respectively. Thus, the definition we have given of deformation 
stack is equivalent to Rim's homogeneous groupoids (see Definition 2.5 and 
Remark 2.6 (b) of \cite{ri1}). 
\end{rem}

To avoid unnecessary obfuscation
involving the use of the opposite category, when working with a deformation
stack $\cS$, we will work freely with ring morphisms and pushforward of
objects of $\cS$. Note the slightly confusing situation that the direction
of a pushforward morphism in $\cS$ will be opposite to the direction of 
the ring homomorphism in $\Art(\Lambda,k)$ over which it lies.

\begin{notn} We still use the notation $\zeta_0$, $\zeta_V$, and 
$\zeta_{\epsilon}$ in the context of deformation stacks as in Notation
\ref{n:triv-objs}; we assume we have fixed a choice of object $\zeta_0$,
and of pushforwards $\zeta_V$ and $\zeta_{\epsilon}$.
\end{notn}

Although these are not uniquely-defined objects, they are 
unique up to unique isomorphism. Note that more generally, the
definition of a category fibered in groupoids together with the
hypothesis that $\cS_k$ is trivial implies that every object $\eta$
in $\cS$ has a unique morphism $\zeta_0 \to \eta$.

To each deformation stack $\cS$, we have the associated functor of 
isomorphism classes $F_{\cS}$. By virtue of the hypothesis that
$\cS_k$ is trivial, $F_{\cS}$ is a predeformation functor, and it is easy 
to check that conditions (i) and (ii) above imply:

\begin{prop}\label{prop:def-functor} Suppose $F_{\cS}$ is the functor 
associated to a deformation stack. Then $F_{\cS}$ is a deformation functor.
\end{prop}

Our philosophy is that although being associated to a deformation stack is
in a literal sense stronger than being a deformation functor, we expect that
``in nature'' any deformation functor is in fact the functor associated to 
a deformation stack, and that furthermore in any given case, the proof that 
it is a deformation functor will include a proof that the natural groupoid 
is a deformation stack. As evidence for this philosophy, we have the
following simple rephrasing of Lemma 1.4.4 of Olsson \cite{ol1}:

\begin{prop} Suppose $\cS$ is a deformation problem obtained via 
restriction around a point of an algebraic stack. Then $\cS$ is a
deformation stack.
\end{prop}

However, deformation stacks are far more general, arising for instance from 
deformations of any scheme (without any polarization, and without 
properness hypotheses), as well as in a variety of other contexts, 
discussed in \S \ref{subsec:first-exs} below. 

Moving beyond Schlessinger's conditions (H1) and (H2), another advantage of 
the stack perspective is a sharp understanding of (H4),
already observed (in necessarily imprecise form) by Schlessinger in 
Remark 2.15 of \cite{sc2}, and made more precise by Rim in Proposition
2.7 \cite{ri1}.

\begin{prop}\label{h4} For every surjection $A \twoheadrightarrow B$ in
$\Art(\Lambda,k)$, we have that the natural map
\begin{equation}\label{eq:h4-map} 
F_{\cS}(A \times _B A) \to F_{\cS}(A) \times_{F_{\cS}(B)} F_{\cS}(A)
\end{equation}
is bijective if and only if for every object $\eta \in \cS_A$, the natural 
map from 
$\Aut(\eta)$ to $\Aut(\eta|_B)$ is surjective.

In particular, $F_{\cS}$ satisfies Schlessinger's (H4) if and only if every 
automorphism can be extended over any tiny extension.
\end{prop}

\begin{proof} Fix $A \twoheadrightarrow B$ a surjection in 
$\Art(\Lambda,k)$, 
and an object $\eta \in \cS_{A}$. The stack conditions on $\cS$ given in
Definition \ref{def:def-stack} mean that given objects 
$\eta,\eta' \in \cS_{A}$
such that $\eta|_B \cong \eta'|_B$, fixing an isomorphism gives a 
bijection between isomorphism classes of objects
$\tilde{\eta} \in \cS_{A \times _B A}$ such that $\tilde{\eta}$ restricts
to $\eta$ and $\eta'$ under the two projection maps, and 
$\Aut(\eta|_B)/\left<\Aut(\eta)|_B,\Aut(\eta')|_B\right>$.

In particular, if $\Aut(\eta) \twoheadrightarrow \Aut(\eta|_B)$ and
$\Aut(\eta') \twoheadrightarrow \Aut(\eta'|_B)$, we have that $\tilde{\eta}$
is uniquely determined by $\eta$ and $\eta'$, so we have the asserted 
bijectivity for the natural map \eqref{eq:h4-map}. Conversely, if we 
consider the case $\eta=\eta'$, bijectivity of \eqref{eq:h4-map} 
is equivalent to uniqueness of $\tilde{\eta}$, which implies that 
$\Aut(\eta)$ surjects onto $\Aut(\eta|_B)$.
\end{proof}

Thus, Schlessinger's conditions (H1), (H2), and (H4) all have very natural
interpretations in the context of deformation stacks. We will ultimately
prove the statement asserted in Theorem \ref{main-rep} that there
is a simple sufficient condition on a (geometric) deformation stack to
obtain (H3) as well.

We conclude with definitions of automorphism, tangent and obstruction 
spaces for deformation stacks. While the last two may be defined in terms
of the associated functors, of course the automorphism space is not. In
fact, we work with successive obstruction spaces: the intuition for 
successive obstructions taking value in vector spaces $V_1,\dots,V_n$ is 
that given a 
small extension $A' \to A$ with kernel $I$, and an object $\eta$ over 
$A$, we obtain an element of $V_1 \otimes_k I$ giving a ``first 
obstruction'' to lifting $\eta$ to $A'$, and if that element vanishes, we 
have an element of $V_2 \otimes_k I$ giving a ``second obstruction,'' and 
so forth, and there exists a lift of $\eta$ to $A'$ if and only if every 
obstruction vanishes. However, for the sake of simplicity we make a 
definition focusing on the first non-zero obstruction.

\begin{defn}\label{d:tangent-obs} The {\bf infinitesimal automorphism 
group} $A_{\cS}$ of $\cS$ is the group of automorphisms of 
$\zeta_{\epsilon}$ in $\cS_{k[\epsilon]}$.

The {\bf tangent space} $T_{\cS}$ of $\cS$ is defined to
be $T_{F_{\cS}}$, the tangent space of the associated deformation functor.

Given $k$-vector spaces $V_1,\dots,V_n$, a {\bf successive obstruction
theory} for $\cS$ taking values in $V_1,\dots,V_n$ consists of the
data, for each small extension $A' \to A$ in $\Art(\Lambda,k)$ with
kernel $I$, and $\eta \in \cS_A$, of an $m \in \{1,\dots,n\}$ and
an element $\ob_{\eta,A'} \in V_m \otimes_k I$ such that:
\begin{ilist}
\itm $\ob_{\eta,A'} \neq 0$ unless $m=n$;
\itm there exists $\eta' \in \cS_{A'}$ such that $\eta'|_A=\eta$ if and only
if $m=n$ and $\ob_{\eta,A'}=0$.
\end{ilist}
Furthermore, we impose the following functoriality condition: suppose
we have another small extension $B' \twoheadrightarrow B$ with kernel $J$, 
and $\vp: A' \to B'$ such that $\vp(I) \subseteq J$, so
that $\vp$ also induces maps $I \to J$ and $A \to B$. 
Then we require that
$\ob_{\eta|_B,B'}=\vp(\ob_{\eta,A'})$ if the latter is non-zero or if
$m=n$, and otherwise 
$\ob_{\eta|_B,B'} \in V_{m'} \otimes _k J$ for $m'>m$. 

In the special case that $n=1$, we say that $V_1$ is an {\bf obstruction
space} for $\cS$.
\end{defn}

While $A_{\cS}$ is {\it a priori} only a group, we will see in Corollary
\ref{cor:def-aut-vs} below that it is actually also a $k$-vector space, 
with addition agreeing with composition.

\subsection{Deformations stacks as stacks}\label{sec:groth-tops}

For those familiar with stack theory, particularly as it is typically used
in moduli space theory, our definition of deformation stack has little to 
do with any usual notion of stack. We now clarify the situation by giving 
an equivalent definition of deformation stack which is stated as a standard
stack condition, but applied to a collection of covers which do not 
satisfy the hypotheses of a Grothendieck topology. This discussion may be
considered purely philosophical, and will not be used in anything which
follows.

The first step is to rephrase the squares of Definition \ref{def:def-stack} 
in a more visibly descent-theoretic form.

\begin{defn} We say a square 
$$\xymatrix{{A' \otimes_A A''} & {A'} \ar[l] \\
{A''} \ar[u] & {A} \ar[u]^{p'}  \ar[l]^{p''}}$$
in $\Art(\Lambda,k)$ is a {\bf Schlessinger square} if the following 
conditions are satisfied:
\begin{ilist}
\itm $A \to A' \times A''$ is injective;
\itm $p'$ is surjective; 
\itm $p''(\ker p')$ is closed under multiplication in $A''$.
\end{ilist}
\end{defn}

From a scheme theory point of view, the first condition says that 
$\Spec A' \coprod \Spec A''$ maps scheme-theoretically surjectively to
$\Spec A$, while the second of course means that $\Spec A' \to \Spec A$
is a closed immersion. However, the third condition remains somewhat
mysterious, from a geometric point of view. Nonetheless, given a 
Schlessinger square, it makes geometric sense to think of $A'$ and $A''$
as covering $A$, with $A' \otimes_A A''$ as the intersection. We will
sometimes describe a square of affine schemes as a Schlessinger square
if its opposite square is a Schlessinger square.

The point of this definition is that the squares we considered earlier are 
completely equivalent to Schlessinger squares:

\begin{lem}\label{equiv-diags} Suppose we have $p':A \to A'$ and 
$p'':A \to A''$ defining a Schlessinger square.
Then we have that 
$$A \cong A' \times _{A' \otimes_A A''} A'',$$
and $A'' \to A' \otimes_A A''$ is surjective.

Conversely, given $q':B' \to B$, $q'':B'' \to B$, with $q''$ surjective,
then we have that 
$$B\cong B' \otimes_{B' \times_B B''} B'',$$
and furthermore the maps $B' \times_B B'' \to B'$ and 
$B' \times_B B'' \to B''$ define a Schlessinger square.
\end{lem}

\begin{proof} For both assertions, the key point is that given 
ring homomorphisms $p':C \to C'$ and $p'':C \to C''$, with $p'$
surjective, we can easily check from the universal property of the
tensor product that the natural map $C'' \to C' \otimes _C C''$, which 
factors through $C''/(p''(\ker p'))$, in fact induces an isomorphism
$C''/(p''(\ker p')) \risom C' \otimes _C C''$.
\end{proof}

This lemma should not be surprising, as cofibered product of schemes 
corresponds to union in many Grothendieck
topologies. The second part of the lemma says that the surjectivity 
condition imposed by Schlessinger implies that his $A$, which is in
principal arbitrary, is in fact the (co)intersection of $A'$ and $A''$ over
$A' \times_A A''$.

We therefore have a definition of deformation stack which looks much closer
to imposing a stack condition with respect to a collection of covers. 
Although the covers considered look quite different from the topologies
normally considered in moduli theory, they are similar to those considered 
by Voevodsky for the ``proper cdh-structure'' in \cite{vo3}, so it is 
natural to ask whether our definition of deformation stack is equivalent 
to a stack condition for some Grothendieck topology on $\Art(\Lambda,k)$. 
We first review what it means to have a stack condition with respect to
a given cover on a category fibered in groupoids.

Given a category $\cS$ fibered in groupoids over a category 
$C$, and an arbitrary family of morphisms $\{U_i \to T\}$ in $C$, 
we have restriction maps on morphisms and objects of $\cS$. 

Indeed, given $\eta,\eta' \in \cS_T$, and a morphism $\vp:\eta \to \eta'$
over $\id_T$, we obtain via pullbacks morphisms $\vp_i:\eta|_{U_i} \to
\eta'|_{U_i}$, and (iii) of Proposition \ref{prop:cfgs-basic} implies
that if we pull back $\vp_i$ and $\vp_j$ to $U_{i,j}:=U_i \times_T U_j$,
we have equality. That is, we have a map:
\begin{equation}\label{mor-map} \Mor_{\id_T}(\eta,\eta') \to
\left\{\{\vp_i\in \Mor_{\id U_i}(\eta|_{U_i},\eta'|_{U_i})\}_i: 
\vp_i|_{U_{i,j}}=\vp_j|_{U_{i,j}}
\ \forall i,j\right\},\end{equation}

Similarly, given an object $\eta \in \cS_T$, we obtain pullback objects
$\eta|_{U_i}$, and (i) and (ii) of Proposition \ref{prop:cfgs-basic} imply
that for any $i,j$ both $(\eta|_{U_i})|_{U_{i,j}}$ and 
$(\eta|_{U_j})|_{U_{i,j}}$ are pullbacks $\eta|_{U_{i,j}}$, and hence
related by a unique isomorphism $\vp_{i,j}$, which must then satisfy the
cocycle condition for any $i,j,\ell$ after pullback to 
$U_{i,j,\ell}:=U_i \times_T U_j \times_T U_\ell$. If, as is often done, we 
were to fix choices of pullbacks for every object, we would obtain a map:
\begin{equation}\label{obj-map} \Obj_T \to 
\left\{\begin{aligned}(\{\eta_i \in \Obj_{U_i}\}_i,
\{\vp_{i,j} \in 
\Mor_{\id_T}(\eta_i|_{U_{i,j}},\eta_j|_{U_{i,j}})\}_{i, j}): \\
\vp_{j,\ell}|_{U_{i,j,\ell}} \circ \vp_{i,j}|_{U_{i,j,\ell}} = 
\vp_{i,\ell}|_{U_{i,j,\ell}} \ \forall i, j, \ell \end{aligned}\right\}.
\end{equation}

A stack condition on $\cS$ relative to $\{U_i \to T\}$ is, roughly speaking,
a bijectivity condition on these two maps. More precisely, we have the
following definition.

\begin{defn} Given a category $\cS$ fibered in groupoids over a category 
$C$, and an arbitrary family of morphisms $\{U_i \to T\}$ in $C$, 
we say that $\cS$ satisfies the {\bf stack condition} relative to
$\{U_i \to T\}$ if the following holds:
\begin{ilist}
\itm (Morphisms form a sheaf) The map \eqref{mor-map} is always a 
bijection.
\itm (Objects have effective descent) The map \eqref{obj-map} is 
always surjective up to isomorphism on the right, where isomorphisms
are required to commute with the $\vp_{i,j}$.
\end{ilist}
\end{defn}

\begin{rem}\label{eff-desc} Even in this more general context, it is easy 
to check that given condition (i), the map \eqref{obj-map} is automatically 
injective if one mods out on the left by isomorphisms and on the right by 
isomorphisms commuting with the $\vp_{i,j}$. Thus, assuming condition (i) 
we have that condition (ii) is equivalent to a bijectivity statement for 
\eqref{obj-map} after modding out by the appropriate isomorphisms.
\end{rem}

In order to examine whether our deformation stacks are in fact stacks
for some Grothendieck (pre)topology, we are forced to consider two 
questions: first, is the condition we impose on the covers coming from
Schlessinger squares equivalent to the stack condition; and second, is
imposing the stack condition relative to those covers equivalent to 
imposing the stack condition relative to some Grothendieck
topology?

We address the first question first. Looking at the definitions, we see
that if we are given a Schlessinger square, it appears that we are 
imposing a stronger descent condition in the definition of deformation
stack, since, considering for example the condition on objects, we
require that a pair of objects on $A'$ and $A''$ descend as long as they
are isomorphic over $A' \otimes_A A''$, with no condition on isomorphisms
on the self-products $A' \otimes _A A'$ and $A'' \otimes _A A''$ 
satisfying the cocycle condition. Although technical, this is a serious
issue: for instance, in the \'etale topology it is crucial to use the
weaker descent condition, as the stronger one will never be satisfied.

In the Zariski topology, this is not a problem, since every covering
family $\{U_i \to T\}$ consists only of immersions, so 
$U_i \times _T U_i=U_i$. Similarly, in our situation the hypothesis that
$A \to A'$ is surjective means that $A' \otimes_A A' = A'$, so the only
difficulty arising in considering what happens on $A'' \otimes_A A''$.
We address this with the following lemma, which closely follows
Voevodsky's theory of regular cd-structures, introduced in \cite{vo2}.

\begin{lem} Let $A \to A', A \to A''$ be two morphisms in $\Art(\Lambda,k)$,
defining a Schlessinger square, and set $A'''=A' \otimes_A A''$. Let 
$I=\ker (A'' \otimes_A A'' \to A''' \otimes_{A'} A''')$, and let 
$\delta=\ker(A'' \otimes _A A'' \to A'')$.

Then we have $I \cap \delta=0$, so that the maps 
$$A'' \otimes _A A'' \to (A'' \otimes_A A'')/I \cong A''' \otimes_{A'} A'''$$
and
$$A'' \otimes _A A'' \to (A'' \otimes _A A'')/\delta \cong A''$$ 
define a Schlessinger square.
\end{lem}

\begin{proof} We first observe that both maps 
$$A'' \otimes _A A'' \to (A'' \otimes_A A'')/I \text{ and } 
A \to (A'' \otimes _A A'')/\delta$$
are surjective, so they define a Schlessinger square if and only if
$$A'' \otimes _A A'' \to (A'' \otimes_A A'')/I \times 
(A'' \otimes _A A'')/\delta$$
is injective, which is equivalent to the assertion that $I \cap \delta=0$. 

We now claim that any element of $I$ is of the form $\sum_i z_i w_i$,
with $w_i \in A'' \otimes_A A''$, and $z_i \in \ker(A \to A')$. 
Indeed, we have $A''' \otimes_{A'} A'''=
(A' \otimes_A A'') \otimes_{A'} (A' \otimes_A A'')=
A' \otimes_A (A'' \otimes_A A'')$. Since $A \to A'$ is surjective, $I$ is
simply the ideal generated by $\ker (A \to A')$, as claimed.

Now let $w=\sum_i z_i w_i$ be as above; we want to show that if 
$w \in \delta$, then we must have $w=0$. Expanding terms as necessary,
write $w_i=w'_i \otimes w''_i$ for each $i$. Since 
$z_i \in \ker(A \to A')$, by the conditions on a 
Schlessinger square we have that $z_i w'_i$ is the image of some 
$u_i \in A$, so we have $w = \sum_i 1 \otimes u_i w''_i=
1 \otimes \sum_i u_i w''_i$. Now, if $w \in \delta$ we have 
$\sum_i u_i w''_i=0$, so it follows that $w=0$. Thus $I \cap \delta=0$,
as desired.
\end{proof}

\begin{cor}\label{descent-equiv} A category $\cS$ fibered in groupoids over 
$\Art(\Lambda,k)$ is a deformation stack if and only if $\cS_k$ is trivial,
and $\cS$ satisfies the stack condition for every cover associated to a 
Schlessinger square.
\end{cor}

\begin{proof} It is easy to check from the definitions and Lemma
\ref{equiv-diags} that if $\cS$ is a deformation stack, then it satisfies 
a stack condition for every Schlessinger square. Thus, we need only prove 
the converse.

We want to see that for every Schlessinger square
$$\xymatrix{{X'''}\ar[r]^{q''}\ar[d]^{q'} & {X''} \ar[d] \\
{X'} \ar[r] & {X},}$$
and every pair of objects $\eta, \eta' \in \cS_X$, a pair of morphisms
$\vp' \in \Mor(\eta|_{X'},\eta'|_{X'})$ and 
$\vp'' \in \Mor(\eta|_{X''},\eta'|_{X''})$ which agree on $X'''$ 
automatically agree also on $X' \times_X X'$ and $X'' \times_X X''$, and
similarly for objects. 

Suppose we have such morphisms $\vp'$ and $\vp''$. Agreement on 
$X' \times_X X'$ is trivial; since $A \to A'$ is assumed surjective, we have 
that $X' \times_X X' = X'$, with both projections being the identity map.
For $X'' \times_X X''$, we invoke the lemma, finding that we have another
Schlessinger square defined by the maps
$$X'' \overset{\Delta}{\to} X'' \times_X X'' \text{ and }
X''' \times_{X'} X''' \to X'' \times_X X'',$$ 
where $\Delta$ denotes the diagonal map.
We wish to show that $p_1^*(\vp'')$ and $p_2^*(\vp'')$ agree on 
$X'' \times_X X''$. They certainly agree on $X''$. On the other hand, on
$X''' \times_{X'} X'''$ they give $(q'' \times q'')^* p_1^*(\vp'')$ 
and $(q'' \times q'')^* p_2^*(\vp'')$, which are the same as 
$p_i^* (q'')^* (\vp'')$ for $i=1,2$. Since we assumed that $\vp',\vp''$
agreed on $X'''$, these give $p_i^*(q')^*(\vp')$ for $i=1,2$, which are
equal since the product is over $X'$. We thus conclude that
$p_1^*(\vp'')=p_2^*(\vp'')$ after restriction to either $X''$ or
$X''' \times _{X'} X'''$, so the stack condition implies agreement
on $X'' \times_X X''$, as desired.

The situation is similar for objects: suppose we are given 
$\eta' \in \cS_{X'}$, and $\eta'' \in \cS_{X''}$, together with an 
isomorphism $\vp \in \Mor(\eta'|_{X'''}, \eta''|_{X'''})$. We want to
produce isomorphisms 
$\vp' \in \Mor_{X' \times_X X'}(p_1^* \eta',p_2^* \eta')$
and
$\vp'' \in \Mor_{X'' \times_X X''}(p_1^* \eta'',p_2^* \eta'')$
satisfying the cocycle condition, which in our case means that
$\vp' \circ \vp' = \vp'$, $\vp'' \circ \vp'' = \vp''$, 
$\vp \circ \vp'=\vp$, and $\vp'' \circ \vp = \vp$, all interpreted
suitably on the triple products. As before, there is nothing to do on
$X'$, since the products are canonically isomorphic to $X'$ itself, and
we can take $\vp'$ to be the identity. Also as before, we construct
$\vp''$ using the Schlessinger square provided by the lemma, and
using the fact that we already established the desired statement for
morphisms. Here we take $\vp''$ to restrict to the identity on $X''$ and
to $p_2^* \vp \circ p_1^* \vp^{-1}$ on $X''' \times_{X'} X'''$. 
We then have to check that $\vp'' \circ \vp'' = \vp''$
on $X'' \times_X X'' \times_X X''$. One first uses the same argument
as in the above lemma to see that 
$X''' \times_{X'} X''' \times_{X'} X''' \to X'' \times_X X'' \times_X X''$, 
and $X'' \times _X X''\overset{\id \times \Delta}{\to} 
X'' \times_X X'' \times_X X''$ together form a Schlessinger square, 
and it then suffices to check the
identity of morphisms after restriction to these, which is straightforward.

Finally, the identity $\vp'' \circ \vp = \vp$ on 
$X' \times_X X'' \times_X X''$ is easily checked via the observation
that $X' \times_X X'' \times_X X'' = X''' \times _{X'} X'''$.

Thus, if $\cS$ satisfies the condition for every cover arising from a
Schlessinger square, it is a deformation stack, completing the proof of 
the corollary.
\end{proof}

We see we have answered the first question positively, expressing
deformation stacks in terms of a descent condition relative to a certain
collection of covers. However, it turns
out that the second question, on the relationship between Schlessinger
squares and Grothendieck topologies, is less straightforward. The obvious
Grothendieck topology to consider is the one generated by covers coming
from Schlessinger squares. The problem with this arises from the condition 
that covers be stable under pullbacks, because it turns out that 
Schlessinger squares are not stable under pullback due to the 
condition on scheme-theoretic surjectivity.

The following example of Schlessinger squares failing to 
be preserved under pushforward is based on a suggestion of Eisenbud.

\begin{ex}\label{pullback-fails} We let $A=k[x,t]/(x^2,t^2,xt)$, 
$A'=k$, and $A''=k[y,t]/(y^2,t^2)$. The map $A \to A'$ is the obvious
quotient map, while the map $A \to A''$ is defined by sending $t$ to 
itself and $x$ to $yt$. We note that $A \to A'$ is surjective and 
$A \to A''$ is injective, and that $\ker(A \to A')=(x,t)$ maps to the
ideal $(t) \subseteq A''$. Thus, we obtain a Schlessinger square. 

Now let $B=k[x]/(x^2)$. We take the map $A \to B$ obtained by sending $t$ 
to $0$. If we push forward our square along $A \to B$, we obtain
$B'=k$, $B''=k[y]/(y^2)$, and the map $B \to B''$ sends $x$ to $0$. 
Therefore the condition $B \hookrightarrow B' \times B''$ fails, and
the pushforward is not a Schlessinger square. 
\end{ex}

We see from this example that we should not expect any deformation 
problem to satisfy a stack condition for the cover associated to this 
pushforward square: indeed, if the problem were prorepresentable, its 
tangent space would be automatically forced to be $0$ by such a condition. 
Thus, we see that it is not reasonable to hope to impose a stack condition
for the Grothendieck topology generated by Schlessinger squares, and
in particular such a condition is far stronger than the condition of 
being a deformation stack. 

The previous example might appear somewhat contrived, so one might 
naturally ask whether there is a ``good'' class of Schlessinger 
squares, stable under pullback, and which are nonetheless sufficient
for our purposes. However, we find that it is simple to modify the 
previous example to apply to even the most basic Schlessinger squares.

\begin{ex} Let $A=k[x,y]/(x^2,xy,y^2)$, equipped with the canonical quotient 
maps to $A'=k[x]/(x^2)$ and $A''=k[y]/(y^2)$. This is one of the most
fundamental Schlessinger squares, obtained from the projection maps of
$k[\epsilon] \times_k k[\epsilon]$. Let $B=k[x,t]/(x^2,t^2)$, and take
the map $A \to B$ which sends $x \mapsto x$ and $y \mapsto xt$. Pushing
forward our Schlessinger square as before, we check that scheme-theoretic
surjectivity once again fails to be preserved.
\end{ex}

These examples strongly suggest that the theory of deformation stacks
does in fact lie outside the usual framework of Grothendieck topologies.

\subsection{Examples}\label{subsec:first-exs}

We discuss several examples of deformation problems which naturally 
constitute deformation stacks.

\begin{ex}\label{ex:sheaves} Deformations of sheaves. Given a scheme 
$X_{\Lambda}$ over 
$\Spec \Lambda$, let $X:=X_{\Lambda}|_{\Spec k}$, and suppose we 
have $\cE$, a quasicoherent $\cO_X$-module. We consider the deformation 
problem classifying flat deformations of $\cE$ over $X_{\Lambda}|_A$ for 
different $A \in \Art(\Lambda,k)$. Our deformation stack 
$\cDef_{X_{\Lambda}}(\cE)$ as a category has objects consisting of triples 
$(A,\cE_A,\vp_A)$, where $A \in \Art(\Lambda,k)$, $\cE_A$ is an 
$\cO_{X_{\Lambda}|_A}$-module, flat over $A$, and $\vp_A:\cE_A \to \cE$ 
a morphism of $\cO_{X_A}$-modules such that the map
$\vp_A|_k:\cE_A \otimes_A k \to \cE$ induced by $\vp_A$ and $k \to \cE$ is 
an isomorphism. Morphisms 
$(A',\cE_{A'},\vp_{A'}) \to (A,\cE_A,\vp_A)$ consist of pairs $(f,\vp)$, 
with $f:A \to A'$, and $\vp:\cE_A \to \cE_{A'}$ a homomorphism of 
$\cO_{X_A}$-modules such that 
$\vp_{A'} \circ \vp = \vp_A$. The 
functor to $\Art(\Lambda,k)$ is then the forgetful one.

It is clear that $\cDef(X_{\Lambda},\cE)$ is a category fibered in groupoids 
over $\Art(\Lambda,k)$, but the fact that it is a deformation stack (like 
the fact that the associated functor satisfies Schlessinger's (H1) and (H2))
requires some justification. However, we see that essentially the same 
argument
will work. Indeed, if we fix $A' \to A$, $A'' \twoheadrightarrow A$, and
write $B=A' \times _A A''$, we need to check two facts. First, if we are 
given any objects $(B,\cE_B,\vp_B), (B,\cE_B',\vp_B')$, we abbreviate them
by $\cE_B$ and $\cE_B'$, and their restrictions to $A,A',A''$ similarly. We
then want that
$$\Mor_{X_B}(\cE_B,\cE_B')=\Mor_{X_{A'}}(\cE_{A'},\cE_{A'}') 
\times _{\Mor_{X_A}(\cE_A,\cE_A')} \Mor_{X_{A''}}(\cE_{A''},\cE_{A''}').$$
It is enough to work locally, and the desired
identity for morphisms of modules over $B,A',A'',A$ follows trivially from
the fact that the modules in question are free, by Corollary 
\ref{cor:mod-free}.
It then remains to check that a morphism of modules over $B$ is in fact
a morphism over $\cO_{X_B}$ if and only if the restrictions are morphisms
over $\cO_{X_{A'}}$ and $\cO_{X_{A''}}$, which similarly follows because
equality of elements of a free $B$-module may be checked after restriction
to $A'$ and $A''$.

Second, if we are given objects $(A',\cE_{A'},\vp_{A'})$ and 
$(A'',\cE_{A''},\vp_{A''})$, together with an isomorphism of their 
restrictions to $A$, we need to have an object $(B,\cE_B,\vp_B)$
inducing the given objects after pullback to $A'$ and $A''$, with the
given isomorphism after pullback to $A$.
Set $\cE_{B}=\cE_{A'} \times_{\cE_A} \cE_{A''}$. To check
the desired flatness and isomorphisms, it is enough to work locally on 
$X$, and the desired statements then follow from Lemma \ref{lem:schl2} below.
Using the natural map
$\cO_{X_B} \to \cO_{X_{A'}} \times_{\cO_{X_A}} \cO_{X_{A''}}$, we induce
an $\cO_{X_B}$-module structure on $\cE_{B}$, and quasicoherence can then 
be checked from the fact that module fiber product commutes with 
localization, so we obtain the desired descent condition.
\end{ex}

\begin{ex} Deformations of schemes. Given a scheme $X$ over $\Spec k$,
we consider the deformation problem classifying flat deformations of $X$
over $\Spec A$ for different $A \in \Art(\Lambda,k)$. Our deformation
stack $\cDef(X)$ has objects $(A,X_A,\vp_A)$, where $A \in \Art(\Lambda,k)$,
$X_A$ is a scheme flat over $A$, and $\vp_A:X \to X_A$ is a morphism over
the closed imbedding $\Spec k \to \Spec A$, inducing an isomorphism after
restriction to $\Spec k$. Morphisms 
$(A', X_{A'},\vp_{A'}) \to (A,X_A,\vp_A)$ consist of $(f,\vp)$, with 
$f:A \to A'$, and $\vp:X_{A'} \to X_A$ a morphism over $f$, 
inducing an isomorphism after restriction to $\Spec A'$, and with 
$\vp_A=\vp \circ \vp_{A'}$.

The proof that this gives a deformation stack is similar to the case of
sheaves. Indeed, the construction proceeds in the same fashion, except
that one has to work with algebras and morphisms of algebras, and it is
also necessary to check that the fiber product construction of the sheaf
case yields a scheme when applied to $\cO_{X_A},\cO_{X_{A'}},\cO_{X_{A''}}$.
As with checking quasicoherence above, this is simply a matter of fiber 
product of rings commuting with localization.
\end{ex}

\begin{ex} Deformations of quotient sheaves. Given a scheme $X_{\Lambda}$ 
over $\Spec \Lambda$, and a quasicoherent $\cO_{X_{\Lambda}}$-module 
$\cE_{\Lambda}$, set $X=X_{\Lambda}|_k$ and $\cE:=\cE_{\Lambda}|_k$. If 
we have $\cE \twoheadrightarrow \cF$ a quasicoherent quotient module of 
$\cE$, we have the deformation problem classifying deformations 
of $\cF$ as a quotient of $\cE_{\Lambda}|_A$ for $A \in \Art(\Lambda,k)$. 
Our deformation stack $\cDef_{X_{\Lambda},\cE_{\Lambda}}(\cF)$ has objects 
$(A,\cF_A,\vp_A)$, with $\cF_A$ flat over $A$, and a surjective map 
$\vp_A:\cE_{\Lambda}|_A \twoheadrightarrow \cF_A$ inducing an isomorphism
$\cF_A\otimes_A \Spec k \risom \cF$. Morphisms 
$(A',\cF_{A'},\vp_{A'}) \to (A, \cF_A, \vp_A)$ consist of $f:A \to A'$ such 
that the induced map 
$(\cE_{\Lambda} \otimes_{\Lambda} A) \otimes_A A' \to 
\cE_{\Lambda}\otimes_{\Lambda} A'$
yields an isomorphism
$\cF_A \otimes_A A' \risom \cF_{A'}$. 

Because there is at most one morphism lying over a given morphism of
$\Art(\Lambda,k)$, we need only check that given
$(A',\cF_{A'},\vp_{A'})$ and $(A'',\cF_{A''},\vp_{A''})$ both restricting to
some $(A,\cF_A,\vp_A)$, there exists a unique $(B,\cF_{B},\vp_B)$ restricting 
to the given objects over $A'$ and $A''$. But by flatness and the argument
of Example \ref{ex:sheaves}, we can take 
$\cF_{B}=\cF_{A'} \times_{\cF_{A}} \cF_{A''}$, and $\vp_B$ induced by
the natural map (which, despite the lack of flatness hypothesis on $\cE$, 
one checks is surjective) 
$\cE_B \twoheadrightarrow \cE_{A'} \times_{\cE_A} \cE_{A''}$. 
It is then easy to check that the construction of $(B,\cF_{B},\vp_B)$ 
provides an inverse to the natural map from quotients over $B$ to
pairs of quotients over $A'$ and $A''$ agreeing on $A$, so we obtain the
desired assertion.
\end{ex}

\begin{ex} Deformations of subschemes. Given a scheme $X_{\Lambda}$ over
$\Spec \Lambda$, set $X=X_{\Lambda}|_{\Spec k}$. Given $Z \subseteq X$
a closed subscheme, we have the deformation problem classifying deformations
of $Z$ as a subscheme of $X$. The deformation stack $\cDef_{X_{\Lambda}}(Z)$ 
has objects $(A,Z_A)$, with $A\in \Art(\Lambda,k)$, and $Z_A$ a closed 
subscheme of $X_{\Lambda}|_A$, flat over $A$, and such that $Z_A|_{\Spec k}=Z$.
Morphisms $(A', Z_{A'}) \to (A,Z_A)$ consist of $f:A \to A'$ such that 
$Z_{A'}|_A$ maps isomorphically 
to $Z_{A}$ under the natural map 
$(X_{\Lambda}|_{A'})|_A \to X_{\Lambda}|_{A}$.

In fact, this is a special case of deformations of quotient sheaves,
obtained by setting $\cE_{\Lambda}=\cO_{X_{\Lambda}}$ and $\cF=\cO_Z$; 
the kernel of any map of $\cO_{X_A}$-modules $\cO_{X_A} \to \cO_{Z_A}$ is
necessarily an ideal, so $\cO_{Z_A}$ inherits a unique algebra structure,
and the flatness of $Z$ over $A$ is equivalent to our condition for 
quotient sheaves.
\end{ex}

\begin{ex} Deformations of morphisms. Given a pair of schemes 
$X_{\Lambda},Y_{\Lambda}$ locally of finite type over $\Lambda$, 
with $X_{\Lambda}$ flat and $Y_{\Lambda}$ separated, set 
$X=X_{\Lambda}|_{\Spec k}, Y=Y_{\Lambda}|_{\Spec k}$. If we have also a
morphism $f:X \to Y$, we have the deformation problem classifying 
deformations of $f$. The deformation stack 
$\cDef_{X_{\Lambda},Y_{\Lambda}}(f)$ has objects $(A,f_A)$ with
$A \in \Art(\Lambda,k)$, and $f_A:X_{\Lambda}|_A \to Y_{\Lambda}|_A$
such that $f_A|_{\Spec k}=f$. Morphisms $(A',f_{A'}) \to (A,f_A)$ consist 
of $g:A \to A'$ such that $f_A|_{A'}=f_{A'}$.

We claim that with our hypotheses on $X$ and $Y$,
deformations of morphisms are a special case of deformations of closed 
subschemes, by considering the graph of the morphism. Because $Y$ is
separated, the graph is a closed subscheme, and it remains to check that 
any deformation of a graph
$\Gamma \subseteq X \times_k Y$ will still have $p_1:\Gamma_A \to X_A$ an
isomorphism, which follows from the flatness of $f$ by Corollary 17.9.5
of \cite{ega44}. Note also that the condition that $p_1$ be an isomorphism
immediately implies that for any deformation of $f$, we have $\Gamma_A$
flat over $A$.
\end{ex}

\begin{ex} Deformations of connections. Given a scheme $X_{\Lambda}$
smooth over $\Lambda$, and a quasicoherent $\cO_{X_{\Lambda}}$-module 
$\cE_{\Lambda}$ flat over $\Lambda$,
write $X=X_{\Lambda}|_{k}$ and $\cE=\cE_{\Lambda}|_k$. Suppose we have a 
connection $\nabla$ on $\cE$. We then have the deformation problem
classifying deformations of $\nabla$. The deformation stack 
$\cDef_{X_{\Lambda},\cE_{\Lambda}}(\nabla)$ has objects consisting of
$(A,\nabla_A)$, with $A \in \Art(\Lambda,k)$ and $\nabla_A$ a connection
on $\cE_{\Lambda}\otimes_{\Lambda} A$ such that $\nabla_A|_k=\nabla$. 
Morphisms $(A',\nabla_{A'}) \to (A,\nabla_A)$ consist of $f:A \to A'$ 
such that $\nabla_A|_{A'}=\nabla_{A'}$.

As in the case of quotient sheaves, there are no automorphisms, so we only 
need to check that any $\nabla_{A'}$ and $\nabla_{A''}$
agreeing over $A$ are obtained from a unique connection $\nabla_{B}$ over 
$B:=A' \times_A A''$. As in the case of deformations of sheaves, we have
$\cE_B=\cE_{A'} \times_{\cE_A} \cE_{A''}$, and we also see by smoothness
of $X_{\Lambda}$ that
$\Omega^1_{X_B/B}=\Omega^1_{X_{A'}/A'} 
\times_{\Omega^1_{X_A/A}} \Omega^1_{X_{A''}/A''}$. Since a connection
$\nabla_B:\cE_B \to \cE_B \otimes \Omega^1_{X_B/B}$ is in particular a map
of $B$-modules, using the same argument as before we find that we get a
unique map $\nabla_B$ of $B$-modules from $\nabla_{A'}$ and $\nabla_{A''}$.
Finally, checking that such a map is a connection if and only if 
$\nabla_{A'}$ and $\nabla_{A''}$ are connections is the same as checking 
$\cO_X$-linearity of morphisms in the case of deformations of sheaves.
\end{ex}

We mention briefly a number of additional examples of deformation stacks,
which will be studied in more detail in 
\cite{os15}. These involve combining previous examples: deformations
of schemes together with morphisms to a fixed scheme, deformations of 
connections together with the underlying sheaf, and deformations of 
subsheaves together with the ambient sheaf.

Indeed, it seems that every deformation problem of interest can be given
the structure of a deformation stack, and that furthermore, checking that
a problem satisfies Schlessinger's (H1) and (H2) entails checking the
axioms for a deformation stack. Thus, even though the conditions are 
technically more stringent than those of a deformation functor, in practice
they seem neither more restrictive nor harder to check.

\subsection{Additive structure}\label{subs:additive}

Throughout the following, we suppose that $\cS$ is a deformation stack,
and $f:A' \to A$ a small extension with kernel $I$. An important special 
case is when $V$ is a $k$-vector space, and we set $A'=k[V]$, $I=V$. In 
this case, we obtain an addition law on the category $\cS_{k[V]}$ in a 
strong sense. In full generality, we have a certain additive structure 
which we describe and explore. We first make some preliminary definitions.

\begin{notn} We denote by $T_I$ the set of isomorphism classes of
$\cS_{k[I]}$, and by $A_I$ the group $\Aut(\zeta_I)$ in $\cS_{k[I]}$.
\end{notn}

Thus, if $k[I]=k[\epsilon]$, then $T_I=T$ is the tangent space, and
$A_I=A$ is the infinitesimal automorphism group of $\cS$. The 
following notation deals with liftings of objects and morphisms over 
small extensions.

\begin{notn}
For a fixed $\eta \in \cS_A$, we write
$$T_{\eta,A'}:=
\{(\eta',\vp):\eta' \in \cS_{A'},\vp:\eta \to \eta'\text{ over }A' \to 
A\}/\cong,$$
where isomorphisms are applied simultaneously to $\eta'$ and $\vp$.

Similarly, fixing $\eta' \in \cS_{A'}$ and $\vp \in \Aut(\eta'|_A)$, we 
write 
$$A_{\vp,\eta'}:=\{\vp' \in \Aut(\eta'):\vp'|_A=\vp\}.$$
\end{notn}

We define below a categorical addition law 
$$+_{A',I}:\cS_{A'} \times \cS_{k[I]} \to \cS_{A'},$$ 
and our main result is then the following.

\begin{thm}\label{thm:def-aut-ext} The addition law $+_{k[I],I}$ gives 
canonical abelian group structures to $T_I$ and $A_I$, with the addition
on the latter agreeing with composition. These addition laws are 
functorial under linear maps $I \to J$.

If we fix $\eta \in \cS_A$, the addition law $+_{A',I}$ induces a
canonical action of $T_I$ on $T_{\eta,A'}$, making $T_{\eta,A'}$ into
a pseudotorsor for $T_{I}$.

Lastly, fixing $\eta' \in \cS_{A'}$, the addition law $+_{A',I}$ induces a 
canonical group isomorphism $A_I \risom A_{\id,\eta'}$. Therefore for
any $\vp \in \Aut(\eta'|_A)$, we have a canonical action of $A_I$ on
$A_{\vp,\eta'}$, making $A_{\vp,\eta'}$ a pseudotorsor for $A_{I}$.

Moreover, both pseudotorsor structures are functorial in the following 
sense: given another small extension $B' \to B$ with kernel $J$, and a 
homomorphism $A' \to B'$ which induces maps $I \to J$ and therefore
$A \to B$, then the additive structures commute with the induced 
restriction maps.
\end{thm}

Recall that we say that a set $S$ is a {\bf pseudotorsor} over a group $G$ 
if $S$ is either empty or a torsor over $G$. This theorem can be seen
as saying that if $\cS$ is a deformation stack with an obstruction space,
then we automatically obtain the additive part of a generalized 
tangent-obstruction theory in the sense of Fantechi and G\"ottsche; in
fact, we will see later that we obtain the full theory; see Remark 
\ref{rem:tan-obs} below.

The first step in defining $+_{A',I}$ is the following basic proposition, 
which essentially rephrases the deformation stack axioms in the form of 
operations on objects and morphisms. Indeed, this proposition is the only 
place in which we use the deformation stack axioms. Although we state the 
proposition quite generally, we will mainly be interested in the case 
$n=2$, with $\mathbf{A} = A' \times_k k[I]$ or $k[I] \times_k k[I]$. We 
will however use other cases in the proofs of certain technical statements 
below.

\begin{prop}\label{prop:ast-op} Fix $n \geq 2$, as well as
$A'_1,\dots,A'_n, A_1,\dots,A_{n-1} \in \Art(\Lambda,k)$ with surjective
maps $A'_i \to A_i, A'_{i+1} \to A_i$ for $i=1,\dots,n-1$, and set
$$\mathbf{A}=
A'_1 \times_{A_1} A'_2 \times_{A_2} \dots \times_{A_{n-1}} A'_n.$$

Given $\eta'_i \in \cS_{A'_i}$ for $i=1,\dots,n$ and $\eta_i \in \cS_{A_i}$
for $i=2,\dots,n-1$, with maps $\eta_i \to \eta'_i$ and 
$\eta_i \to \eta'_{i+1}$ lying over the given maps $A'_i \to A_i$ and
$A'_{i+1} \to A_i$ respectively, there exists an object 
$$\eta'_1 \ast_{\eta_1} \cdots \ast_{\eta_{n-1}} \eta'_n \in \cS_{\mathbf A}$$ 
having maps $q_i$ from $\eta'_i$, which lie over the projection 
maps $p_i$, and commuting with the given maps from each $\eta_i$. 
Furthermore, the tuple 
$(\eta'_1 \ast_{\eta_1} \dots \ast_{\eta_{n-1}} \eta'_n, q_1, \dots, q_n)$ 
is unique up to unique isomorphism. 

If we are given additionally $\mu'_i \in \cS_{A'_i}$ with maps from the
given $\eta_i$ as above, and morphisms $\vp'_i:\eta'_i \to \mu'_i$ 
in $\cS_{A'_i}$ commuting with the maps to the $\eta_i$,
there is a unique morphism 
$$\vp'_1 \ast_{\eta_1} \dots \ast_{\eta_{n-1}} \vp'_n:
\eta'_1 \ast_{\eta_1} \cdots \ast_{\eta_{n-1}} \eta'_n \to 
\mu'_1 \ast_{\eta_1} \cdots \ast_{\eta_{n-1}} \mu'_n$$
in $\cS_{\mathbf A}$ recovering the $\vp'_i$ under projection. 
The $\ast$ operation on morphisms commutes with composition.
If the choice
of $\eta'_1 \ast_{\eta_1} \cdots \ast_{\eta_{n-1}} \eta'_n$ or 
$\mu'_1 \ast_{\eta_1} \cdots \ast_{\eta_{n-1}} \mu'_n$ 
is changed, the resulting 
$\vp'_1 \ast_{\eta_1} \dots \ast_{\eta_{n-1}} \vp'_n$ changes by the 
corresponding unique isomorphism(s).

Finally, if for any $i$ we have $A_i=k$, then $\eta_i$ and the maps
$\eta_i \to \eta'_i$ and $\eta_i \to \eta'_{i+1}$ are irrelevant, so
we omit $\eta_i$ from our notation.
\end{prop}

\begin{proof}
The existence of 
$\eta'_1 \ast_{\eta_1} \cdots \ast_{\eta_{n-1}} \eta'_n$ together with
the $q_i$ follows immediately from induction on condition (ii) of a 
deformation stack. The existence of a unique isomorphism between any two 
tuples then follows similarly from condition (i). 

The existence and uniqueness of 
$\vp'_1 \ast_{\eta_1} \dots \ast_{\eta_{n-1}} \vp'_n$
similarly follows from condition (i) of the deformation stack axioms.
The fact that $\ast$ commutes with composition then follows formally from 
the definition, as does the behavior under change of the choice of 
$\eta'_1 \ast_{\eta_1} \cdots \ast_{\eta_{n-1}} \eta'_n$ or
$\mu'_1 \ast_{\eta_1} \cdots \ast_{\eta_{n-1}} \mu'_n$.

Finally, if $A_i=k$, because $\cS_{k}$ is trivial any choice 
of $\eta_i$ is uniquely isomorphic to $\zeta_0$, and there is always a 
unique morphism from $\eta_i$, so the conditions on restricting to 
$\eta_i$ and commuting with the given maps to $\eta_i$ are always satisfied
automatically, and we are justified in omitting $\eta_i$ from our notation.
\end{proof}

From the existence and uniqueness of $\eta_1 \ast \eta_2$, we can define
our additive structure using the ring map 
$$\sigma_{A',I}:A' \times_k k[I] \to A'$$
defined by $(x, \pi(x)+i) \mapsto x+y$.
We phrase the following corollary in a slightly stronger form than
is strictly necessary, so that we can conclude immediately that
$\cS_{k[I]}$ is a Picard category in the sense of Deligne (see 
Definition 1.4.2 of Expose XVIII of \cite{sga43}). This will not be
used here, but will be convenient in \cite{os15}.

\begin{cor}\label{c:add-struct} Appropriate choices of objects 
$\eta_1 \ast \eta_2$ and $\sigma_* (\eta_1 \ast \eta_2)$ for all
$\eta_1 \in \cS_{A'}, \eta_2 \in \cS_{k[I]}$ as well as, when
$A' \neq k[I]$, for all $\eta_1,\eta_2 \in \cS_{k[I]}$, gives us an 
additive structure consisting of functors
$$+_{A',I}:\cS_{A'} \times \cS_{k[I]} \to \cS_{A'}$$
and 
$$+_{k[I],I}:\cS_{k[I]} \times \cS_{k[I]} \to \cS_{k[I]}$$ 
and an isomorphism of functors 
$$\alpha:(\us +_{A',I} \us) +_{A',I} \us \risom 
\us +_{A',I} (\us +_{k[I],I} \us)$$
such that:
\begin{ilist}
\itm for any $\eta \in \cS_{k[I]}$, the functor 
$\us +_{A',I} \eta : \cS_{A'} \to \cS_{A'}$ is an equivalence of
categories;
\itm $\zeta_I \in \cS_{k[I]}$ is an identity object, in the sense that
$\us +_{A',I} \zeta_I: \cS_{A'} \to \cS_{A'}$ is the identity functor;
\itm $+_{k[I],I}$ is commutative, in the sense that 
$+_{k[I],I} = +_{k[I],I} \circ \sw$, where
$\sw: \cS_{k[I]} \times \cS_{k[I]} \to \cS_{k[I]} \times \cS_{k[I]}$
is the functor switching factors;
\itm for any $\eta_1,\eta_2,\eta_3,\eta_4 \in \cS_{k[I]}$, we have the 
identities
$$(\id_{\eta_1}+\alpha_{\eta_2, \eta_3,\eta_4}) \circ 
\alpha_{\eta_1,\eta_2+\eta_3,\eta_4} \circ
(\alpha_{\eta_1,\eta_2,\eta_3}+\id_{\eta_4}) =
\alpha_{\eta_1,\eta_2,\eta_3+\eta_4} \circ
\alpha_{\eta_1+\eta_2,\eta_3,\eta_4}$$
and
$$\alpha_{\eta_1,\eta_2,\eta_3} \circ \alpha_{\eta_3,\eta_1,\eta_2} 
= \alpha_{\eta_1,\eta_3,\eta_2}.$$
\end{ilist}

In the above and hereafter, we denote by
$\alpha_{\eta_1,\eta_2,\eta_3}: (\eta_1+\eta_2)+\eta_3 
\risom \eta_1 +(\eta_2 + \eta_3)$ the isomorphism obtained from $\alpha$.

Such an addition law is canonical up to unique isomorphism of functors.
It is also functorial in the following sense: given a small extension
$B' \to B$ with kernel $J$ and a map $A' \to B'$ inducing $I \to J$
and therefore $A \to B$, if
we have chosen $+_{A',I}$ and $+_{B', J}$ as above, then given
$\eta_1 \in \cS_{A'}$ and $\eta_2 \in \cS_{k[I]}$, and choices of
pushforwards $\eta_1|_{B'}$ and $\eta_2|_{k[J]}$, there is a
canonical map $\eta_1|_{B'}+\eta_2|_{k[J]} \to \eta_1+\eta_2$ 
lying over $A' \to B'$. In addition, under the induced identification
$(\eta_1+\eta_2)|_{B'} = \eta_1|_{B'}+\eta_2|_{k[J]}$, addition
of morphisms commutes with pushforward under $A' \to B'$.
\end{cor}

In what follows, we will frequently drop the subscripts on $\alpha$ and
$+$ when there is little likelihood of confusion. 

\begin{proof} We first claim that if we make arbitrary choices of 
$\eta_1 \ast \eta_2$ and $\sigma_*(\eta_1 \ast \eta_2)$ for every 
$\eta_1 \in \cS_{A'}, \eta_2 \in \cS_{k[I]}$, and likewise for
all $\eta_1,\eta_2 \in \cS_{k[I]}$, we obtain functors 
$+_{A',I}$ and $+_{k[I],I}$, with the isomorphism $\alpha$, and 
satisfying (i) as well as the first identity of (iv). The functors are 
given on objects by the choice
of $\sigma_*(\eta_1 \ast \eta_2)$, and by Proposition \ref{prop:ast-op}
we then also obtain an addition operation on morphisms. Compatibility
with composition of $\ast$ and pushforward implies that we obtain
functors, as asserted, and the fact that both $\ast$ and pushforward are
unique up to unique isomorphism implies that any two choices of $+_{A',I}$
differ by unique isomorphism. To check (i), for any $\eta \in \cS_{k[I]}$, 
we construct a functor $\us -_{A',I} \eta$ which is ``inverse''
to $\us +_{A',I} \eta$, in the sense that the composition on either
side is fully faithful, and satisfies 
$(\eta'+\eta)-\eta\cong(\eta'-\eta)+\eta \cong \eta'$ 
for all $\eta' \in \cS_{A'}$. We define $\us -_{A',I}\eta$ using the 
map $\delta:A' \times_k k[I] \to A'$ given by $\delta(x,\pi(x)+i)=x-i$. 
For every $\eta' \in \cS_{A'}$, choose a pushforward 
$\delta_*(\eta' \ast \eta)$, and define the functor on morphisms 
by $\vp \mapsto \delta_*(\vp \ast \id_{\eta})$. We then see that
\begin{align*}
(\us-\eta)+\eta & = (\delta_*(\us \ast \eta))+\eta \\
& \real \sigma_*(\delta \times \id)_* (\us \ast \eta \ast \eta) \\
& = \sigma_*(\delta \times \id)_* (\id \times \Delta)_* (\us \ast \eta) \\
& = p_{1 *} (\us \ast \eta),
\end{align*}
and
\begin{align*}
(\us+\eta)-\eta & = (\sigma_*(\us \ast \eta))-\eta \\
& \real \delta_*(\sigma \times \id)_* (\us \ast \eta \ast \eta) \\
& = \delta_*(\sigma \times \id)_* (\id \times \Delta)_* (\us \ast \eta) \\
& = p_{1 *} (\us \ast \eta),
\end{align*}
where $\Delta:k[I] \to k[I]\times_k k[I]$ is the diagonal map,
and $F \real G$ denotes that the functor $F$ is a particular 
realization of $G$, where $G$ is defined only up to unique isomorphism.
This proves the desired assertions.

The associativity isomorphism $\alpha$ is obtained similarly, except
that we are forced to keep track of our unique isomorphisms. Given
$\eta_1,\eta_2,\eta_3$, if we choose any $\eta_1 \ast \eta_2 \ast \eta_3$
with projection maps from $\eta_1,\eta_2,\eta_3$, we see that our prior
choices of $\eta_1 \ast \eta_2$ and $\eta_2 \ast \eta_3$, with
projection maps, induce unique maps 
$\eta_1 \ast \eta_2 \to \eta_1 \ast \eta_2 \ast \eta_3$ and
$\eta_2 \ast \eta_3 \to \eta_1 \ast \eta_2 \ast \eta_3$ over $p_{12}$
and $p_{23}$ respectively: indeed, any 
$p_{12*}(\eta_1 \ast \eta_2 \ast \eta_3)$ gives a choice of 
$\eta_1 \ast \eta_2$, which then differs by a unique isomorphism from our
given choice, and composing the pushforward with the isomorphism gives us 
the desired map. The same works for $p_{23}$. A similar argument produces
unique maps 
$(\eta_1+\eta_2)\ast \eta_3 \to \eta_1 \ast \eta_2 \ast \eta_3$
and 
$\eta_1 \ast (\eta_2 + \eta_3) \to \eta_1 \ast \eta_2 \ast \eta_3$
lying over $\sigma \times \id$ and $\id \times \sigma$ respectively,
and compatible with the previously constructed maps. 
We thus see that 
\begin{align*}
(\eta_1+\eta_2)+\eta_3& =\sigma_*(\sigma_*(\eta_1 \ast \eta_2)\ast \eta_3) \\
& \real 
\sigma_*((\sigma \times \id)_* (\eta_1 \ast \eta_2 \ast \eta_3)) \\
& = \sigma^{(2)}_*(\eta_1 \ast \eta_2 \ast \eta_3),
\end{align*} 
where 
$\sigma^{(2)}:A' \times_k k[I] \times_k k[I] \to A'$ is the map defined
by 
$$\sigma^{(2)} (x,\pi(x)+i_1,\pi(x)+i_2) =x+i_1+i_2.$$
But $\eta_1+(\eta_2+\eta_3)$ is then a different realization of 
$\sigma^{(2)}_*(\eta_1 \ast \eta_2 \ast \eta_3)$, so we obtain a 
unique isomorphism $\alpha_{\eta_1,\eta_2,\eta_3}$, which one checks
is independent of the choice of $\eta_1 \ast \eta_2 \ast \eta_3$ and gives 
the desired isomorphism of functors.
The first identity of (iv) is checked similarly, by expressing both
$(\id_{\eta_1}+\alpha_{\eta_2, \eta_3,\eta_4}) \circ 
\alpha_{\eta_1,\eta_2+\eta_3,\eta_4} \circ
(\alpha_{\eta_1,\eta_2,\eta_3}+\id_{\eta_4})$ and 
$\alpha_{\eta_1,\eta_2,\eta_3+\eta_4} \circ
\alpha_{\eta_1+\eta_2,\eta_3,\eta_4}$ as the unique isomorphisms between
two different pushforwards of the quadruple product 
$\eta_1 \ast \eta_2 \ast \eta_3 \ast \eta_4$.

It remains to show that there exist choices of $\eta_1 \ast \eta_2$ and
$\eta_1 + \eta_2$ satisfying (ii), (iii), and the second identity of
(iv) (which makes sense only if (iii) is assumed). We first observe
that given (iii), the second identity of (iv) is automatically satisfied
for the same reason as the first: both 
$\alpha_{\eta_1,\eta_2,\eta_3} \circ \alpha_{\eta_3,\eta_1,\eta_2}$ and
$\alpha_{\eta_1,\eta_3,\eta_2}$ can be realized as the unique isomorphisms
between two choices of pushforwards of $\eta_1 \ast \eta_2 \ast \eta_3$.
Next, for (ii) we first
observe that because $\zeta_I = s_* (\zeta_0)$ where $s:k \to k[I]$ is 
the structure map, we always have (up to unique isomorphism)
$$\eta' + \zeta_I = \sigma_*(\eta' \ast (s_* \zeta_0)) =
\sigma_* (\id \times (s \circ \pi))_* (\eta') = \id_* \eta' =\eta',$$
so if we realize $\eta' \ast \zeta_I$ as $(\id \times s)_* \eta'$ in
such a way that the composition map $\eta' \to \eta' \ast \zeta_I \to \eta'$
over 
$A' \overset{\id \times s}{\to} A' \times_k k[I] \overset{p_1}{\to} A'$
is the identity, we obtain a unique map $\eta' \to \eta' \ast \zeta_I$
over $\sigma$ such that the composition map 
$\eta' \to \eta' \ast \zeta_I \to \eta'$
over 
$A' \overset{\id \times s}{\to} A' \times_k k[I] \overset{\sigma}{\to} A'$
is also the identity. One then checks that with this choice of
$\eta' \ast \zeta_I$ for each $\eta'$, given $\vp:\eta'_1 \to \eta'_2$
we have $\vp \ast \id_{\zeta_I}=(\id \times s)_* \vp$, and therefore
that we have that $\us + \zeta_I$ is the identity functor, as desired.

We now verify that when $A'=k[I]$, we may further make the necessary
choices in such a way that commutativity holds. For this, we place a 
total ordering on the objects of $\cS_{k[I]}$ such that $\zeta_I$ is
minimal. Then, for every $\eta'_1 \geq \eta'_2$, we fix choices of
$\eta'_1 \ast \eta'_2$ and $\sigma_*(\eta'_1 \ast \eta'_2)$ arbitrarily,
except that when $\eta'_2=\zeta_I$ we make our choice as dictated above.
For $\eta'_1 < \eta'_2$, we finally set $\eta'_1 \ast \eta'_2$ to
be a choice of pushforward $\sw_{I *} \eta'_2 \ast \eta'_1$, where
$\sw_{I}:k[I]\times_k k[I] \to k[I] \times_k k[I]$ is the map switching
factors. We then choose $\eta'_1 + \eta'_2 = \eta'_2+\eta'_1$ as an
object, with the map $\eta'_1 \ast \eta'_2 \to \eta'_1+\eta'_2$ over
$\sigma$ chosen so that the map 
$\eta'_2 \ast \eta'_1 \to \eta'_1 \ast \eta'_2$ over $\sw_I$ pushes
forward under $\sigma$ to the identity map. One checks via a diagram
chase that such a choice makes $+_{A',I}$ commutative not only on the
level of objects but also of morphisms: given $f_i:\eta'_i \to \mu'_i$,
if we denote by $\sw_{I*}$ both pushforward (iso)morphisms 
$\eta'_1 \ast \eta'_2 \to \eta'_2 \ast \eta'_1$ and 
$\mu'_1 \ast \mu'_2 \to \mu'_2 \ast \mu'_1$ over $\sw_I$, one checks
that $\sw_{I*}^{-1} \circ (f_2 \ast f_1) \circ \sw_{I*}$ satisfies
the characterization of $f_1 \ast f_2$, which yields the desired result.

Finally, our functoriality assertion follows from basic properties
of categories fibered in groupoids, and the uniqueness of $\ast$ in
Proposition \ref{prop:ast-op}. Indeed, if we have $\eta_1 \ast \eta_2$
with maps from $\eta_1$ and $\eta_2$, and choices of 
$\eta_1|_{B'} \to \eta_1$ and $\eta_2|_{k[J]} \to \eta_2$
one checks directly from the definition of a category fibered in groupoids
that for any choice of 
$(\eta_1 \ast \eta_2)|_{B' \times_k k[J]} \to \eta_1 \ast \eta_2$
there exist unique maps from $\eta_1|_{B'}$ and $\eta_2|_{k[J]}$
to $(\eta_1 \ast \eta_2)|_{B' \times_k k[J]}$ making everything commute.
From Proposition \ref{prop:ast-op} we then obtain a unique isomorphism
$(\eta_1 \ast \eta_2)|_{B' \times_k k[J]} \to 
\eta_1|_{B'} \ast \eta_2|_{k[J]}$, and composing with this 
isomorphism gives us a map 
$\eta_1|_{B'} \ast \eta_2|_{k[J]} \to \eta_1 \ast \eta_2$ over
$A' \times_k k[I] \to B' \times_k k[J]$. Using the definition
of a category fibered in groupoids once again, we obtain the desired unique 
map $\eta_1|_{B'} + \eta_2|_{k[J]} \to \eta_1 +\eta_2$ over
$A' \to B'$, commuting with the prechosen pushforward maps
$\eta_1 +\eta_2 \to \eta_1 \ast \eta_2$ and 
$\eta_1|_{B'} + \eta_2|_{k[J]} \to \eta_1|_{B'} \ast \eta_2|_{k[J]}$. 
One then checks by chasing definitions
that the desired functoriality for morphisms follows from the commutativity
of all the involved maps.
\end{proof}

This says, roughly speaking, that $\cS_{k[V]}$ is always an ``abelian
group category'', and that $\cS_{k[I]}$ acts naturally on $\cS_{A'}$. We
next show that in fact the action of $\cS_{k[I]}$ factors through 
a collection of much ``smaller'' auxiliary categories $\cS_{\eta,A'}$ 
associated to $\cS_{A'}$, making each of the $\cS_{\eta,A'}$ 
into ``pseudotorsor categories'' for $\cS_{k[I]}$.

\begin{defn} Given $\eta \in \cS_{A}$, let $\cS_{\eta,A'}$ be the
category whose objects consist of pairs $(\eta',\vp)$ with
$\eta' \in \cS_{A'}$, and $\vp:\eta \to \eta'$ lying over $A' \to A$,
and with morphisms $(\eta'_1,\vp_1)\to (\eta'_2,\vp_2)$ given by
$\psi:\eta'_1 \to \eta'_2$ in $\cS_{A'}$ with $\vp_2 \circ \psi=\vp_1$.
\end{defn}

Our basic result is then the following:

\begin{prop}\label{p:torsor-cat} For any $\eta \in \cS_{A}$, the functor
$+_{A',I}$ of Corollary \ref{c:add-struct} induces a natural functor 
(which we still denote by $+_{A',I}$) on $\cS_{\eta,A'}$,
with any two choices of $+_{A',I}$ differing by an isomorphism
of functors. The associativity isomorphism $\alpha$ induces a new
isomorphism of functors which we still denote by $\alpha$, and
$\us +_{A',I} \zeta_I$ is still the identity functor. 

Furthermore, if we choose any object 
$(\eta',\vp) \in \cS_{\eta,A'}$, the induced functor
$$(\eta',\vp)+_{A',I} \us: \cS_{k[I]} \to \cS_{\eta,A'}$$
is an equivalence of categories.
\end{prop}

\begin{proof} To define the functor induced by $+_{A',I}$ on $\cS_{\eta,A'}$, 
suppose we have $(\eta'_1,\vp_1)\in \cS_{\eta,A'}$, and 
$\eta_2 \in \cS_{k[I]}$. We set 
$(\eta'_1,\vp_1)+\eta_2=(\eta'_1+\eta_2,\vp_{12})$, where $\vp_{12}$ is
uniquely determined from the axioms of a category fibered in groupoids
by the condition that the composed map
$\eta \overset{\vp_1}{\to} \eta'_1 \overset{q_1}{\to} \eta'_1 \ast \eta_2$ 
over $A' \times_k k[I] \overset{p_1}{\to} A' \to A$ agrees with
$\eta \overset{\vp_{12}}{\to} \eta'_1+\eta_2 \to \eta'_1 \ast \eta_2$ 
over $A' \times_k k[I] \overset{\sigma}{\to} A' \to A$. One checks easily
that the addition of morphisms given by $+_{A',I}$ makes our new
addition law on $\cS_{\eta,A'}$ into a functor.
Moreover, we already know that two different choices of $+_{A',I}$ differ 
by a unique isomorphism of functors on $\cS_{A'}$, and it is routine to
check that the relevant isomorphisms in $\cS_{A'}$ are in fact isomorphisms 
in $\cS_{\eta,A'}$, so it follows formally that we obtain isomorphisms of 
any two choices of the new addition functor as well.
Similarly, we check directly from our construction that
the isomorphisms in $\cS_{A'}$ given by $\alpha$ are isomorphisms in
$\cS_{\eta,A'}$, and that $\zeta_I$ still acts as the identity.

To prove the equivalence assertion, we again construct the ``inverse''
functor $-(\eta',\vp)+\us$ explicitly. Given 
$(\eta',\vp) \in \cS_{\eta,A'}$, 
and any $(\eta'_1,\vp_1)\in \cS_{\eta,A'}$, we obtain 
$-(\eta',\vp)+(\eta'_1,\vp_1) \in \cS_{k[I]}$ as follows: we consider 
$\eta' \ast_{\eta} \eta'_1$,
and write $\delta:A' \times_A A' \to k[I]$ for the map
given by $(x,y) \mapsto \pi(x)+(y-x)$. We then set 
$-(\eta',\vp)+(\eta'_1,\vp_1)=\delta_*(\eta' \ast_{\eta} \eta'_1)$. 
Here we fix any choices of $\eta' \ast_{\eta} \eta'_1$ and the resulting
pushforwards.  Given also $(\eta'_2,\vp_2)$, if we have a morphism 
$\psi:\eta'_1 \to \eta'_2$ commuting with the $\vp_i$, we obtain 
$\id_{\eta'} \ast_{\eta} \psi$,
and pushing forward under $\delta$ gives us a morphism 
$-(\eta',\vp)+(\eta'_1,\vp_1) \to -(\eta',\vp)+(\eta'_2,\vp_2)$.
Compatibility of $\ast_{\eta}$ and pushforward with composition 
demonstrates that this defines a functor from $\cS_{k[I]}$ to 
$\cS_{\eta,A'}$.

It remains to check that $- (\eta',\vp)+\us$ is ``inverse'' to
$(\eta',\vp)+\us$, in the same sense as in the proof of 
Corollary \ref{c:add-struct} (i). But we see that
\begin{align*}
(\eta',\vp)+(-(\eta',\vp)+\us) 
& = 
(\eta',\vp)+(\delta_*((\eta',\vp) \ast_{\eta} \us)) \\
& = \sigma_*((\eta',\vp) \ast (\delta_*((\eta',\vp) \ast_{\eta} \us))) \\
& \real \sigma_*(\id \times \delta)_* 
((\eta',\vp) \ast (\eta',\vp) \ast_{\eta} \us) \\
& = \sigma_*(\id \times \delta)_* 
((\Delta \times \id)_* ((\eta',\vp) \ast_{\eta} \us)) \\
& = p_{2*} ((\eta',\vp) \ast_{\eta} \us),
\end{align*}
and
\begin{align*}
-(\eta',\vp)+((\eta',\vp)+\us)
& = -(\eta',\vp)+(\sigma_*((\eta',\vp) \ast \us)) \\
& = \delta_*((\eta',\vp)\ast_{\eta} (\sigma_*((\eta',\vp) \ast \us))) \\
& \real \delta_*(\id \times \sigma)_*
((\eta',\vp)\ast_{\eta} (\eta',\vp) \ast \us) \\
& = \delta_*(\id \times \sigma)_*
((\Delta' \times \id)_* ((\eta',\vp) \ast \us)) \\
& = p_{2*} ((\eta',\vp) \ast \us),
\end{align*}
where $\Delta:A' \to A' \times_k A'$ and $\Delta':A' \to A' \times_A A'$
are the diagonal maps. Thus we conclude the desired statement.
\end{proof}

As indicated by Theorem \ref{thm:def-aut-ext}, this has many consequences, 
on the level of isomorphism classes of objects and on the level of 
(auto)morphisms. We next make a closer study of the behavior of our 
additive structure on morphisms. Here we mention that several of our 
statements, including that automorphism groups are abelian and canonically 
identified with one another, were stated by Grothendieck in a similar 
setting in \S 1.5 of \cite{gr4}. He did not include any arguments,
saying that the statements were ``without doubt well known to category 
theorists.''

\begin{cor}\label{cor:mors-group} Given $\eta \in \cS_A$, we have:
\begin{ilist}
\itm given two morphisms 
$f:\eta_1 \to \eta_2$, $g:\eta_2 \to \eta_3$ in $\cS_{k[I]}$, we
have the identity $g+_{k[I],I} f= g \circ f +_{k[I],I} \id_{\eta_2}$;
\itm in $\Aut(\zeta_I)$, composition agrees with $+_{k[I],I}$, and
in particular $\Aut(\zeta_I)$ is an abelian group;
\itm given $\eta' \in \cS_{\eta,A'}$, the map 
$\Aut(\zeta_{I}) \to \Aut(\eta')$
induced by $\vp \mapsto \id_{\eta'}+_{A',I} \vp$ is an isomorphism of
groups under composition. In particular, $\Aut(\eta')$ is abelian;
\itm when $\eta' \cong \eta''$ in $\cS_{\eta,A'}$, the identification
of $\Aut(\eta')$ with $\Aut(\eta'')$ obtained via the above mutual
isomorphisms with $\Aut(\zeta_I)$ agrees with the usual identification
obtained by conjugating with a choice of isomorphism $\eta' \risom \eta''$;
\itm given $\eta_1 \in \cS_{\eta,A'}$, $\eta_2,\eta_3 \in \cS_{k[I]}$, and 
$\vp_i \in \Aut(\eta_i)$ for $i=1,2,3$, 
$$(\vp_1+\vp_2)+\vp_3=\vp_1+(\vp_2+\vp_3)$$
under the above identification of
$\Aut((\eta_1+\eta_2)+\eta_3)$ with $\Aut(\eta_1+(\eta_2+\eta_3))$;
\itm given $\eta_1 \in \cS_{\eta,A'}$ and $\eta_2 \in \cS_{k[I]}$, and 
$\vp_1,\vp_2 \in \Aut(\zeta_I)$, we have 
$$(\id_{\eta_1+\eta_2})+(\vp_1+\vp_2)=
(\id_{\eta_1}+\vp_1)+(\id_{\eta_2}+\vp_2)$$
as elements of $\Aut(\eta_1+\eta_2)$.
\end{ilist}

Note that in the above, automorphism groups for objects of $\cS_{\eta,A'}$
are the automorphism groups in $\cS_{\eta,A'}$, and not the automorphism
groups of the underlying object of $\cS_{A'}$. That is, we consider
automorphisms which restrict to the identity in $\Aut(\eta)$.
\end{cor}

\begin{proof} For (i), using commutativity and the fact that 
$+_{k[I],I}$ is a functor, we have
$$g \circ f + \id_{\eta_2}=(g+\id_{\eta_2})\circ(f+\id_{\eta_2})=
(g+\id_{\eta_2}) \circ (\id_{\eta_2}+f)=(g \circ \id_{\eta_2}) + 
(\id_{\eta_2} \circ f) = g + f.$$
Statement (ii) then immediately follows from (i), since 
$+_{k[I],I} \id_{\zeta_I}$ acts as the identity on morphisms.

(iii) follows immediately from Proposition \ref{p:torsor-cat}, since if
$\eta'+_{A',I}:\cS_{\eta,A'} \to \cS_{k[I]}$ is an equivalence of
categories, it necessarily induces isomorphisms on automorphism groups.

To see (iv) we suppose we are given $f:\eta' \risom \eta''$,
and $\vp' \in \Aut(\eta')$. Let $\vp \in \Aut(\zeta_I)$ be the unique
automorphism with $\id_{\eta'}+_{A',I} \vp=\vp'$, and set
$\vp'' \in \Aut(\eta'')$ to be $\id_{\eta''}+_{A',I} \vp$. We wish to
see that $\vp' = f^{-1} \circ \vp'' \circ f$. But using again that
$+_{k[I],I}$ is a functor, we have 
\begin{align*}
f^{-1} \circ \vp'' \circ f 
& = (f^{-1}+_{A',I} \id_{\zeta_I}) \circ (\id_{\eta''} +_{A',I} \vp) 
\circ (f+_{A',I} \id_{\zeta_I}) \\
& =(f^{-1} \circ \id_{\eta''} \circ f)+_{A',I} \vp \\
& = \id_{\eta'} +_{A',I} \vp = \vp',
\end{align*}
as desired.

(v) follows trivially from (iv), since we have 
$$(\vp_1+\vp_2)+\vp_3=\alpha_{\eta_1,\eta_2,\eta_3}^{-1} \circ
(\vp_1+(\vp_2+\vp_3)) \circ \alpha_{\eta_1,\eta_2,\eta_3}.$$
Finally, (vi) then follows by repeated application of (v), after
expanding $\id_{\eta_1+\eta_2}=\id_{\eta_1}+\id_{\eta_2}$. 
\end{proof}

In particular, we see that the ``linearization'' imposed by the 
deformation stack condition means that the infinitesimal automorphism 
group $A_{\cS}$ is necessarily abelian.

We can now easily give the proof of Theorem \ref{thm:def-aut-ext}.

\begin{proof}[Proof of Theorem \ref{thm:def-aut-ext}] We first address
$T_I$. The operation $+_{k[I],I}$ is commutative with identity 
$\zeta_I$ by construction. That it is associative on the level of 
isomorphism classes is clear from the existence of the isomorphisms
$\alpha_{x,y,z}$. Moreover, the existence of inverses follows from the
fact that addition of a given object gives an equivalence of categories.
The structure is canonical because any two choices of the 
functor $+_{k[I],I}$ differ by isomorphisms, so induce the same 
operation on $T_I$. The desired statements for $A_I$ are simply 
Corollary \ref{cor:mors-group} (ii). Finally, functoriality for both 
$T_I$ and $A_I$ follows from the functoriality statement of Corollary 
\ref{c:add-struct}. 

The assertions for $T_{\eta,A'}$ then follow easily from Proposition 
\ref{p:torsor-cat}. Indeed, $T_{\eta,A'}$ is simply the set of
isomorphism classes of $\cS_{\eta,A'}$, so we easily see that we have
a map $T_{\eta,A'} \times T_I \to T_{\eta,A'}$, with $\zeta_I$ acting
as the identity, and the existence of $\alpha$ implies that this map
indeed defines a group action.
Lastly, if $T_{\eta,A'}$ is non-empty, we see that it is in fact a
torsor from the assertion of Proposition \ref{p:torsor-cat} that
addition induces an equivalence of categories.

We next address the corresponding statements for automorphisms.
We first note that the general case follows immediately from the
case that $\vp=\id$, since for any $\vp$ we have that $A_{\vp,\eta'}$
is naturally a pseudotorsor for $A_{\id,\eta'}$. But the isomorphism
$A_I \risom A_{\id,\eta'}$ is precisely Corollary \ref{cor:mors-group}
(iii). 

Functoriality in both cases follows immediately from the functoriality
assertion of Corollary \ref{c:add-struct}. 
\end{proof}

\subsection{Linear structure}

Given a $k$-vector space $V$, having analyzed the additive structures 
associated to $\cS_{k[V]}$ for a deformation stack $\cS$,
we may similarly define scalar multiplication maps over $k$, making
$T_V$ and $A_V$ into $k$-vector spaces. Although it is 
possible to do this on a categorical level as we did with addition,
expressing the proper conditions for associativity and distributivity 
isomorphisms becomes substantially more complicated, so we settle for
the simpler task of defining vector space structures on $T_V$ and 
$A_V$.

\begin{defn} Given $\lambda \in k$, let $m_{\lambda}:k[V] \to k[V]$ be
defined by $a+v \mapsto a+\lambda v$, for $a \in k, v \in V$. We then
define maps 
$m_{\lambda *}:T_{V} \to T_{V}$ 
and
$m_{\lambda *}:A_{V} \to A_{V}$ 
by pushing forward under $m_{\lambda}$.
\end{defn}

These maps define vector space structures:

\begin{cor}\label{cor:def-aut-vs} The sets $T_{V}$ and $A_{V}$ carry
canonical vector space structures defined by $+_{k[V],V}$ and the 
$m_{\lambda}$ maps, with the addition law on
$A_{V}$ agreeing with composition. 
\end{cor}

\begin{proof} We have already seen in Theorem \ref{thm:def-aut-ext} 
that $+_{k[V],V}$ gives canonical abelian group structures on $T_V$ and 
$A_V$ and agrees with composition for $A_V$, so it remains to check
that the scalar multiplication maps are canonically defined, and
satisfy the appropriate vector space axioms. That they are canonically
defined for $T_V$ is clear, since any two choices of $m_{\lambda *} \eta$
differ by unique isomorphism. For $A_V$, we define $m_{\lambda *}$ more
precisely by choosing a pushforward map $\zeta_V \to \zeta_V$ over
$m_{\lambda}$ (which we observe we can do). This pushforward map is
well-defined only up to composition by an element of $A_V$, but the
resulting pushforward of automorphisms is then conjugated by that element,
which leaves the map on $A_V$ unchanged because $A_V$ is abelian.

We next check the vector space axioms on $T_V$. It is clear that for
$\eta \in T_V$, we have $m_{1 *} \eta=\eta$, since $m_1$ is the
identity map. Similarly, since 
$m_{\lambda \lambda'} = m_{\lambda} \circ m_{\lambda'}$, we have
\begin{align*}
m_{\lambda \lambda' *} \eta & = (m_{\lambda} \circ m_{\lambda'})_* \eta \\
& = m_{\lambda *} (m_{\lambda' *} \eta).\end{align*} 
Next, given given $\lambda,\lambda' \in k$ and $\eta \in T_V$, we have
\begin{align*}
(m_{\lambda *} \eta)+(m_{\lambda' *} \eta) 
& = \sigma_* ((m_{\lambda *} \eta)\ast (m_{\lambda' *} \eta)) \\
& = \sigma_* (m_{\lambda,\lambda' *} (\eta \ast \eta)) \\
& = \sigma_* (m_{\lambda,\lambda' *} (\Delta_* \eta)) \\ 
& = m_{\lambda+\lambda' *} (\eta),
\end{align*}
where $\Delta:k[V] \to k[V] \times_k k[V]$ is the diagonal imbedding
and $m_{\lambda,\lambda'}$ is the map 
$k[V]\times_k k[V] \to k[V] \times_k k[V]$ given by $m_{\lambda}$
on the first factor and $m_{\lambda'}$ on the second.

Similarly, given
$\lambda \in k$ and $\eta,\eta' \in T_V$, we have 
\begin{align*}
m_{\lambda *}(\eta + \eta')& =m_{\lambda *}(\sigma_* (\eta \ast \eta')) \\
& =\sigma_*(m_{\lambda,\lambda *} (\eta \ast \eta')) \\
& =\sigma_*((m_{\lambda *} \eta) \ast (m_{\lambda *} \eta')) \\
& = (m_{\lambda *}\eta)+(m_{\lambda *} \eta').
\end{align*}

Finally, the vector space axioms on $A_V$ follow from precisely the same
arguments. 
\end{proof}

The following lemma gives an important functoriality statement.

\begin{lem}\label{lem:linear} A linear map $V \to W$ induces a linear map 
$T_{V} \to T_{W}$ under pushforward. Moreover, in this way we 
obtain a linear map $\Hom(V,W) \to \Hom(T_{V},T_{W})$.

Finally, the same is true if we replace $T_{V}$ and $T_{W}$ by 
$A_{V}$ and $A_{W}$.
\end{lem}

\begin{proof} Indeed, the linearity of the map $V \to W$ implies that the
induced maps commute with the $\sigma$ and $m_{\lambda}$ used to define
the vector space structure, and hence give linear maps $T_V \to T_W$ by 
simple diagram chases.

It is clear from the definitions that the induced map 
$\Hom(V,W) \to \Hom(T_V,T_W)$ commutes with scalar multiplication, so we
need only check that it also commutes with addition. That is, given 
$\eta \in T_V$, and $\vp_1,\vp_2 \in \Hom(V,W)$, we need to check that
$$(\vp_1+\vp_2)_*(\eta)=\vp_{1*} (\eta)+\vp_{2*} (\eta).$$
But this is easy enough: write $\eta'$ for the object of
$F_{\cS}(k[W] \times_k k[W])$ inducing $\vp_{1*}(\eta)$ and $\vp_{2*}(\eta)$
under the projection maps, so that 
$\vp_{1*}(\eta)+\vp_{2*}(\eta)=\sigma_* (\eta')$. Then one need only note
that $\eta'$ is the pushforward of $\eta$ under the map 
$k[V] \to k[W] \times_k k[W]$ induced by $\vp_1$ and $\vp_2$, and that the
composition of this map with $\sigma$ is precisely $\vp_1+\vp_2$.

Finally, the argument in the case of $A_V,A_W$ proceeds 
in exactly the same fashion.
\end{proof}

We can now prove the following explicit descriptions of $T_V$ and $A_V$.

\begin{prop}\label{prop:tensors}
For any $k$-vector space $V$ we have canonical $k$-vector space 
isomorphisms 
$$T_{V} \cong T_{\cS} \otimes_k V$$
and
$$A_{V} \cong A_{\cS} \otimes_k V,$$
functorial in $V$.
\end{prop}

\begin{proof}
We can define a natural map 
$T_{\cS} \otimes_k V \to T_{V}$ induced by sending $\eta \otimes v$ to the 
image of $\eta$ under $\vp_v:k[\epsilon] \to k[V]$ defined by 
$\epsilon \mapsto v$. One checks using Lemma \ref{lem:linear} that this 
gives for fixed $v$ a linear map $T_{\cS} \to T_{V}$, and for fixed 
$\eta$ a linear map $V \to T_{V}$, so we obtain a well-defined linear
map from the tensor product.
To check that we have an isomorphism, we need only check bijectivity.
We choose a basis $v_1,\dots,v_d$ of $V$, and factor our map through the
isomorphism 
$T_{\cS} \otimes V \risom T_{(v_1)} \times \dots \times T_{(v_d)}$, 
using the isomorphisms $T_{\cS} \risom T_{(v_i)}$ induced by $\vp_{v_i}$.

The map $T_{(v_1)} \times \dots \times T_{(v_d)} \to T_{V}$ is 
then induced
by the inclusions $(v_i) \to V$ together with the $d$-fold summation map 
$\sigma_*:
\underbrace{T_{V} \times \dots \times T_{V}}_{d\text{ times}} 
\to T_{V}$; 
it thus suffices to check that this composition map is a bijection.
But considering the diagram
$$\xymatrix{{}&{T_{V}} \\
{T_{(v_1)\oplus \cdots \oplus (v_d)}} \ar[ru]^{\sim} \ar[r] \ar[d]^{\sim} 
& {T_{V \oplus \cdots \oplus V}} \ar[u]^{\sigma_*} \ar[d]^{\sim} \\
{T_{(v_1)}\times \dots \times T_{(v_d)}} \ar[r] & 
{T_{V} \times \dots \times T_{V}},}$$
one checks from the definitions that everything commutes, and inductively
uses the deformation stack axioms to show that the lower vertical arrows
are isomorphisms. Our map
$T_{(v_1)} \times \dots \times T_{(v_d)} \to T_{V}$ is obtained 
from the diagram by inverting the appropriate isomorphism, so we see that 
it is a bijection, as desired.

Finally, to check functoriality, if we have $\phi:V \to W$, which by abuse
of notation we also consider as a map $k[V] \to k[W]$, we know from
Lemma \ref{lem:linear} and the above that all maps in question are
linear, so it suffices to check agreement on elements of the form
$\eta \otimes v$. We therefore need to see that 
$\vp_{\phi(v) *} \eta = \phi_* (\vp_{v *} \eta)$. But we have
$\phi \circ \vp_v = \vp_{\phi(v)}$, so this is trivial.

We next move on to the case of $A_{V}$; the argument is largely the 
same as above. Given a tensor $\vp \otimes v$ with 
$\vp \in A_{\cS}$ and $v \in V$, we
can use the induced map $i_v:k[\epsilon] \to k[V]$ sending $\epsilon$ to
$v$ to obtain $i_{v*} \vp \in \Aut(\zeta_V)$. This map can be extended 
additively to give a map $A_{\cS} \otimes_k V \to \Aut(\zeta_V)$, and it
follows from Lemma \ref{lem:linear} that this map
is well-defined and $k$-linear.
Finally, one checks that this map is an isomorphism by choosing a basis
$v_1,\dots,v_d$ of $V$, and exploiting the description 
$$k[V] \cong k[v_1] \times _k k[v_2] \times_k \cdots \times_k k[v_d],$$
as in the case of $T_V$. The argument for functoriality is likewise the
same.
\end{proof}

Applying the proposition and Theorem \ref{thm:def-aut-ext}, we see
that by writing any $A$ as a series of small extensions of $k$, we
inductively conclude the following.

\begin{cor}\label{cor:triv-aut-def} If $T_{\cS}=0$, then there is at most 
one isomorphism class in $\cS_A$ over any $A$.

If $A_{\cS}=0$, then $\Aut(\eta)=\{\id\}$ for every $\eta \in \cS$.
In particular, in this case (H4) is satisfied for $F_{\cS}$.
\end{cor}

Note that the last assertion follows from Proposition \ref{h4}. 

One can often study $A_{\cS}$ more directly than would be the case
for more general automorphism groups in $\cS$, and in particular one
can sometimes show that $\cS$ has no non-trivial automorphisms even
in cases where automorphisms could {\it a priori} exist. See
Corollary \ref{cor:rigid-schemes} below for an example.

\begin{rem}\label{rem:tan-obs} Proposition \ref{prop:tensors}, together 
with Theorem \ref{thm:def-aut-ext}, also shows that if $\cS$ is a 
deformation stack with an obstruction space, then $F_{\cS}$ has a 
generalized tangent-obstruction theory in the sense of Fantechi and 
G\"ottsche, Definition 6.1.21 of \cite{f-g1}.
\end{rem}

\section{Geometric deformation stacks}\label{sec:gd-stacks}

We now move on to the focus of our main theorems: deformation problems 
which are sufficiently geometric to be naturally associated to a scheme 
$X$. After defining geometric deformation stacks and studying their
most basic properties, we use our results on deformation stacks to
prove Theorems \ref{main-desc} and \ref{main-rep}. We then discuss
examples and special cases of these theorems.

\subsection{Definitions and basic properties}

Many geometric deformation problems arise with more context than just a
deformation functor or stack: specifically, such problems are naturally
associated to some scheme, and given a morphism of schemes, we can pull
back to obtain another deformation problem. We do not pursue this point
of view in such generality: we are interested mainly in being able to
restrict deformation problems to open subsets, so that we can formally
talk about how to describe deformation problems globally in terms of
their local behavior.

We let $X$ be a scheme over a field $k$, and we denote by $\Zar(X)$ 
the category of the Zariski topology on $X$.

\begin{defn} Let $\cS$ be a category fibered in groupoids over
$\Zar(X)\times \Art(\Lambda,k)^{\opp}$. We say that $\cS$ is a 
{\bf geometric deformation stack} or {\bf gd-stack} if it satisfies:
\begin{ilist}
\itm For each $A \in \Art(\Lambda,k)$, $\cS_{(\us,A)}$ is a stack for the 
Zariski topology on $X$;
\itm For each $U \in \Zar(X)$, $\cS_{(U,\us)}$ is a deformation stack. 
\end{ilist}
\end{defn}

Here $\cS_{(\us,A)}$ is the category fibered in groupoids over $\Zar(X)$
obtained as the subcategory of $\cS$ consisting of objects lying over
$(U,A)$ for some $U \in \Zar(X)$, and whose morphisms lie over 
$(\iota,\id_A)$ where $\iota$ is any morphism in $\Zar(X)$. We define
$\cS_{(U,\us)}$ similarly. We assume we have chosen objects 
$\zeta_{0,U} \in \cS_{U,k}$ for each $U$, and 
$\zeta_{V,U} \in \cS_{k[V],U}$ a pushforward of $\zeta_{0,U}$ for every 
finite-dimensional $k$-vector space $V$.

We take a moment to consider the extra structure offered by a category
fibered in groupoids $\cS$ over a product category $C_1 \times C_2$.
Given an object $\eta$ over $(T_1,T_2)$, and a morphism $f_1:T'_1 \to T_1$
in $C_1$, we will often write $f_1^*(\eta)$ as an abbreviation for
$(f_1,\id_{T_2})^*(\eta)$, and similarly for morphisms in $C_2$. We can
define pullbacks of morphisms somewhat more generally than we did initially:

\begin{defn} Let $\eta,\eta'$ be objects of $\cS$ over $(T_1,T_2)$ and
$(T_1',T_2)$ respectively, and let $\vp:\eta' \to \eta$ be a morphism
lying over $(f_1,\id_{T_2})$. Then given $f_2:T_2' \to T_2$, we can
define the {\bf pullback} $f_2^*(\vp):f_2^*(\eta') \to f_2^*(\eta)$ over
$(f_1,\id_{T_2'})$ as the morphism obtained by applying (ii) of the
definition of a category fibered in groupoids to the composed morphism 
$f_2^*(\eta') \to \eta' \overset{\vp}{\to} \eta$ and $f_2^*(\eta) \to \eta$.
We make a similar definition for pullbacks under $f_1$ of morphisms lying 
over $(\id_{T_1},f_2)$.
\end{defn}

The following basic properties of pullback are checked from the definitions,
with the first coming from the fact that 
$(f_1,\id) \circ (\id f_2)=(f_1,f_2)=(\id,f_2) \circ (f_1,\id)$, and the
rest checked directly as in Proposition \ref{prop:cfgs-basic}.

\begin{prop}\label{prop:cfgs-prod-basic} We have
\begin{ilist}
\itm Pullbacks in different variables commute: given $\eta$ in 
$\cS_{(T_1,T_2)}$ and $f_i:T_i' \to T_i$ for $i=1,2$, we have 
$f_1^*(f_2^*(\eta))=f_2^*(f_1^*(\eta))$.
\itm Pullback of morphisms commutes with composition: given 
$T_2'' \overset{f_2'}{\to} T_2' \overset{f_2}{\to} T_2$ in $C_2$,
and $f_1:T_1' \to T_1$ in $C_1$, and $\eta,\eta'$ in $\cS$ over 
$(T_1,T_2),(T_1',T_2)$ respectively, with a morphism $\vp:\eta' \to \eta$
over $(f_1,\id)$, we have $f_2'^*(f_2^*(\vp))=(f_2' \circ f_2)^*(\vp)$.
\itm Pullback of morphisms commutes with composition: given
$f_2:T_2' \to T_2$ in $C_2$, and 
$T_1'' \overset{f_1'}{\to} T_1' \overset{f_1}{\to} T_1$ in $C_1$,
and $\eta,\eta',\eta''$ in $\cS$ over $(T_1,T_2),(T_1',T_2),(T_1'',T_2)$
respectively, with morphisms $\vp: \eta' \to \eta$ and $\vp':\eta'' \to \eta'$
over $(f_1,\id)$ and $(f_1',\id)$ respectively, we have
$f_2^*(\vp \circ \vp')= f_2^*(\vp) \circ f_2^*(\vp')$.
\end{ilist}
\end{prop}

We use the additive structure of \S \ref{subs:additive}
on each open subset $U$ of $X$ to define
an additive structure on $\cS$, which in the case that $A'=k[I]$ gives
a Picard stack in the sense of Deligne (see Definition 1.4.5 of Expose 
XVIII of \cite{sga43}).

\begin{prop}\label{prop:add-struct-stack} Appropriate choices of objects 
$\eta_1 \ast \eta_2$ and $\sigma_* (\eta_1 \ast \eta_2)$ for all
$U \subseteq X$ and
$\eta_1 \in \cS_{U,A'}, \eta_2 \in \cS_{U,k[I]}$ as well as, when
$A' \neq k[I]$, for all $\eta_1,\eta_2 \in \cS_{U,k[I]}$, gives us an 
additive structure consisting of functors
$$+_{A',I}:\cS_{(\us,A')} \times_{\Zar(X)} \cS_{(\us,k[I])} \to 
\cS_{(\us,A')}$$ 
and 
$$+_{k[I],I}:\cS_{(\us,k[I])} \times_{\Zar(X)} \cS_{(\us,k[I])} \to 
\cS_{(\us,k[I])}$$ 
and an isomorphism of functors 
$$\alpha: (\us +_{A',I} \us) +_{A',I} \us \risom
\us +_{A',I} (\us +_{k[I],I} \us)$$
such that:
\begin{ilist}
\itm for any $U\in \Zar(X)$ and $\eta \in \cS_{U,k[I]}$, the functor 
$\us +_{A',I} \eta : \cS_{U,A'} \to \cS_{U,A'}$ is an equivalence of
categories;
\itm $+_{k[I],I}$ is commutative, in the sense that 
$+_{k[I],I} = +_{k[I],I} \circ \sw$, where
$\sw: \cS_{(\us,k[I])} \times_{\Zar(X)} \cS_{(\us,k[I])} \to 
\cS_{(\us,k[I])} \times_{\Zar(X)} \cS_{(\us,k[I])}$ 
is the functor switching factors;
\itm for every $U \in \Zar(X)$, we have that
$\us +_{A',I} \zeta_{I,U}: \cS_{U,A'} \to \cS_{U,A'}$ is the identity 
functor;
\itm for any $U \in \Zar(X)$ and 
$\eta_1,\eta_2,\eta_3,\eta_4 \in \cS_{U,k[I]}$, we have the identities
$$(\id_{\eta_1}+\alpha_{\eta_2, \eta_3,\eta_4}) \circ 
\alpha_{\eta_1,\eta_2+\eta_3,\eta_4} \circ
(\alpha_{\eta_1,\eta_2,\eta_3}+\id_{\eta_4}) =
\alpha_{\eta_1,\eta_2,\eta_3+\eta_4} \circ
\alpha_{\eta_1+\eta_2,\eta_3,\eta_4}$$
and
$$\alpha_{\eta_1,\eta_2,\eta_3} \circ \alpha_{\eta_3,\eta_1,\eta_2} 
= \alpha_{\eta_1,\eta_3,\eta_2}.$$
\end{ilist}

In the above and hereafter, we denote by
$\alpha_{\eta_1,\eta_2,\eta_3}: (\eta_1+\eta_2)+\eta_3 
\risom \eta_1 +(\eta_2 + \eta_3)$ the isomorphism obtained from $\alpha$.

Such an addition law is canonical up to unique isomorphism of functors.
\end{prop}

Here, the category $\cS_{(\us,A')} \times_{\Zar(X)} \cS_{(\us,k[I])}$ 
has objects consisting of pairs of objects in $\cS_{(\us,A')}$ over a 
given $U \in \Zar(X)$, and morphisms pairs of morphisms lying over a given
inclusion in $\Zar(X)$.

\begin{proof} Indeed, we claim that we obtain the desired addition law
from any choice of addition law as in Corollary \ref{c:add-struct}, 
chosen independently for each $U \subseteq X$. We need to see that
this defines the desired functor
$+_{A',I}:\cS_{(\us,A')} \times \cS_{(\us,k[I])} \to \cS_{(\us,A')}$. 
We clearly
obtain a functor on the level of objects, and also for morphisms over a 
fixed $U \in \Zar(X)$, simply by restricting to $U$. 

What remains to define is addition of arbitrary pairs of morphisms: given
$U' \subseteq U$, $\eta_1 \in \cS_{U',A'}$, $\eta_2 \in \cS_{U',k[I]}$,
$\mu_1 \in \cS_{U,A'}$ and $\mu_2 \in \cS_{U,k[I]}$, together
with $\vp_i:\eta_i \to \mu_i$ over $U' \to U$, we need to define
$\vp_1 +_{A',I} \vp_2$ so that $+_{A',I}$ is a functor. We have already
chosen $\eta_1 \ast \eta_2$ and $\mu_1 \ast \mu_2$, and their pushforwards
$\eta_1 + \eta_2$ and $\mu_1 + \mu_2$, and we claim there exists a unique
$\vp_1 \ast \vp_2:\eta_1 \ast \eta_2 \to \mu_1 \ast \mu_2$ such that
the restrictions to $A'$ and $k[I]$ recover $\vp_1$ and $\vp_2$ 
respectively. We first show that there exists some map
$\vp:\eta_1 \ast \eta_2 \to \mu_1 \ast \mu_2$ lying over
$U' \to U$. Indeed, if we make any choice of pullback
$(\mu_1 \ast \mu_2)|_{U'} \to \mu_1 \ast \mu_2$,
we can use the $\vp_i$ as identifications $\eta_i=\mu_i|_{U'}$, and
we therefore obtain maps $(\mu_1 \ast \mu_2)|_{U'} \to \eta_i$
over $A' \times_k k[I] \to A'$ and $A' \times_k k[I] \to k[I]$.
This realizes $(\mu_1 \ast \mu_2)|_{U'}$ as $\eta_1 \ast \eta_2$,
so we obtain a (unique) isomorphism with our chosen $\eta_1 \ast \eta_2$,
and composing the isomorphism with the pullback map gives the desired
map $\vp$. 

We now show that we have a unique $\vp_1 \ast \vp_2$. The 
definition of a category fibered in groupoids implies that morphisms 
$\eta_1 \ast \eta_2 \to \mu_1 \ast \mu_2$ over $U' \to U$ are in 
bijection with $\Aut_{U'}(\eta_1 \ast \eta_2)$,
and similarly for maps $\eta_i \to \mu_i$ and $\Aut_{U'}(\eta_i)$.
Thus, if $\vp|_{A'}$ differs from $\vp_1$ by $\psi_1 \in \Aut_{U'}(\eta_1)$,
and $\vp|_{k[I]}$ differs from $\vp_2$ by $\psi_2 \in \Aut_{U'}(\eta_2)$,
by the deformation stack axioms there is a unique 
$\psi \in \Aut_{U'}(\eta_1 \ast \eta_2)$ restricting to $\psi_1$ and 
$\psi_2$, and setting $\vp_1 \ast \vp_2 = \vp \circ \psi$ yields the
desired unique map. We then set $\vp_1 + \vp_2$ to be the restriction 
under $\sigma$ of $\vp_1 \ast \vp_2$, given the choices of $\eta_1+\eta_2$
and $\mu_1 + \mu_2$. The uniqueness of this construction immediately
implies that it commutes with composition, so we obtain the desired
functor $+_{A',I}$.

We then obtain the isomorphism $\alpha$ from Corollary \ref{c:add-struct},
since any given $\alpha_{\eta_1,\eta_2,\eta_3}$ occurs over a single $U$.
For the same reason, properties (i), (iii), and (iv) follow directly from 
Corollary \ref{c:add-struct}. Finally, (ii) is automatically satisfied
on the level of objects, and one checks that it also works for morphisms
exactly as before. Since our choices are made exactly as in Corollary
\ref{c:add-struct}, the assertion that two such choices differ by a
unique isomorphism also follows from Corollary \ref{c:add-struct}.
\end{proof}

We obtain the following corollary on additivity of pullbacks:

\begin{cor}\label{cor:pullbacks-add} Given $U' \subseteq U$ in $\Zar(X)$,
and objects $\eta_1 \in \cS_{U,A'}$, $\eta_2 \in \cS_{U,k[I]}$, and choices 
of pullbacks
$\eta_1|_{U'} \to \eta_1$, $\eta_2|_{U'} \to \eta_2$, the 
$+_{A',I}$ functor gives a
canonical pullback map $\eta_1|_{U'}+\eta_2|_{U'} \to \eta_1 +\eta_2$.

Given also $\mu_1 \in \cS_{U,A'}$, $\mu_2 \in \cS_{U,k[I]}$, and 
choices of pullbacks $\mu_1|_{U'} \to \mu_1$, $\mu_2|_{U'} \to \mu_2$, if
we use the induced pullback map 
$\mu_1|_{U'}+\mu_2|_{U'} \to \mu_1 +\mu_2$, then for any 
$f_1:\eta_1 \to \mu_1$ in $\cS_{U,A'}$ and $f_2:\eta_2 \to \mu_2$ in
$\cS_{U,k[I]}$, we have the indentity
$$(f_1+f_2)|_{U'}=f_1|_{U'}+f_2|_{U'}.$$
\end{cor}

\begin{proof}
It is clear from the definitions that $+_{A',I}$ gives a canonical
pullback map $\eta_1|_{U'}+\eta_2|_{U'} \to \eta_1 +\eta_2$, simply
as the sum of $\eta_1|_{U'} \to \eta_1$ and $\eta_2|_{U'} \to \eta_2$.
For addition of morphisms, the commutative diagrams defining $f_1|_{U'}$
and $f_2|_{U'}$ yield a commutative diagram in 
$\cS_{(\us,A')} \times_{\Zar(X)} \cS_{(\us,k[I])}$,
and the fact that $+_{A',I}$ is a functor then shows that
$f_1|_{U'}+f_2|_{U'}$ satisfies the commutativity condition defining
$(f_1+f_2)|_{U'}$.
\end{proof}

Our observations that the functor associated to a deformation stack is
always a deformation functor, and that (H4) can be understood concretely
in terms of automorphisms, reduces representability questions down to
(H3), which is to say, to understanding tangent spaces. Thus, we will
be mainly interested in describing tangent and obstruction spaces to a 
gd-stack. 

\begin{defn} Given a gd-stack $\cS$, we define the {\bf tangent space} 
$T_{\cS}$ of $\cS$, a {\bf successive obstruction theory} for $\cS$ taking 
values in $V_1,\dots,V_n$, and an {\bf obstruction space} for $\cS$ 
to be the tangent space of $\cS_{(X,\us)}$, a successive obstruction theory
for $\cS_{(X,\us)}$, and an obstruction space for $\cS_{(X,\us)}$
respectively.
\end{defn}

\begin{lem}\label{lem:restrict-linear} Let $\cS$ be a gd-stack, 
$U_1 \subseteq U_2 \in \Zar(X)$, and $V$ a $k$-vector space. 
The maps $T_{\cS_{(U_2,\us)},V} \to T_{\cS_{(U_1,\us)},V}$ and
$A_{\cS_{(U_2,\us)},V} \to A_{\cS_{(U_1,\us)},V}$ induced by pullback 
from $U_2$ to $U_1$ are linear maps of $k$-vector spaces. 
\end{lem}

\begin{proof} That the maps are additive follows from Corollary
\ref{cor:pullbacks-add}. They also commute with scalar multiplication
because it is defined in terms of pushforwards along ring homomorphisms, 
and since we are simply working over the product category 
$\Zar(X) \times \Art(\Lambda,k)^{\opp}$, such pushforwards commute with
the restriction from $U_2$ to $U_1$.
\end{proof}

As mentioned earlier, our main results on tangent and obstruction spaces
will take the form of describing the spaces in terms of local information.
The main tools which will allow this are the following sheaves of
$k$-vector spaces on $X$:

\begin{defn} Given a gd-stack $\cS$, denote by $\cT_{\cS}$ and $\cA_{\cS}$
respectively the {\bf tangent sheaf} of $\cS$, given as the sheafification 
of the presheaf $U \mapsto T_{\cS_{(U,\us)}}$, and the {\bf infinitesimal 
automorphism sheaf} $U \mapsto A_{\cS_{(U,\us)}}$. By Lemma 
\ref{lem:restrict-linear}, these are both sheaves of $k$-vector spaces.

Similarly, for each $k$-vector space $V$, let $\cT_{\cS,V}$ and $\cA_{\cS,V}$
be the sheaves of $k$-vector spaces given by the sheafication of
$U \mapsto T_{\cS_{(U,\us)},V}$, and by $U \mapsto A_{\cS_{(U,\us)},V}$ 
respectively.
\end{defn}

In order to study lifts of objects over small extensions, it will also
be convenient to introduce the following:

\begin{notn} Given a small extension $A' \to A$ with kernel $I$, and 
$\eta \in \cS_{(X,A)}$, fix choices of restrictions $\eta|_U$ for all 
$U \in \Zar(X)$, and denote by $\cS_{\eta,A'}$ the stack over $\Zar(X)$
whose objects consist of triples $(U,\eta',\vp)$, with $U \in \Zar(X)$,
$\eta' \in \cS_{U,A'}$, and $\vp:\eta|_U \to \eta'$ lying over $A' \to A$. 
Morphisms $\psi:(U_1,\eta'_1,\vp_1) \to (U_2,\eta'_2,\vp_2)$ consist
of a morphism $\eta'_1 \to \eta'_2$ in $\cS_{(\us,A')}$ which lies over an
inclusion $U_1 \to U_2$ and which satisfies
$\psi \circ \vp_1 = \vp_2 \circ r$, where $r:\eta|_{U_1} \to \eta|_{U_2}$
is the canonical map.

We also let $T_{\eta,A',U}$ be the set of isomorphism classes of
$(\cS_{\eta,A'})_U$, and $\cT_{\eta,A'}$ be the sheaf of 
$\cT_{I}$-pseudotorsors obtained as the sheafification of the 
presheaf given by $U \mapsto T_{\eta,A',U}$.
\end{notn}

Finally, for the statements of our main theorems in full generality, it
will be helpful to have the notion of a local obstruction sheaf also:

\begin{defn} We say that a sheaf of $k$-vector spaces $\cOb$ on $X$ is
a {\bf local obstruction sheaf} if we have also the data, for each
$U \in \Zar(X)$, each small extension $A' \twoheadrightarrow A$
in $\Art(\Lambda,k)$ with kernel $I$, and each $\eta \in \cS_{(U,A)}$, of 
an element $\ob_{\eta,A'} \in \cOb(U)\otimes _k I$, compatible with 
restriction in $\Zar(X)$ and satisfying the functoriality condition of an 
obstruction space given in Definition \ref{d:tangent-obs}, and such that 
$\ob_{\eta,A'}=0$ if and only if there exists an open cover of $U$ by 
$\{U_i\}$, and $\eta'_i \in \cS_{(U_i,A')}$ such that 
$\eta'_i|_A \cong \eta|_{U_i}$ for all $i$.
\end{defn}

\subsection{Proof of the main theorems}

We now move on to the proofs of our main theorems. Given our results on
deformation stacks, the proofs follow the usual proofs used in standard
examples of deformation problems, and are fundamentally quite intuitive: 
we pass from local to global using open covers and gluing along 
intersections. The main construction for both the tangent space exact 
sequence and the final obstruction space in Theorem \ref{main-desc} is 
the following.

\begin{prop}\label{prop:very-short} Fix a small extension $A' \to A$,
with kernel $I$, and $\eta \in \cS_{(X,A)}$. Then the
map of sets induced by sheafification of $U \mapsto T_{\eta,A',U}$ can be 
extended to the following very short exact sequence:
$$T_{\eta,A',X} \to \Gamma(X,\cT_{\eta,A'}) \to H^2(X, \cA_I),$$
meaning that we have maps of sets, with the kernel of the second
map equal to the image of the first. This sequence is functorial in that 
for any small extension $B' \to B$ with kernel $J$, and map 
$\phi:A' \to B'$ mapping $I$ into $J$, the above maps commute with the 
maps induced by $\phi$, restriction of $\eta$ to $B$, and $I \to J$.

Furthermore, the constructed map
$\Gamma(X,\cT_{\eta,A'}) \to H^2(X, \cA_I)$ is compatible with addition,
and in the case that $A'=k[I]$ and $A=k$, our construction yields a
$k$-linear map $\Gamma(X,\cT_{I}) \to H^2(X, \cA_I)$.
\end{prop}

Before giving a more elementary proof along the lines of standard 
\v{C}ech cohomology arguments, we mention a stack-theoretic
proof of the proposition: given a section $s \in H^0(X,\cT_{\eta,A'})$,
one can define a substack $\cS_s$ of $\cS_{\eta,A'}$ consisting
of the full subcategory containing every object locally isomorphic
to the objects defining $s$. This is then an $\cA_I$-gerbe over $X$,
whose class is an element of $H^2(X, \cA_I)$, giving us the desired
map. Exactness is a consequence of the fact that the class of a
gerbe is trivial if and only if the gerbe has a global section, which
is to say if and only if $s$ actually came from an element of $T_{\eta,A'}$.

Our more elementary proof follows the ideas of classical \v{C}ech 
cohomology, but cannot rely solely on them, and the 
reason is quite simple: given a presheaf $\cF$ on a topological space $X$ 
with sheafification $\tilde{\cF}$, and a global section 
$\rho \in \tilde{\cF}(X)$, there is not any reason to think that
we can find an open cover $\{U_i\}$ of $X$ on which $\rho$ is the
image of some $\{\rho_i \in \cF(U_i)\}$ which agree on the $U_i \cap U_j$.
Certainly, $\rho$ is the image of some $\rho_i$, but we see that 
{\it a priori} we have equality of $\rho_i$ and $\rho_j$ only locally
on $U_i \cap U_j$, and it is not necessarily possible to refine the $U_i$ 
to obtain the desired equality on all of $U_i \cap U_j$. See Example
\ref{ex:hyper-cech} below for one such presheaf. We therefore work with 
the simplest generalization of \v{C}ech cohomology in the direction of 
hypercovers, which is enough to describe the sheaf cohomology group $H^2$ 
on an arbitrary topological space with an arbitrary sheaf of groups. We 
give a self-contained account of this approach to cohomology in Appendix
\ref{app:hyper-cech} below.

\begin{proof} We construct the map 
$\Gamma(X,\cT_{\eta,A'}) \to H^2(X, \cA_I)$
as follows: given $\rho \in \Gamma(X,\cT_{\eta,A'})$, let $\{U_i\}$ be an 
open cover small enough so that $\rho$ is represented by sections 
$(\rho_i,\vp_i) \in T_{\eta,A',U_i}$, with $\rho_i \in \cS_{(U_i,A')}$ and
$\vp_i:\eta|_{U_i}\to \rho_i$ lying over $A' \to A$. For each 
$i_0 \neq i_1$, let 
$\{U_{i_0,i_1,j}\}_{j}$ be an open cover of $U_{i_0,i_1}$ such that 
on each $U_{i_0,i_1,j}$ there exists an isomorphism 
$\vp_{i_0,i_1,j}:
\rho_{i_0}|_{U_{i_0,i_1,j}} \to \rho_{i_1}|_{U_{i_0,i_1,j}}$ commuting
with $\vp_{i_0}|_{U_{i_0,i_1,j}}$ and $\vp_{i_1}|_{U_{i_0,i_1,j}}$.
We choose such $\vp_{i_0,i_1,j}$, imposing that 
$\vp_{i_0,i_1,j}=\vp_{i_1,i_0,j}^{-1}$, and we also set the convention
that for $i_0=i_1$, we take the one-set cover $U_{i_0}$, with 
$\vp_{i_0,i_0,\_}=\id$.
We then obtain a $2$-cochain $\rho'$ with coefficients in $\cA_I$, in the 
sense of Appendix \ref{app:hyper-cech}, uniquely characterized via
Theorem \ref{thm:def-aut-ext} by the identity, for 
$\bi=(i_0,i_1,i_2)$ and $\bj=(j_{1,2},j_{0,2},j_{0,1})$,
$$\id_{\rho_{i_0}}+\rho'_{\bi,\bj}= 
\vp_{i_2,i_0,j_{0,2}}|_{U_{\bi,\bj}} \circ
\vp_{i_1,i_2,j_{1,2}}|_{U_{\bi,\bj}} \circ
\vp_{i_0,i_1,j_{0,1}}|_{U_{\bi,\bj}}.$$ 
Here and subsequently, to minimize notational clutter we generally leave 
as implicit restrictions of identity maps to the appropriate open subsets.

We now check that the $\rho'$ we have constructed is a cocycle. Given
$\bi=(i_0,i_1,i_2,i_3)$ and $\bj=(j_{m,n})_{0 \leq m < n \leq 3}$, denote 
by $\bj_{\ell}$ the triple of $j_{m,n}$ with $m,n \neq \ell$. To check
that $(d \rho')_{\bi,\bj}=0$, we want to see that
$$\rho'_{i_0,i_1,i_2,\bj_3}|_{U_{\bi,\bj}}-
\rho'_{i_0,i_1,i_3,\bj_2}|_{U_{\bi,\bj}}+
\rho'_{i_0,i_2,i_3,\bj_1}|_{U_{\bi,\bj}} - 
\rho'_{i_1,i_2,i_3,\bj_0}|_{U_{\bi,\bj}}=0,$$ 
or equivalently,
$$\rho'_{i_0,i_2,i_3,\bj_1}|_{U_{\bi,\bj}} \circ
\rho'_{i_0,i_1,i_2,\bj_3}|_{U_{\bi,\bj}} =
\rho'_{i_0,i_1,i_3,\bj_2}|_{U_{\bi,\bj}} \circ
\rho'_{i_1,i_2,i_3,\bj_0}|_{U_{\bi,\bj}}.$$
By Corollary \ref{cor:mors-group} (iii) this is equivalent to checking
that 
$$\id_{\rho_{i_0}}+(\rho'_{i_0,i_1,i_2,\bj_3}|_{U_{\bi,\bj}} 
\circ \rho'_{i_0,i_2,i_3,\bj_1}|_{U_{\bi,\bj}}) =
\id_{\rho_{i_0}}+(\rho'_{i_0,i_1,i_3,\bj_2}|_{U_{\bi,\bj}}
\circ \rho'_{i_1,i_2,i_3,\bj_0}|_{U_{\bi,\bj}}).$$
But using Corollary \ref{cor:pullbacks-add} and Corollary 
\ref{cor:mors-group} (iii) we have
$$\id_{\rho_{i_0}}+(\rho'_{i_0,i_1,i_2,\bj_3}|_{U_{\bi,\bj}}
\circ \rho'_{i_0,i_2,i_3,\bj_1}|_{U_{\bi,\bj}}) = 
(\id_{\rho_{i_0}} + \rho'_{i_0,i_1,i_2,\bj_3})|_{U_{\bi,\bj}} \circ
(\id_{\rho_{i_0}} + \rho'_{i_0,i_2,i_3,\bj_1})|_{U_{\bi,\bj}}$$
and using also Corollary \ref{cor:mors-group} (iv) we find
\begin{multline*} \id_{\rho_{i_0}} +
(\rho'_{i_0,i_1,i_3,\bj_2}|_{U_{\bi,\bj}}\circ
\rho'_{i_1,i_2,i_3,\bj_0}|_{U_{\bi,\bj}}) \\
= (\id_{\rho_{i_0}} + \rho'_{i_0,i_1,i_3,\bj_2}|_{U_{\bi,\bj}}) \circ
(\id_{\rho_{i_0}} + \rho'_{i_1,i_2,i_3,\bj_0}|_{U_{\bi,\bj}}) \\
= (\id_{\rho_{i_0}} + \rho'_{i_0,i_1,i_3,\bj_2})|_{U_{\bi,\bj}} 
\circ (\vp^{-1}_{i_0,i_1,j_{0,1}}|_{U_{\bi,\bj}} \circ 
(\id_{\rho_{i_1}} + \rho'_{i_1,i_2,i_3,\bj_0})|_{U_{\bi,\bj}} 
\circ \vp_{i_0,i_1,j_{0,1}}|_{U_{\bi,\bj}}),
\end{multline*}
and finally expanding both sides in terms of the various $\vp_{\bi,\bj}$ 
we obtain the desired identity.

A similar argument shows that our constructed map to $H^2(X, \cA_I)$ is 
well-defined: first, if we modify our $\vp_{i_0,i_1,j}$ by automorphisms 
$\gamma_{i_0,i_1,j}$ of the $\rho_{i_0}|_{U_{i_0,i_1,j}}$, we change the 
resulting cocycle by a coboundary. Indeed, if 
$\delta_{i_0,i_1,j} \in \cA_I(U_{i_0,i_1,j})$ is determined by
$\id_{\rho_{i_0}}+\delta_{i_0,i_1,j}
=\gamma_{i_0,i_1,j}$ for all $i_0,i_1,j$, and if we fix triples
$\bi=(i_0,i_1,i_2)$, $\bj=(j_{0,1},j_{0,2},j_{1,2})$, we assert that
precomposing the $\vp_{i_0,i_1,j}$ by $\gamma_{i_0,i_1,j}$ has the effect 
of modifying $\rho'_{\bi,\bj}$ by 
$\delta_{i_0,i_1,j_{0,1}}|_{U_{\bi,\bj}}
-\delta_{i_0,i_2,j_{0,2}}|_{U_{\bi,\bj}}
+\delta_{i_1,i_2,j_{1,2}}|_{U_{\bi,\bj}}$.
To check this, using Corollary \ref{cor:mors-group} (iv) and
Corollary \ref{cor:pullbacks-add} we see that 
$$\id_{\rho_{i_0}}+
\delta_{i_1,i_2,j_{1,2}}|_{U_{\bi,\bj}}
=\vp_{i_0,i_1,j_{0,1}}^{-1}|_{U_{\bi,\bj}} \circ
\gamma_{i_1,i_2,j_{1,2}}|_{U_{\bi,\bj}} \circ 
\vp_{i_0,i_1,j_{0,1}}|_{U_{\bi,\bj}}$$ 
and
$$\id_{\rho_{i_0}}+
(-\delta_{i_0,i_2,j_{0,2}})|_{U_{\bi,\bj}}
=\vp_{i_2,i_0,j_{0,2}}|_{U_{\bi,\bj}} \circ
\gamma_{i_2,i_0,j_{0,2}}|_{U_{\bi,\bj}} \circ 
\vp_{i_2,i_0,j_{0,2}}^{-1}|_{U_{\bi,\bj}},$$ 
and we then obtain the required identity by expanding 
\begin{multline*}
\id_{\rho_{i_0}}+
(\rho'_{\bi,\bj}+\delta_{i_0,i_1,j_{0,1}}|_{U_{\bi,\bj}}
-\delta_{i_0,i_2,j_{0,2}}|_{U_{\bi,\bj}}
+\delta_{i_1,i_2,j_{1,2}}|_{U_{\bi,\bj}}) \\
= \id_{\rho_{i_0}}+
((-\delta_{i_0,i_2,j_{0,2}}|_{U_{\bi,\bj}})\circ
\rho'_{\bi,\bj} \circ \delta_{i_1,i_2,j_{1,2}}|_{U_{\bi,\bj}}
\circ \delta_{i_0,i_1,j_{0,1}}|_{U_{\bi,\bj}}).
\end{multline*}
Next, if we replace our choices of $(\rho_{i},\vp_{i})$ by different 
but isomorphic choices, we simply fix isomorphisms $\alpha_{i}$ on each 
$U_{i}$, and replace each $\vp_{i_0,i_1,j}$ by 
$\alpha_{i_1}^{-1}|_{U_{i_0,i_1,j}} \circ \vp_{i_0,i_1,j} 
\circ \alpha_{i_0}|_{U_{i_0,i_1,j}}$,
and using Corollary \ref{cor:mors-group} once more we verify that
$\rho'$ remains unchanged. Finally, given any two choices
of the $(\rho_{i},\vp_i)$ on any two covers $\{U_i\}$, if we take a common
refinement and refine further if necessary so that our choices are
isomorphic on each $\{U_i\}$, we see that our choices differ by a
coboundary, so our map is well defined.

It then follows that we have the desired functoriality, as well: if we carry
out the construction for $\rho$, using Propositions 
\ref{prop:cfgs-prod-basic} (i) and \ref{prop:cfgs-basic} (iv) and the 
functoriality of addition of automorphisms in Theorem \ref{thm:def-aut-ext} 
we see that pushforward to $B'$ commutes with each step of the 
construction, so that the image of $\rho|_{B'}$ is $\rho'|_{J}$, as desired.

We next check exactness at $\Gamma(X,\cT_{\eta,A'})$: the kernel of the 
constructed map consists precisely of the $\rho$ for which, after possible
refinement of the cover, there exist a collection of
$\rho''_{i_0,i_1,j} \in \cA_I(U_{i_0,i_1,j})$ such that $\rho'=d\rho''$.
Given any collection of $\rho''$, one checks that if we set
$\vp'_{i_0,i_1,j}=
\vp_{i_0,i_1,j}\circ (\id_{\rho_{i_0}}- \rho''_{i_0,i_1,j}),$
then $\rho'=d\rho''$ is equivalent to the $\vp'$ satisfying the cocycle
condition, and the existence of such $\vp'$ is in turn equivalent to 
being able to glue the $\rho_{i_0}$ compatibly to obtain an element of 
$T_{\eta,A'}$ inducing $\rho$. Since the $\rho''$ can be recovered from
the $\vp'$, this proves exactness.

Finally, we check the additivity and linearity assertions.
For additivity, we use the additivity of pullbacks given in 
Corollary \ref{cor:pullbacks-add}: suppose we have 
$\rho^1\in \Gamma(X,\cT_{\eta,A'})$ and $\rho^2\in \Gamma(X,\cT_I)$, which 
we have chosen to be represented by $\{\rho^1_i\}$ and $\{\rho^2_i\}$ 
respectively on a common cover $\{U_i\}$, with isomorphisms
$\{\vp_{i_0,i_1,j}^1\}$ and $\{\vp_{i_0,i_1,j}^2\}$ on common refinements
$U_{i_0,i_1,j}$. The resulting cocycles ${\rho^1}'$ and ${\rho^2}'$ are
characterized by the condition that on each $U_{\bi,\bj}$,
for $m=1,2$ we have
$${\rho^m}'_{i_0,i_1,i_2,\bj}+\id_{\rho_{i_0}^m}=
\vp_{i_2,i_0,j_{0,2}}^m \circ \vp_{i_1,i_2,j_{1,2}}^m 
\circ \vp_{i_0,i_1,j_{0,1}}^m.$$
We can then represent $\rho^1+\rho^2$ by $\{\rho^1_i+\rho^2_i\}$ on
$\{U_i\}$, and choose our isomorphisms on $U_{i_0,i_1,j}$ to be given
by $\vp_{i_0,i_1,j}^1+\vp_{i_0,i_1,j}^2$. The desired assertion then
comes down to checking that 
\begin{multline*}
(\id_{\rho_{i_0}^1|_{U_{\bi,\bj}}} + {\rho^1}'_{i_0,i_1,i_2,\bj}) +
(\id_{\rho_{i_0}^2|_{U_{\bi,\bj}}} + {\rho^2}'_{i_0,i_1,i_2,\bj}) \\
= \id_{\rho_{i_0}^1|_{U_{\bi,\bj}} + \rho_{i_0}^2|_{U_{\bi,\bj}}}+
({\rho^1}'_{i_0,i_1,i_2,\bj} + {\rho^2}'_{i_0,i_1,i_2,\bj}).
\end{multline*}
But this is Corollary \ref{cor:mors-group} (vi). 

Lastly, in the case that $A'=k[I]$, our map is also $k$-linear: one need
only verify that it commutes with scalar multiplication, which is
easily checked directly: we may represent $m_{\lambda *} \rho$
by $m_{\lambda *} \rho_i$ on each $U_i$, with isomorphisms
$m_{\lambda *} \vp_{i_0,i_1,j}$, and we get the desired identity
from the fact that $m_{\lambda *}$ commutes with addition of
automorphisms and composition of morphisms.
\end{proof}

We are now ready to prove Theorem \ref{main-desc}, whose statement we
now recall:

\begin{thm*} Let $\cS$ be a gd-stack. 
Then the tangent space $T_{\cS}$ of $\cS$ fits into an exact sequence of 
$k$-vector spaces
$$0 \to H^1(X,\cA) \to T_{\cS} \to H^0(X,\cT) \to H^2(X,\cA),$$
and if we are given a local obstruction sheaf $\cOb$ for $\cS$, 
we have successive obstructions lying in 
$H^0(X,\cOb)$, $H^1(X,\cT)$, and $H^2(X,\cA)/H^0(X,\cT)$.
\end{thm*}

\begin{proof}[Proof of Theorem \ref{main-desc}] We first address the 
exact sequence for $T_{\cS}$. The second 
map is obtained from the sheafification map, and is therefore automatically
$k$-linear, with
its kernel consisting precisely of deformations which are locally trivial.
The first map is constructed by considering such locally trivial
deformations: given a \v{C}ech 1-cocycle $\rho$ of $\cA$ in terms of an 
open cover $\{U_i\}$, on each $U_i$ we can take the trivial deformation, 
and by the stack condition for $\cS$ over $\Zar(X)$, we can use $\rho$ to
glue the deformations to obtain an $\eta_{\rho} \in T_{\cS}$, which
moreover will be trivial if and only if $\rho$ can be simultaneously
trivialized by automorphisms of the $U_i$; i.e., if and only if $\rho$ is 
a coboundary. This gives the first inclusion, and it is clear
that the image is precisely the set of locally trivial deformations, so
we obtain exactness at $T_{\cS}$. One also checks easily that the inclusion 
map is $k$-linear using the properties given in Proposition 
\ref{prop:cfgs-prod-basic} and the fact from Proposition 
\ref{prop:add-struct-stack} and Corollary \ref{cor:pullbacks-add} that 
$+_{k[\epsilon],(\epsilon)}$ is a functor, and we can use it to add 
pullback maps.

The last map is obtained from Proposition \ref{prop:very-short} in the
case $V=(\epsilon)$, completing the desired exact sequence. 

We now consider obstructions. We therefore fix a small extension
$A' \to A$ in $\Art(\Lambda,k)$ with kernel $I$, and an object 
$\eta \in \cS_{(X,A)}$. By
definition, we get an element $\ob^1_{\eta,A'} \in H^0(X,\cOb)\otimes I$, 
which is zero if and only if $\eta$ can be extended to $A'$ locally. If 
$\ob^1_{\eta,A'} \neq 0$, we set this as our obstruction. 

If it is zero,
there exists some cover $\{U_i\}$ of $X$ and liftings $\eta'_i$ of 
$\eta|_{U_i}$ from $A$ to $A'$. The next question is whether these liftings
can be chosen so that they are isomorphic on the $U_{i_0,i_1}$. By 
Theorem \ref{thm:def-aut-ext} and Proposition \ref{prop:tensors}, the 
choices of the $\eta'_i$ together with given pullback maps 
$\eta|_{U_i} \to \eta'_i$ over $A' \to A$ (up to simultaneous 
isomorphism) are torsors under $T_{\cS_{(U_i,\us)}} \otimes I$, so if we 
restrict to the $U_{i_0,i_1}$, taking differences of $\eta'_{i_0}$ and 
$\eta'_{i_1}$ give elements 
$\eta'_{i_0,i_1} \in T_{\cS_{(U_{i_0,i_1},\us)}} \otimes I$. Moreover, 
by Lemma \ref{lem:restrict-linear} and general properties of torsors
we see that given $i_0,i_1,i_2$ we have the cocycle 
condition 
$$\eta'_{i_0,i_1}|_{U_{i_0,i_1,i_2}}+\eta'_{i_1,i_2}|_{U_{i_0,i_1,i_2}}
+\eta'_{i_2,i_0}|_{U_{i_0,i_1,i_2}}=0,$$
so we obtain an element $\ob^2_{\eta,A'} \in H^1(X,\cT) \otimes I$.
Furthermore, this element vanishes if and only if 
the choices of lifts $\eta'_i$ can be modified simultaneously so that they
are all isomorphic on every $U_{i_0,i_1}$, giving us an element of
$\Gamma(X,\cT_{\eta,A'})$.
Thus, if $\ob^2_{\eta,A'} \neq 0$, we use it as our obstruction.
We observe that this is independent of choices: modifying our choices
of the $\eta'_i$ changes the constructed cocycle by a coboundary,
while different choices of the cover $\{U_i\}$ may be compared by
restriction to a common refinement. In particular, $\ob^2_{\eta,A'}=0$
if and only if $\Gamma(X,\cT_{\eta,A'})$ is non-empty.

Finally, if $\ob^2_{\eta,A'}=0$, by choosing lifts $\eta'_i$ which agree on
$U_{i_0,i_1}$ we obtain an element of $\Gamma(X,\cT_{\eta,A'})$, and by 
Proposition \ref{prop:very-short} 
above and using the isomorphism $\cA_I \risom \cA \otimes_k I$ of
Proposition \ref{prop:tensors}, we obtain an element of
$H^2(X,\cA) \otimes I$. We then set $\ob^3_{\eta,A'}$ to be its image in 
$$(H^2(X, \cA)/H^0(X,\cT)) \otimes I=
(H^2(X,\cA) \otimes I)/(H^0(X,\cT) \otimes I).$$ 
Moreover, because $\Gamma(X,\cT_{\eta,A'})$ is a torsor over 
$H^0(X,\cT_I)=H^0(X,\cT)\otimes I$, it follows from the additivity
in Proposition \ref{prop:very-short} that 
$\ob^3_{\eta,A'}$ is independent of the choice of element of 
$\Gamma(X,\cT_{\eta,A'})$. 
It thus follows from the exactness in Proposition \ref{prop:very-short}
that if $T_{\eta,A'}$ is non-empty, we have $\ob^3_{\eta,A'}=0$. 
Conversely, if $\ob^3_{\eta,A'}=0$, we have an element of
$\mu \in \Gamma(X,\cT_{\eta,A'})$ whose image inside $H^2(X,\cA_I)$
agrees with the image of some 
$\rho \in H^0(X,\cT) \otimes I = H^0(X,\cT_I)$. 
Thus, by the additivity and exactness in Proposition \ref{prop:very-short} 
we see that $\rho-\mu$ maps to $0$ in $H^2(X,\cA_I)$ and hence is the
image of an element of $T_{\eta,A'}$, meaning that an $\eta'$ lifting
$\eta$ exists, as desired.

It remains only to check that the obstruction theory we have constructed
is functorial, and this is straightforward: for $\ob^1$ this is part of
the definition of a local obstruction sheaf; for $\ob^2$ it follows from 
the functoriality in Theorem \ref{thm:def-aut-ext} and Propositions 
\ref{prop:cfgs-prod-basic} (i) and \ref{prop:tensors}; and for $\ob^3$ we use 
the functoriality in Propositions \ref{prop:very-short} and 
\ref{prop:tensors}.
\end{proof}

\begin{rem} We observe that in fact in the proof of Theorem \ref{main-desc}
we did not use the existence of a local obstruction sheaf, but merely
of a single vector space which measures the local obstruction to extending
a given global object. However, whenever one has such a space one expects
to obtain such spaces under restriction to every $U \subseteq X$, and
therefore to obtain a local obstruction sheaf, so the sheaf terminology
seems more natural.
\end{rem}

We are now able to say quite a bit in the context of Schlessinger's theory 
of representability and hulls, as well. Indeed, recall the statement of 
Theorem \ref{main-rep}:

\begin{thm*} Let $\cS$ be a gd-stack on a scheme $X$. Then:

\begin{ilist}
\itm the associated functor $F_{\cS_{(X,\us)}}$ satisfies Schlessinger's 
(H1) and (H2), and satisfies (H4) if and only if for each tiny extension
$A' \to A$ in $\Art(\Lambda,k)$, and each object $\eta \in A'$, the natural
map
$$\Aut(\eta) \to \Aut(\eta|_A)$$
is surjective;
\itm if $X$ is proper, and the sheaves $\cA_{\cS}$ and $\cT_{\cS}$ both
carry the structure of coherent $\cO_X$-modules, then $F_{\cS_{(X,\us)}}$
satisfies Schlessinger's (H3), so has a hull $R$;
\itm if further we have a local obstruction sheaf $\cOb$ for $\cS$,
and it carries the structure of a coherent $\cO_X$-module, then
$$h^0(X,\cT)+h^1(X,\cA)- h^0(X,\cOb) - h^1(X,\cT)- h^2(X,\cA) 
\leq \dim R -\dim \Lambda \leq \dim T_{\cS},$$ 
and if the first inequality is an equality and $\Lambda$ is regular, $R$ 
is a local complete intersection ring. If we have
$$h^0(X,\cOb) = h^1(X,\cT) = 
\dim T_{\cS} + h^2(X,\cA)-h^1(X,\cA)-h^0(X,\cT) = 0,$$ 
then $R$ is smooth over
$\Lambda$.
\end{ilist}
\end{thm*}

\begin{proof}[Proof of Theorem \ref{main-rep}] Putting together
Proposition \ref{prop:def-functor}, Proposition \ref{h4}, Theorem 
\ref{thm:schl} and Theorem \ref{main-desc}, we obtain everything except 
the last part of the desired statement. Noting that from the first
exact sequence of Theorem \ref{main-desc}, the dimension of the third
obstruction space $H^2(X, \cA)/H^0(X,\cT)$ is given by
$\dim T_{\cS} + h^2(X,\cA)-h^1(X,\cA)-h^0(X,\cT)$, the only ingredient 
still missing is the below theorem, entirely in the realm of classical 
deformation theory.
\end{proof}

\begin{thm}\label{dim-ests} Let $F$ be a deformation functor satisfying (H3), 
so that it has a hull $(R,\xi)$. Suppose that $F$ has a successive 
obstruction theory taking values in finite-dimensional spaces 
$V_1,\dots,V_m$. Then we have:
$$\dim T_F - \sum_i \dim V_i \leq \dim R-\dim \Lambda \leq \dim T_F,$$
and if the first inequality is an equality and $\Lambda$ is regular, 
then $R$ is a local complete intersection ring. If further $\dim V_i=0$ 
for all $i$, then $R$ is smooth over $\Lambda$.
\end{thm}

The theorem is essentially due to Mori, and largely follows the argument 
presented in Proposition 2.A.11 of \cite{h-l}. We include the argument here 
partly for convenience, and partly because our statement is more general, 
and requires some slight modifications.

We use a simple lemma to reduce to the prorepresentable case:

\begin{lem} Suppose that $F_1,F_2$ are functors $\Art(\Lambda,k)\to \Set$,
and we have a morphism $f:F_1 \to F_2$ which is formally smooth, and 
a successive obstruction theory for $F_2$ taking values in $V_1,\dots,V_m$.
Then $f$ induces a successive obstruction theory for $F_1$ taking 
values in $V_1,\dots,V_m$.
\end{lem}

\begin{proof} Given $\eta \in F_1(A)$, and a small extension $A' \to A$,
we can define the obstruction $\ob_{\eta,A'}$ to be simply
$\ob_{f(\eta),A'}$. The smoothness of $f$ then implies that $\eta$ can
be lifted to $A'$ if and only if $f(\eta)$ can be lifted to $A'$, so
the main conditions for an obstruction theory are satisfied, and it
remains only to check functoriality, which follows from functoriality of
the obstruction theory given for $F_2$ together with the required 
functoriality of $f$.
\end{proof}

\begin{proof}[Proof of the theorem] By the definition of a hull 
$h_R|_{\Art(\Lambda,k)} \overset{\xi}{\to} F$ is formally smooth, and 
gives an isomorphism 
of tangent spaces. By the lemma, we have a successive obstruction theory 
taking values in the same $V_i$ for $h_R|_{\Art(\Lambda,k)}$ as well, so 
it is enough to 
prove the theorem in the case that $F$ is prorepresentable. 

In this case, we work explicitly: if we write $d= \dim T_F$, Schlessinger's 
construction of $R$ in the proof of 2.11 of \cite{sc2}
is as a quotient of $B:=\Lambda[[t_1,\dots,t_d]]$ by some
ideal $J$, so to prove the theorem, it is enough to show that the number of 
generators of $J$ is bounded above by $\sum_i \dim V_i$. 

By the Artin-Rees lemma, 
$J \cap \fm_{B}^n \subseteq J \fm_{B}$ 
for some $n$. We now set 
$A^{(0)}=\Lambda[[t_1,\dots,t_d]]/(\fm_{B} J + \fm_{B}^n)$, 
and $A=R/\fm_R^n=\Lambda[[t_1,\dots,t_d]]/(J+\fm_{B}^n)$, so that 
we get a small extension
$$0 \to I \to A^{(0)} \to A \to 0$$
with $I=(J+\fm_{B}^n)/(\fm_{B} J + \fm_{B}^n) =
J/\fm_{B} J$. From the natural map $R \to A$ we obtain an object
$\xi_A \in F(A)$, with an obstruction 
$\ob_{\xi_A,A^{(0)}} \in V_i \otimes_k I$ to 
extending $\xi_A$ to $A^{(0)}$, for some $i$. If $i=1$, we can write
$\ob_{\xi_A,A^{(0)}}=\sum_{j=1}^{\dim V_1} v_{1,j} \otimes \bar{x}_{1,j}$, 
where the
$v_{1,j}$ form a basis for $V_1$, and the $\bar{x}_j$ are the images in $I$
of elements $x_{1,j} \in J$. If $i>0$, we declare the $x_{1,j}$ all to be $0$. 
We then consider the ring $A^{(1)}:=A^{(0)}/(x_{1,1},\dots,x_{1,\dim V_1})$; 
this surjects onto $A$ with kernel $I^{(1)}$, and
we again have an obstruction $\ob_{\xi_A,A^{(1)}}$ to extending $\xi_A$ to
$A^{(1)}$. 

If $m=1$, we stop. Otherwise, by the functoriality of the 
obstruction, we see that $\ob_{\xi_A,A^{(1)}} \in V_{i'} \otimes_k I^{(1)}$ 
for $i'>1$: indeed, we could only have $i'=1$ if we had 
before $i=1$, in which case the functoriality implies, since we modded out
by the $x_{1,j}$, that $\ob_{\xi_A,A^{(1)}}=0$, which is only allowed if 
$m=1$. We thus can write
$\ob_{\xi_A,A^{(1)}}=\sum_{j=1}^{\dim V_2} v_{2,j} \otimes \bar{x}_{2,j}$
as before (again, setting all $x_{2,j}=0$ if $i'>2$), 
and set $A^{(2)}=A^{(1)}/(x_{2,j},\dots,x_{2,\dim V_2})$. We repeat this
process until we have constructed $A'':=A^{(m)}$, which we see immediately
is obtained from $A'$ by modding out by (at most) $\sum_i V_i$ elements. 

We note that again by the functoriality of obstructions, we will necessarily
have $\ob_{\xi_A,A''}=0$. Thus, $\xi_A$ may be lifted to $A''$, and because
$F=h_R|_{\Art(\Lambda,k)}$, this means we can lift the map $R \to A$ to a 
map $R \to A''$.
We wish to show that this implies
\begin{equation}\label{lift-inc} J \subseteq 
\fm_{B} J + (\{x_{i,j}\}_{i,j}) + \fm_{B}^n,\end{equation}
which is equivalent 
to the stronger assertion that we have a lifting which commutes with the
natural quotient maps from $\Lambda[[t_1,\dots,t_d]]$ to $R$ and to $A''$.
Now, if we are given any lifting, we have 
$$\xymatrix{{B=\Lambda[[t_1,\dots,t_d]]} \ar[r]\ar@{-->}[d]^{\vp} &
{R} \ar[d] \ar[dr] & {} \\
{B=\Lambda[[t_1,\dots,t_d]]} \ar[r] & {A''} \ar[r] & {A},}$$
and we can fill in the dashed arrow $\vp$ to make the diagram commute
by choosing appropriate values for $\vp(t_i)$, $i=1,\dots,d$. By hypothesis, 
$\vp$ commutes with the maps to $A$, so must be the identity modulo
$J+\fm_{B}^n$. In particular, we conclude that $\vp$ induces the
identity map on $\fm_{B}/\fm_{B}^2$, so is an isomorphism,
and then that $\vp^{-1}(J)\subseteq J+\fm_{B}^n$, so that 
$J \subseteq \vp(J) + \vp(\fm_{B}^n)=\vp(J)+\fm_{B}^n$.
But by commutativity of the maps to $R$ and $A''$, we see
$\vp(J) \subseteq \fm_{B} J + (\{x_{i,j}\}_{i,j}) + \fm_{B}^n$, 
and putting these together gives (\ref{lift-inc}).

Since we had originally 
$J \cap \fm_{B}^n \subseteq \fm_{B} J$,
we finally conclude that $J$ is contained in, hence equal to
$\fm_{B} J + (\{x_{i,j}\}_{i,j})$. By Nakayama's lemma, we
conclude that $J$ is generated by $\{x_{i,j}\}_{i,j}$, as desired.
\end{proof}

\begin{rem} We observe that Theorem \ref{dim-ests} only uses a weaker
functoriality than what we impose in the definition of an obstruction
theory: namely, that the obstruction elements for $A' \to A$ be functorial 
for restriction to any intermediate small extensions obtained from
quotients of $A'$. However, in practice it appears that the stronger
functoriality is always satisfied, and the stronger functoriality makes
our definition compatible with definitions used in other contexts, for
instance by Artin \cite{ar3}.
\end{rem}

\subsection{Special cases and examples}

Specializing to the case of locally unobstructed gd-stacks, we 
immediately find two special cases of Theorem \ref{main-desc} in
which we obtain descriptions of tangent and obstruction spaces in terms
of standard sheaf cohomology. The first is when we have no
non-trivial infinitesimal automorphisms.

\begin{cor}\label{h0h1} Let $\cS$ be a gd-stack which is locally 
unobstructed, and has trivial infinitesimal automorphisms (i.e., $\cA=0$). 
Then the presheaf given by $U \mapsto T_{\cS_{(U,\us)}}$ is already
a sheaf, and we have $T_{\cS}=H^0(X,\cT)$, with obstructions lying in 
$H^1(X,\cT)$.
\end{cor}

The second case is the situation mentioned in the introduction.

\begin{cor}\label{cor:h1h2} Let $\cS$ be a gd-stack which is locally 
unobstructed, and has locally trivial deformations (i.e., $\cT=0$).
Then we have $T_{\cS}=H^1(X,\cA)$, and obstructions lie in $H^2(X,\cA)$.
\end{cor}

We remark that there is a situation, less general than that of 
Theorem \ref{main-desc}, in which the tangent and obstruction spaces are
described as hypercohomology groups of a two-term complex. This situation
simultaneously generalizes the two cases above, and is examined in
\cite{os15}.

We next return to the examples examined earlier. We first point out that 
all examples discussed in \S \ref{subsec:first-exs} in fact have natural 
structure of gd-stacks, and that given that each example gives a 
deformation stack over every open set, checking the gd-stack conditions 
is a mere formality, because we need only work with the Zariski topology.
We will take for granted the various well-known descriptions of automorphism,
local deformation, and local obstruction sheaves. However, we will then be 
able to conclude a number of the tangent and obstruction space descriptions 
as formal consequences of Theorem \ref{main-desc}.

\begin{ex} Deformations of sheaves. Given $X_{\Lambda}$ flat over
$\Lambda$, and a coherent $\cO_X$-module $\cE$, where 
$X=X_{\Lambda}|_{\Spec k}$, let $\cGDef_{X_{\Lambda}}(\cE)$ be the 
associated gd-stack of deformations of $\cE$. The automorphism sheaf is 
then $\cHom(\cE,\cE)$, the local deformation sheaf is $\cExt^1(\cE,\cE)$,
and $\cExt^2(\cE,\cE)$ is a local obstruction sheaf. Indeed, this 
follows from Proposition 3.1.5 of Chapter IV of \cite{il1}.

In particular, if $\cE$ is locally free, the latter two sheaves are $0$,
and by Corollary \ref{cor:h1h2} we have that the tangent space is
$H^1(X,\cEnd(\cE))$, with obstructions lying in $H^2(X,\cEnd(\cE))$.
\end{ex}

\begin{ex}\label{ex:scheme-def-aut-tang-obs} Deformations of schemes. 
Let $X$ be a scheme over $k$, and $\cGDef(X)$ the associated gd-stack of 
deformations of $X$. The automorphism sheaf, local deformation sheaf,
and obstruction sheaf are described by the Lichtenbaum-Schlessinger $T^i$
sheaves on $X$; see \S 2.4 as well as 4.3.3 and 4.3.4 of \cite{l-s3}. We 
do not describe these sheaves in general, but remark that 
$T^0=\cHom(\Omega^1_{X/k},\cO_X)$, and when $X$ is a 
local complete intersection scheme, we have $T^2=0$, and if further $X$ is 
generically smooth,
we have $T^1=\cExt^1(\Omega^1_{X/k},\cO_X)$. The assertion on $T^2$ is 
3.2.2 of \cite{l-s3}, and we sketch how to see both assertions, following 
the notation of {\it loc.\ cit.} 

We suppose that $U=\Spec B$ is an affine open
subset of $X$, with $B$ realized as the quotient of a polynomial ring $R$
by an ideal $I$.
The $T^i$ on $U$ are defined 
as the cohomology of $\Hom(L^{\bullet},B)$, where $L^{\bullet}$ is a complex 
of $B$-modules constructed as follows. Let $F$ be a free $R$-module, with a
surjection $j:F \twoheadrightarrow I$, with kernel $J$. Let $J_0$ be the
submodule ({\it a priori} of $F$, but in fact of $J$) generated by all
elements of the form $j(f_1)f_2-j(f_2)f_1$ for $f_1,f_2 \in F$. The
cotangent complex is then defined as 
$$J/J_0 \to F/J_0 \otimes_R B \to \Omega^1_{R/k} \otimes_R B.$$
Note that $J/J_0$ can be given a $B$-module structure, and that
$F/J_0 \otimes _R B = F \otimes_R B = F/IF$ (see also \S 3.1 of
\cite{ha2}). 
By Proposition II.8.4A of \cite{ha1}, the map 
$I/I^2 \to \Omega^1_{R/k} \otimes_R B$ and hence our map
$F/J_0 \otimes_R B \to \Omega^1_{R/k} \otimes_R B$ has cokernel 
$\Omega^1_{B/k}$, giving the desired description of $T^0$ (of course,
one can also compute the infinitesimal automorphisms directly). 

In the case that $X$ is a local complete intersection scheme, we suppose
we have chosen $U$ small enough that it may be realized as a complete
intersection inside of $\Spec R$. Now, the claim is that in this case, 
the complex $L^{\bullet}$ in fact consists of two terms: 
$I/I^2 \to \Omega^1_{R/k} \otimes_R B$. We choose $F$ to be a free module 
generated by a minimal set of generators in $I$, which necessarily form an 
$R$-sequence. We therefore need to check that $J_0=J$, and $F/IF=I/I^2$.
The first equality follows from exactness of the Koszul complex for the
generators of $I$ (see Theorem 16.5(i) of \cite{ma1}), while the second
follows from Theorem II.8.21A(e) of \cite{ha1}.

We then have that $T^2$ is the second cohomology of a two-term complex,
and hence equal to $0$, as asserted. Furthermore, if $X$ is generically
smooth, is it also necessarily locally integral. The complex 
$I/I^2 \to \Omega^1_{R/k} \otimes_R B$ is generically injective by Theorem
II.8.17 of \cite{ha1}. But since both terms in the complex are locally
free, the map is in fact injective, and therefore gives a locally free
resolution of $\Omega^1_{B/k}$, so we conclude that $T^1$ computes
$\Ext^1(\Omega^1_{B/k}, B)$, as desired. See also Theorem 4.4 of 
\cite{vi2} for an exposition treating the global case.

Finally, if further $X$ is smooth, then also $T^1=0$, and $T^0=T_X$, the 
tangent sheaf of $X$. We thus get from Corollary \ref{cor:h1h2} that the 
tangent space is $H^1(X,T_X)$, with obstructions in $H^2(X,T_X)$.
\end{ex}

\begin{ex} Deformations of quotient sheaves. Given $\cE_{\Lambda}$
coherent on some $X_{\Lambda}$, write $\cE$ and $X$ for the
restrictions to $\Spec k$. Given also a quotient $\cF$ of $\cE$, let
$\cGDef_{X_{\Lambda},\cE_{\Lambda}}(\cF)$ be the associated gd-stack of 
deformations of $\cF$ as a quotient of $\cE$. Assume further that 
$\cE_{\Lambda}$ is 
flat over $\Lambda$. Write $\cG=\ker(\cE \twoheadrightarrow \cF)$
Then the automorphism sheaf is $0$, the local deformation sheaf is given
by $\cHom(\cG,\cF)$, and we can take $\cExt^1(\cG,\cF)$ as a
local obstruction sheaf. These follow from Lemma 2.5 and subsequent 
discussion of Olsson-Starr \cite{o-s1}, noting that under the flatness 
hypothesis on $\cE_{\Lambda}$, the first obstruction discussed there 
always vanishes.

In particular, if $X_{\Lambda}$ is flat over $\Lambda$, and $\cE$ and 
$\cF$ (and therefore also
$\cG$) are locally free, then the local obstruction sheaf is $0$, and
we have from Corollary \ref{h0h1} that the tangent space is 
$H^0(X,\cHom(\cG,\cF))$, with obstructions lying in $H^1(X,\cHom(\cG,\cF))$.
\end{ex}

\begin{ex} Deformations of subschemes. Let $X_{\Lambda}$ be a scheme flat
over $\Lambda$, and $X$ its restriction to $\Spec k$. Given $Z \subseteq X$
a closed subscheme with ideal sheaf $\cI_Z$, let $\cGDef_{X_{\Lambda}}(Z)$ 
be the associated gd-stack of deformations of $Z$ inside $X$. Then it
follows directly from the above case of quotient sheaves that the 
automorphism sheaf is $0$, the local deformation sheaf is 
$\cHom_X(\sI_Z,\cO_Z)=\cHom_Z(\sI_{Z}/\sI_{Z}^2,\cO_Z)$, and the local 
obstructions lie in $\cExt^1_X(\sI_Z,\cO_Z)$.

In particular, if $Z$ is a local complete intersection inside $X$,
in the sense that it is locally cut out by regular sequences,
then the local deformation 
sheaf is the normal bundle $\cN_{Z/X}$, and we claim that there are no 
local obstructions. Indeed, this may be shown directly by observing that
any liftings of local equations for $Z$ inside $X$ will yield a local 
deformation of $Z$;
see for instance \cite{vi2}, Lemma 2.7 and the preceding discussion. 
Thus, from Corollary 
\ref{h0h1} we see that the tangent space is $H^0(Z,\cN_{Z/X})$ and the 
obstructions lie in $H^1(Z,\cN_{Z/X})$.
\end{ex}

\begin{ex} Deformations of morphisms. Given $X_{\Lambda},Y_{\Lambda}$
locally of finite type over $\Lambda$, with $X_{\Lambda}$ flat and
$Y_{\Lambda}$ separated, and $X$ and $Y$ the respective restrictions to
$\Spec k$, suppose $f:X \to Y$ is a morphism, and let 
$\cGDef_{X_{\Lambda},Y_{\Lambda}}(f)$ be the gd-stack of deformations of 
$f$. Then from the previous example we see that the automorphism sheaf is 
trivial, the local deformation sheaf is $\cHom(f^* \Omega^1_{Y/k},\cO_X)$, 
and local obstructions lie in
$\cExt^1_{X \times_k Y} (\sI_{\Gamma(f)},\cO_{\Gamma(f)})$.

In particular, if $Y$ is smooth, then $\Gamma(f)$ is a local complete
intersection, so we have by the previous example and Corollary \ref{h0h1} 
that the tangent space is $H^0(X,f^* T_Y)$, with obstructions in 
$H^1(X,f^* T_Y)$.
\end{ex}

\begin{ex} Deformations of connections. Suppose we have a scheme 
$X_{\Lambda}$ smooth over $\Lambda$, and a locally free 
$\cO_{X_{\Lambda}}$-module $\cE_{\Lambda}$, and write $X$ and $\cE$ for
the restrictions to $\Spec k$. Given also a connection $\nabla$ on $\cE$, 
let $\cGDef_{X_{\Lambda},\cE_{\Lambda}}(\nabla)$
be the gd-stack of deformations of $\nabla$. Then the automorphism sheaf
is $0$ by definition, the local deformation sheaf is 
$\cHom(\cE,\cE \otimes \Omega^1_{X/k})$, and local obstructions vanish. 
Indeed, both the last two statements follow from the fact that $\cE$ is
locally free, so that connections may be expressed locally explicitly in
terms of matrices with coefficients in $\Omega^1_{X/k}$.
By Corollary \ref{h0h1}, we find that the tangent space is 
$H^0(X,\cEnd(\cE)\otimes \Omega^1_{X/k})$, with obstructions lying in 
$H^1(X,\cEnd(\cE) \otimes \Omega^1_{X/k})$.
\end{ex}

We note as a consequence of Example \ref{ex:scheme-def-aut-tang-obs}
a non-trivial example of a deformation stack with trivial automorphisms 
and therefore satisfying (H4):

\begin{cor}\label{cor:rigid-schemes} Suppose that $X/k$ is a scheme with 
$\Hom(\Omega^1_{X/k},\cO_X)=0$. Then the deformation stack $\Def_X$ has 
trivial automorphisms, and therefore satisfies Schlessinger's (H4).
\end{cor}

\begin{proof} By Example \ref{ex:scheme-def-aut-tang-obs}, we see that the 
first-order infinitesimal automorphisms vanish, and the remaining 
assertions follow from Corollary \ref{cor:triv-aut-def}.
\end{proof}

\begin{ex} Suppose $X$ is a smooth, proper curve of genus at least $2$.
Then $\Def_X$ has trivial automorphisms, and satisfies Schlessinger's 
(H4).
\end{ex}

\appendix
\section{Two lemmas of Schlessinger}

In order to be as self-contained as possible, we include here the statements
of two lemmas of Schlessinger, which play a key role in checking both
Schlessinger's criteria and the deformation stack condition in several 
important examples. Although the proofs are not difficult, we do not
reproduce them here.

\begin{alem}\label{lem:schl1} (Schlessinger, Lemma 3.3 of \cite{sc2})
Let $A$ be a ring, with nilpotent ideal $J$, and $u:M \to N$ a homomorphism
of $A$-modules, with $N$ flat over $A$. If $\bar{u}:M/JM \to N/JN$ is
an isomorphism, then $u$ is an isomorphism.
\end{alem}

It follows easily that:

\begin{acor}\label{cor:mod-free} Let $M$ be a flat module over an Artin 
local ring $A$. Then $M$ is free.
\end{acor}

\begin{alem}\label{lem:schl2} (Schlessinger, Lemma 3.4 of \cite{sc2})
Consider a commutative diagram
$$\xymatrix{{N} \ar[rrr]^{p''}\ar@{-}[ddd]\ar[dr]^{p'} & {} & {} & {M''}
\ar@{-}[ddd]\ar[dr]^{u''} & {} \\
{} & {M'} \ar[rrr]^{u'}\ar@{-}[ddd] & {} & {} &{M}\ar@{-}[ddd] \\
{} & {} & {} & {} & {} \\
{B} \ar[rrr]\ar[dr] & {} & {} & {A''} \ar[dr] & {} \\
{} & {A'} \ar[rrr] & {} & {} & {A}}$$
of compatible ring and module homomorphisms, where $B=A' \times_A A''$,
$N=M' \times _M M''$, and $M'$ and $M''$ are flat over $A'$ and $A''$,
respectively. Suppose
\begin{ilist}
\itm $A'' \twoheadrightarrow A$, with nilpotent kernel,
\itm $u'$ induces $M' \otimes_{A'} A \risom M$, and similarly for $u''$.
\end{ilist}

Then  $N$ is flat over $B$, and $p'$ induces $N \otimes_B A' \risom M'$,
and similarly for $p''$.
\end{alem}

\section{Hypercovers and $H^2$}\label{app:hyper-cech}

It is very widely known that although \v{C}ech cohomology provides an
effective tool for making sheaf cohomology more concrete, it does not
always agree with (the derived functor version of) sheaf cohomology.
However, less widely known and even less widely used is the version
of \v{C}ech cohomology developed by Verdier in \S 7 of Expose V of
\cite{sga42}, which always agrees with sheaf cohomology (even on an
arbitrary site). The basic idea is quite simple: instead of fixing
an open cover $\{U_i\}$ of a space $X$ and working only on the
intersections of the various $U_i$, one allows further refinements at each
stage, taking covers of each $U_i \cap U_j$, and so forth. 
In order to avoid unnecessary hypotheses in our main theorem,
and because of a lack of suitably down-to-earth references, in this 
appendix we describe how to use the simplest non-trivial aspect of 
hypercovers to describe the sheaf cohomology group $H^2$ on an arbitrary 
topological space. 

We begin by defining the type of covers we work with, using terminology
consistent with that of Beke \cite{be2}:

\begin{adefn} Given a topological space $X$, a {\bf cover $\cU$ of level $2$} 
of $X$ consists of an open cover $\{U_i\}_{i\in I}$ of $X$, together 
with the data of a cover $\{U_{i_0,i_1,j}\}_{j \in J_{i_0,i_1}}$ of 
each $U_{i_0} \cap U_{i_1}$, with $i_0 \neq i_1$.
\end{adefn}

Note that for the sake of simplicity, the data of a cover of level $2$ does 
not include different covers for $U_{i_0} \cap U_{i_1}$ and for 
$U_{i_1} \cap U_{i_0}$; by the same token, we take the cover of
$U_{i_0} \cap U_{i_0}$ to be the single open set $U_{i_0}$, and denote
this by an underscore. For instance, $U_{i_0,i_0,\_}=U_{i_0}$ is the 
unique open set of the cover of $U_{i_0} \cap U_{i_0}$.

Associated to $\cU$ and a sheaf $\cF$, we will define $H^2(\cU,\cF)$, a 
cohomology group (in fact, a module over $\Gamma(X,\cO_X)$) which directly 
generalizes \v{C}ech cohomology, and our main purpose is to prove the 
following very special case of Verdier's work:

\begin{athm} Given any ringed space $(X,\cO_X)$ and an $\cO_X$-module
$\cF$, and any $\cU$ a cover of level $2$ of $X$, there exists a 
natural homomorphism of $\Gamma(X,\cO_X)$-modules
$H^2(\cU, \cF) \to H^2(X, \cF)$.

Furthermore, this map is compatible with refinement, and induces an
isomorphism
$$\limr_{\cU} H^2(\cU,\cF) \risom H^2(X,\cF).$$
\end{athm}

Here, a refinement of a cover of level $2$ is defined in the natural
way. Formally, it is:

\begin{adefn} Let $\cU$ be a cover of level $2$ of $X$. A 
{\bf refinement} $\cU'$ of $\cU$ is another cover of level $2$ of $X$,
together with order-preserving maps of index sets $\pi: I' \to I$ and for 
each $i'_0<i'_1 \in I'$ maps of index sets 
$\pi_{i'_0,i'_1}:J_{i'_0,i'_1} \to J_{\pi(i'_0),\pi(i'_1)}$, such
that $U'_{i'} \subseteq U_{\pi(i')}$ for each $i' \in I'$, and
$U_{i'_0,i'_1,j'} \subseteq U_{\pi(i'_0),\pi(i'_1),\pi_{i'_0,i'_1}(j')}$
for each $j' \in J'_{i'_0,i'_1}$.
\end{adefn}

Note that covers of level $2$ form a directed set under refinement, so it 
makes sense to take direct limits over all covers of level $2$ of a given 
space. 

As an immediate corollary of the theorem, we obtain the same statements for 
cohomology of sheaves of abelian groups on any topological space. We give 
the statement for ringed spaces because the proof is the same, and we will
want to know that our $k$-linear structure is preserved by the 
isomorphism.

The definition of the group $H^2(\cU,\cF)$ is the following. Note
that unlike \cite{ha1}, we make use of unordered Cech cohomology
rather than alternating Cech cohomology. 

\begin{adefn} Given an $\cO_X$-module $\cF$ on $X$, and $\cU$ a 
cover of level $2$ of $X$, we make the following definitions:
$C^1(\cU,\cF)$ is the $\Gamma(X,\cO_X)$-module of $1$-cochains
$$\prod_{i_0,i_1 \in I, j \in J_{i_0,i_1}} \cF(U_{i_0,i_1,j}),$$
and 
$C^2(\cU,\cF)$ is the $\Gamma(X,\cO_X)$-module of $2$-cochains
$$\prod_{\scriptsize\begin{matrix}\bi=(i_0,i_1,i_2)\in I^3 \\
\bj = (j_{1,2}, j_{0,2}, j_{0,1}) \in 
J_{i_1,i_2} \times J_{i_0,i_2} \times J_{i_0,i_1}\end{matrix}} 
\!\!\!\!\!\!\!\!\!\!\!\!\!\!\!\!\!\!\!\!\!\!\!\!\!\!\!\cF(U_{\bi,\bj}),$$
where
$U_{\bi, \bj} := U_{i_1,i_2,j_{1,2}} \cap U_{i_0,i_2,j_{0,2}} 
\cap U_{i_0,i_1,j_{0,1}}$. 

Then $Z^2(\cU,\cF)$ is the submodule consisting of cocycles, which is
to say $\rho \in C^2(\cU,\cF)$ satisfying the condition that for any 
$\bi'=(i_0,i_1,i_2,i_3) \in I$ and 
$\bj':=\{j_{m,m'}\in J_{i_m,i_{m'}}\}_{m,m' \in \{0,1,2,3\}}$ we have
\begin{multline*}0=(d\rho)_{\bi',\bj'}:= 
\rho_{i_1,i_2,i_3,j_{2,3},j_{1,3},j_{1,2}}|_{U_{\bi',\bj'}}
-\rho_{i_0,i_2,i_3,j_{2,3},j_{0,3},j_{0,2}}|_{U_{\bi',\bj'}} \\
+\rho_{i_0,i_1,i_3,j_{1,3},j_{0,3},j_{0,1}}|_{U_{\bi',\bj'}}
-\rho_{i_0,i_1,i_2,j_{1,2},j_{0,2},j_{0,1}}|_{U_{\bi',\bj'}}.
\end{multline*}

Next, we define $B^2(\cU,\cF)$ as the submodule of $Z^2(\cU,\cF)$ 
consisting of coboundaries. These are the elements obtainable
from a $1$-cochain $\rho' \in C^1(\cU,\cF)$ by setting 
$\rho_{\bi, \bj}=(d \rho')_{\bi,\bj}:=
\rho'_{i_1,i_2,j_{1,2}}|_{U_{\bi, \bj}}
-\rho'_{i_0,i_2,j_{0,2}}|_{U_{\bi, \bj}}+
\rho'_{i_0,i_1,j_{0,1}}|_{U_{\bi,\bj}}$.

Finally, $H^2(\cU,\cF)$ is defined to be $Z^2(\cU,\cF)/B^2(\cU,\cF)$.
\end{adefn}

Note that in the special case that each cover $\{U_{i_0,i_1,j}\}_j$
consists of the single set $U_{i_0} \cap U_{i_1}$, we recover precisely
the usual definition of the \v{C}ech cohomology group
$\check{H}^2(\{U_i\}, \cF)$.

We take a very hands-on approach to proving the theorem.
A key lemma is the following:

\begin{alem} Suppose that $\cF$ is flasque. Then $H^2(\cU,\cF)=0$.
\end{alem}

\begin{proof} Suppose we have $\rho \in Z^2(\cU,\cF)$. 
We will construct a certain coboundary $d\rho' \in B^2(\cU,\cF)$, 
and show that subtracting it from $\rho$ removes the dependence on
refinements, so that we obtain an element of the standard \v{C}ech 
cohomology group $\check{H}^2(\{U_i\},\cF)$. Since $\cF$ is
flasque, this element is then a coboundary by standard \v{C}ech 
cohomology theory. 

$\rho'$ is constructed as follows. Fix $i_0, i_1$, and for each
$j,j' \in J_{i_0,i_1}$, consider 
$$\rho_{i_0,i_0,i_1,j,j',\_}
-\rho_{i_0,i_0,i_0,\_,\_,\_}\in \cF(U_{i_0,i_1,j}\cap U_{i_0,i_1,j'}).$$
This yields a \v{C}ech $1$-cochain for the cover 
$\cU_{i_0,i_1}:=\{U_{i_0,i_1,j}\}_{j}$
of $U_{i_0} \cap U_{i_1}$, and in fact we obtain a $1$-cocycle, as can
be verified using
the cocycle condition for $\rho$, evaluated at $\bi=(i_0,i_0,i_0,i_1)$. 
Since $\cF$ is
flasque, the $1$-cocycle can be written as $d\rho^{i_0,i_1}$ for
some $0$-cochain $\rho^{i_0,i_1}$ on $\cU_{i_0,i_1}$. Having done
this for all $i_0, i_1$, we then
define $\rho'$ by $\rho'_{i_0,i_1,j}:=\rho^{i_0,i_1}_{j}$. Then $d\rho'$ 
is our desired coboundary.

We next claim that $\rho''=\rho-d\rho'$ is a standard \v{C}ech 
$2$-cocycle, as desired. It suffices to show that for any
$\bi=(i_0,i_1,i_2)$ and $\bj,\bj'$, we have
$\rho''_{\bi,\bj}|_{U_{\bi,\bj}\cap U_{\bi,\bj'}}
- \rho''_{\bi,\bj'}|_{U_{\bi,\bj}\cap U_{\bi,\bj'}} = 0$.
Expanding the expression we find that the left hand side is given by
$$\rho_{\bi,\bj}-\rho_{\bi,\bj'}+\rho_{i_0,i_0,i_2,j_{0,2},j_{0,2}',\_}
-\rho_{i_0,i_0,i_1,j_{0,1},j_{0,1}',\_}
-\rho_{i_1,i_1,i_2,j_{1,2},j_{1,2}',\_}
+\rho_{i_1,i_1,i_1,\_,\_,\_}.$$
One then checks that the desired statement follows from the cocycle 
condition on $\rho$, evaluated at $(i_0,i_0,i_1,i_2)$, 
$(i_0,i_1,i_1,i_2)$, and $(i_0,i_1,i_1,i_1)$.
\end{proof}

It is now relatively straightforward to prove the theorem.

\begin{proof}[Proof of theorem] We begin by constructing the asserted
natural map. Let 
$$\cF \to \cG_0 \overset{d_0}{\to} \cG_1 \overset{d_1}{\to} \cG_2 
\overset{d_2}{\to} \cG_3 \to \dots$$
be a flasque resolution of $\cF$. We compute everything in terms of
this resolution, making the identification $H^2(X,\cF)=\ker d_2/\im d_1$. 
Given $\rho \in H^2(\cU,\cF)$, and choosing a representative $2$-cocycle
for $\rho$, we obtain its image $\rho^0 \in Z^2(\cU,\cG_0)$, which by the 
lemma is $d \tilde{\rho}{^0}$ for some $\tilde{\rho}{^0}\in C^1(\cU,\cG_0)$. 
Taking the image of $\tilde{\rho}{^0}$ in
$\cG_1$, we obtain $\rho^1\in C^1(\cU,\cG_1)$. We claim that in fact
$\rho^1$ is a standard \v{C}ech $1$-cocycle. We have 
$d\rho^1=d(d_0(\tilde{\rho}{^0}))=d_0(d(\tilde{\rho}{^0}))=d_0(\rho^0)=0$,
since $\rho^0$ consists of sections of $\cF$.
Thus, $\rho^1$ satisfies the cocycle
condition, and is a standard \v{C}ech $1$-cocycle if and only if 
for any $i_0 \neq i_1$, and $j,j' \in J_{i_0,i_1}$, we have 
$$\rho^1_{i_0,i_1,j}|_{U_{i_0,i_1,j}\cap U_{i_0,i_1,j'}}
-\rho^1_{i_0,i_1,j'}|_{U_{i_0,i_1,j}\cap U_{i_0,i_1,j'}}=0.$$
However, setting $\bi=(i_0,i_0,i_1)$ and $\bj=(j,j',\_)$, and evaluating
$$0=(d\rho^1)_{i_0,i_0,i_1,j,j',\_}=
\rho^1_{i_0,i_1,j}|_{U_{\bi,\bj}}-\rho^1_{i_0,i_1,j'}|_{U_{\bi,\bj}}
+\rho^1_{i_0,i_0,\_}|_{U_{\bi,\bj}},$$
we need only see that $\rho^1_{i_0,i_0,\_}=0$, which follows from
evaluating $0=(d\rho^1)_{i_0,i_0,i_0,\_,\_,\_}=\rho^1_{i_0,i_0,\_}$.
We therefore have that $\rho^1$ is a \v{C}ech $1$-cocycle, and again
by standard \v{C}ech cohomology we conclude that 
$\rho^1=d\tilde{\rho}{^1}$ for some \v{C}ech $0$-cochain 
$\tilde{\rho}{^1}$ on $\{U_i\}$ with coefficients in $\cG_1$. 
Finally taking the image of $\tilde{\rho}{^1}$ in $\cG_2$, we see that we 
get a global section $\rho^2$ of $\cG_2$, since the differences on 
$U_{i_0} \cap U_{i_1}$ take values in $\cG_0$. Moreover, this global 
section lies in $\ker d_2$, since it is locally in the image of $d_1$. 
Taking the class of $\rho^2$ in $H^2(X,\cF)$ constructs the desired map 
$H^2(\cU, \cF) \to H^2(X, \cF)$.

It remains to check that this map is a well-defined homomorphism, compatible
with refinement, and bijective in the limit. For the map to be
well defined, we need to note that if $\rho$ is modified by a 
$2$-coboundary of $\cF$, we will have $\tilde{\rho}^0$ modified by a 
$1$-cochain in the image of $\cF$, which then vanishes after taking
the image in $\cG_1$, and similarly if $\tilde{\rho}{^0}$ is modified by a 
\v{C}ech $1$-cocycle (equivalently, coboundary) of $\cG_{0}$, we will 
modify $\tilde{\rho}{^1}$ by a $0$-cochain in the image of $\cG_0$, which 
vanishes in $\cG_2$. Finally, if $\tilde{\rho}{^1}$ is modified by a
$0$-cocycle of $\cG_1$, then its image in $\ker d_2$
is modified by an element of $\im d_1$, leaving the class in 
$H^2(X,\cF)$ unaffected.
We also see easily that the map we have defined
is a module homomorphism, since each of the $d_i$ are homomorphisms.
Compatibility with refinement from $\cU$ to $\cU'$ is clear, since at each 
stage we can choose $\rho$, $\rho^0$, and $\rho^1$ on $\cU'$ to be 
obtained by refinement from $\cU$, and the image $\rho^2$ in $\cG_2$
will then be unchanged. 

Injectivity and surjectivity in the limit then follow by explicit
construction: if we have for some cover $\cU$ of level $2$ an element
$\rho \in H^2(\cU,\cF)$ having image $0$ in $H^2(X,\cF)$, we have that the 
image of $\tilde{\rho}^1$ in $\cG_2$ agrees with the image of some global 
section $s_{\rho}\in \Gamma(\cG_1)$, so that 
$\tilde{\rho}^1_i-s_{\rho}|_{U_i}$ maps to $0$ in $\cG_2$, and after 
possible refinement of $\cU$, is the image of some $0$-cochain of $\cG_0$.
Thus, $\tilde{\rho}^1$ is the sum of a $0$-cocycle of $\cG_1$ and 
$0$-cochain of $\cG_0$, which means that $\rho^1$ is the image of a
$1$-coboundary of $\cG_0$, and hence that $\tilde{\rho}^0$ is the
sum of a $1$-coboundary of $\cG_0$ and a $1$-cochain in the kernel of
$d_0$. The latter cochain is in the image of $\cF$ after refining $\cU$ 
once more, so we find that $\rho^0$ is the image of a $2$-coboundary of
$\cF$, so that we have a refinement $\cU'$ of $\cU$ on which $\rho$ 
vanishes, as desired.

Finally, surjectivity in the limit is proved by
starting with $\xi \in H^2(X,\cF)$, represented as an element of 
$\ker d_2 \subseteq \Gamma(\cG_2)$, and choosing $\{U_i\}$ to be a cover
of $X$ on which $\xi$ is the image of some $\xi_i \in \cG_1(U_i)$. 
Then for each $i_0 \neq i_1$, we let $\{U_{i_0,i_1,j}\}_{j}$ be a cover 
of $U_{i_0}\cap U_{i_1}$ on which $\xi_{i_1}-\xi_{i_0}$ is the image of 
some $\tilde{\rho}{^0}_{i_0,i_1,j} \in \cG_0(U_{i_0,i_1,j})$.
This gives us a $\tilde{\rho}{^0} \in C^1(\cU,\cG_0)$,
and taking $\rho^0=d \tilde{\rho}{^0}$ 
we see that every section $\rho^0_{\bi,\bj}$ is in the kernel of $d_0$,
so after possible further refinement, is the image of some 
$\rho_{\bi,\bj} \in \cF(U_{\bi,\bj})$. We have thus constructed 
$\rho \in Z^2(\cU,\cF)$, with image $\xi \in H^2(X,\cF)$, as desired.
\end{proof}

We conclude with an example demonstrating that all of this was actually
necessary, even for very reasonable (by algebrogeometric standards!)
topological spaces. See also the discussion preceding the proof of
Proposition \ref{prop:very-short}.

\begin{aex}\label{ex:hyper-cech} We let $X$ be the topological space 
underlying the scheme
$\AA^2_k$, where $k$ is any field. We produce a presheaf $\cF$ of 
abelian groups such that the sheafification $\tilde{\cF}$ is the constant
sheaf associated to $\ZZ$, but there is no open cover $\{U_i\}$ of
$X$ with sections $\rho_i \in \cF(U_i)$ such that each $\rho_i$ maps to
$1$ in $\tilde{F}(U_i)$ and $\rho_i|_{U_i \cap U_j}=\rho_j|_{U_i \cap U_j}$
in $\cF(U_i \cap U_j)$.

We construct $\cF$ as follows. Let $V_1=\AA^2_k\smallsetminus\{(0,0)\}$,
and $V_2=\AA^2_k\smallsetminus \{(1,0)\}$. Let $\cV$ be the collection
of open subsets of $\AA^2_k$ whose complement contains an irreducible
curve through $(0,0)$ and $(1,0)$. We then define:
$$\cF(U)=\begin{cases} 
0: & U \not\subseteq V_1 \text{ and } U \not\subseteq V_2; \\
\ZZ \times 0 : & U \subseteq V_1 \text{ and } U \not\subseteq V_2; \\
0 \times \ZZ : & U \subseteq V_2 \text{ and } U \not\subseteq V_1; \\
\ZZ \times \ZZ : & U \subseteq V_1 \cap V_2 \text{ and } U \not\in \cV; \\
\ZZ : & U \in \cV.
\end{cases}$$
Restriction maps $\ZZ \times \ZZ \to \ZZ$ are given by the summation map,
while the remaining restriction maps are the obvious ones. 

We note that
$\cV$ is an open cover of $V_1 \cap V_2$, since we can for instance 
remove the curves $y=0$ or $y=x(x-1)$, whose intersection is precisely
$(0,0)$ and $(1,0)$. It is then easy to check that $\cF$ is a presheaf
with constant sheafification $\tilde{\cF}=\underline{\ZZ}$, and that if 
$\{U_i\}$ is any cover of $X$ with $\rho_i \in \cF(U_i)$ mapping to $1$ in 
$\tilde{\cF}$ for all $i$, if we choose $U_i$ containing $(0,0)$ and $U_j$
containing $(1,0)$, then we must have 
$\rho_i|_{U_i \cap U_j} \neq \rho_j|_{U_i \cap U_j}$. Specifically, we
check that we must have $U_i \subseteq V_2$ but $U_i \not \subseteq V_1$
and $U_j \subseteq V_1$ but $U_j \not\subseteq V_2$, and that
$U_i \cap U_j \not\in \cV$, from which the claim clearly follows.
\end{aex}

\bibliographystyle{hamsplain}
\bibliography{hgen}
\end{document}

%% file: headerstemplate.tex
\usepackage{amssymb,euscript,amsmath, mathrsfs}
\usepackage[dvips]{graphicx}
\usepackage[dvips]{color}

\newcounter{ENUM}
\newcommand{\itm}{\item}
\newenvironment{ilist}{\renewcommand{\theENUM}{\roman{ENUM}}\renewcommand{\itm}{\addtocounter{ENUM}{1}\item[(\theENUM)]}\begin{itemize}\setcounter{ENUM}{0}}{\end{itemize}}

\def\risom{\overset{\sim}{\rightarrow}}

\def\limr{\varinjlim}

\input xy
\xyoption{all}
\CompileMatrices

\def\real{{\;\preccurlyeq\;}}
\def\us{{\,\underline{\;\;\:}\,}}
\def\bi{{\mathbf i}}
\def\bj{{\mathbf j}}

\def\ZZ{{\mathbb Z}}

\def\AA{{\mathbb A}}

\def\cF{{\mathcal F}}
\def\cS{{\mathcal S}}
\def\cG{{\mathcal G}}

\def\cE{{\mathcal E}}
\def\cO{{\mathcal O}}
\def\cI{{\mathcal I}}

\def\cA{{\mathcal A}}
\def\cN{{\mathcal N}}

\def\cT{{\mathcal T}}
\def\cU{{\mathcal U}}
\def\cV{{\mathcal V}}

\def\cEnd{{\mathcal End}}
\def\cHom{{\mathcal Hom}}
\def\cExt{{\mathcal Ext}}

\def\cOb{{\mathcal O}b}
\def\cDef{{\mathcal D}ef}
\def\cGDef{{\mathcal G}{\mathcal D}ef}
\def\sI{{\mathscr I}}

\def\fm{{\mathfrak m}}

\def\vp{\varphi}

\def\Ext{\operatorname{Ext}}
\def\Mor{\operatorname{Mor}}
\def\Hom{\operatorname{Hom}}

\def\Aut{\operatorname{Aut}}

\def\Spec{\operatorname{Spec}}

\def\id{\operatorname{id}}

\def\Def{\operatorname{Def}}

\def\im{\operatorname{im}}

\def\Art{\operatorname{Art}}
\def\Set{\operatorname{Set}}
\def\ob{\operatorname{ob}}
\def\Zar{\operatorname{Zar}}
\def\Obj{\operatorname{Obj}}
\def\opp{\operatorname{opp}}
\def\sw{\operatorname{sw}}

\newtheorem{thm}{Theorem}[subsection]
\newtheorem*{thm*}{Theorem}
\newtheorem{prop}[thm]{Proposition}
\newtheorem{lem}[thm]{Lemma}
\newtheorem{cor}[thm]{Corollary}

\newtheorem{alem}{Lemma}[section]
\newtheorem{acor}[alem]{Corollary}
\newtheorem{athm}[alem]{Theorem}
\theoremstyle{definition}
\newtheorem{defn}[thm]{Definition}

\newtheorem{ex}[thm]{Example}

\newtheorem{adefn}[alem]{Definition}
\newtheorem{aex}[alem]{Example}
\theoremstyle{remark}
\newtheorem{notn}[thm]{Notation}
\newtheorem{rem}[thm]{Remark}

\numberwithin{equation}{subsection}
\numberwithin{figure}{subsection}

%% file: deform-1.bbl
\newcommand{\noopsort}[1]{} \newcommand{\printfirst}[2]{#1}
  \newcommand{\singleletter}[1]{#1} \newcommand{\switchargs}[2]{#2#1}
\providecommand{\bysame}{\leavevmode\hbox to3em{\hrulefill}\thinspace}
\begin{thebibliography}{10}

\bibitem{ar3}
Michael Artin, \emph{Versal deformations and algebraic stacks}, Inventiones
  Mathematicae \textbf{27} (1974), 165--189.

\bibitem{sga42}
Michael Artin, Alexandre Grothendieck, and Jean-Louis Verdier, \emph{Th\'eorie
  des topos et cohomologie \'etale des sch\'emas. {T}ome 2}, SGA, no.~4,
  Spring-Verlag, 1972.

\bibitem{sga43}
\bysame, \emph{Th\'eorie des topos et cohomologie \'etale des sch\'emas. {T}ome
  3}, SGA, no.~4, Spring-Verlag, 1973.

\bibitem{be2}
Tibor Beke, \emph{Higher \v {C}ech theory}, $K$-theory \textbf{32} (2004),
  no.~4, 293--322.

\bibitem{f-g1}
Barbara Fantechi and Lothar G\"ottsche, \emph{Local properties and {H}ilbert
  schemes of points}, Fundamental Algebraic Geometry, Mathematical Surveys and
  Monographs, vol. 123, American Mathematical Society, 2005, pp.~139--178.

\bibitem{ega44}
Alexander Grothendieck and Jean Dieudonn\'e, \emph{{\'E}l\'ements de
  g\'eom\'etrie alg\'ebrique: {IV.} \'{E}tude locale des sch\'emas et des
  morphismes de sch\'emas, quatri\'eme partie}, vol.~32, Publications
  math\'ematiques de l'I.H.\'E.S., no.~2, Institut des Hautes \'Etudes
  Scientifiques, 1967.

\bibitem{gr4}
Alexandre Grothendieck, \emph{Cat\'egories cofibr\'ees additives et complexe
  cotangent relatif}, Lecture Notes in Mathematics, no.~79, Spring-Verlag,
  1968.

\bibitem{ha2}
Robin Hartshorne, \emph{Lectures on deformation theory}, in preparation.

\bibitem{ha1}
\bysame, \emph{Algebraic geometry}, Springer-Verlag, 1977.

\bibitem{h-l}
Daniel Huybrechts and Manfred Lehn, \emph{The geometry of moduli spaces of
  sheaves}, Max-Planck-Institut fur Mathematik, 1997.

\bibitem{il1}
Luc Illusie, \emph{Complexe contangent et d\'eformations. {I}}, Lecture Notes
  in Mathematics, no. 239, Springer-Verlag, 1971.

\bibitem{l-s3}
Stephen Lichtenbaum and Michael Schlessinger, \emph{The cotangent complex of a
  morphism}, Transactions of the AMS \textbf{128} (1967), 41--70.

\bibitem{ma1}
Hideyuki Matsumura, \emph{Commutative ring theory}, Cambridge University Press,
  1986.

\bibitem{o-s1}
Martin Olsson and Jason Starr, \emph{Quot functors for {D}eligne-{M}umford
  stacks}, Communications in Algebra \textbf{31} (2003), no.~8, 4069--4096,
  \mbox{arxiv:math.AG/0204307}.

\bibitem{ol1}
Martin~C.\ Olsson, \emph{Crystalline cohomology of stacks and {H}yodo-{K}ato
  cohomology}, preprint.

\bibitem{os15}
Brian Osserman, \emph{Deformations and automorphisms: triangles of deformation
  problems}, in preparation.

\bibitem{ri1}
Dock~Sang Rim, \emph{Formal deformation theory}, Groupes de monodromie en
  g\'eom\'etrie alg\'ebrique, Lecture Notes in Mathematics, vol. 288,
  Springer-Verlag, 1972, expose {VI}.

\bibitem{sc2}
Michael Schlessinger, \emph{Functors of {A}rtin rings}, Transactions of the AMS
  \textbf{130} (1968), 208--222.

\bibitem{se1}
J.~P. Serre, \emph{Local fields}, Springer-Verlag, 1979.

\bibitem{vi2}
Angelo Vistoli, \emph{The deformation theory of local complete intersections},
  \mbox{arxiv:alg-geom/9703008}.

\bibitem{vo2}
Vladimir Voevodsky, \emph{Homotopy theory of simplicial sheaves in completely
  decomposable topologies},  (2000), preprint.

\bibitem{vo3}
\bysame, \emph{Unstable motivic homotopy categories in nisnevich and
  cdh-topologies},  (2000), preprint.

\end{thebibliography}
